\documentclass[12pt]{article}
\usepackage{amsfonts,amsmath}
\usepackage{amssymb}
\usepackage{latexsym}
\usepackage{epsfig}
\usepackage{graphicx,color,graphics}

\newtheorem{lemma}{Lemma}[section]

\newtheorem{theorem}[lemma]{Theorem}

\newtheorem{corollary}[lemma]{Corollary}

\newcommand{\proof}{{\it Proof. }}

\newcommand{\N}{\ifmmode{{\Bbb N}}\else{\mbox{${\Bbb N}$}}\fi}
\newcommand{\R}{\ifmmode{{\Bbb R}}\else{\mbox{${\Bbb R}$}}\fi}

\begin{document}

\title{\textbf{Gain of regularity for a coupled system of generalized nonlinear Schr\"{o}dinger equations}}

\author{Raul Nina Mollisaca\thanks{Departamento de Matem\'{a}tica, Universidad
de Tarapaca, Arica,
Chile, raulnmollisaca@gmail.com.} ,
Mauricio Sep\'{u}lveda  Cort\'{e}s\thanks{Departamento de Ingenier\'ia Matem\'{a}tica and CI$^2$MA, Universidad
de Concepci\'on, Concepci\'{o}n,
Chile, masepulveda.cortes@gmail.} ,\\
Rodrigo V\'ejar Asem\thanks{Departamento de Matem\'aticas, Universidad
de La Serena, La Serena, Chile, rodrigo.vejar@userena.cl.} , 
Octavio  Vera Villagran\thanks{Departamento de Matem\'{a}tica, Universidad
de Tarapac\'a, Arica,
Chile, opverav@academicos.uta.cl. } }

\date{\today}
\maketitle

\maketitle

\begin{abstract}
\noindent
In this paper we study the smoothness properties of solutions to a one-dimensional coupled nonlinear Schr\"{o}dinger system equations that describes some physical phenomena such as propagation of polarized laser beams in birefringent Kerr medium in nonlinear optics. We show that the equations dispersive nature leads to a gain in regularity for the solution. In particular, if the initial data $(u_{0},\,v_{0})$ possesses certain regularity and sufficient decay as $|x|\rightarrow\infty,$ then the solution $(u(t),\,v(t))$ will be smoother than $(u_{0},\,v_{0})$ for $0 < t \leq T$ where $T$ is the existence time of the solution.
\end{abstract}
2010\noindent{\it Mathematics Subject Classification.} Primary 35Q53; Secondary 47J353

\renewcommand{\theequation}{\thesection.\arabic{equation}}
\setcounter{equation}{0}

\section{Introduction}
We consider the following initial value problem (IVP) for two coupled nonlinear Schr\"{o}dinger equations (NLS)
\begin{eqnarray}
\label{main1} \left\lbrace
\begin{array}{l}
iu_{t} + u_{xx} = (\alpha|u|^{2p} + \beta|u|^{q}|v|^{q + 2})u, \\
iv_{t} + v_{xx} = (\alpha|v|^{2p} + \beta|v|^{q}|u|^{q + 2})v,  \\
u(x,\,0) = u_{0}(x),\qquad v(x,\,0) = v_{0}(x),
\end{array}
\right.
\end{eqnarray}
where $\alpha,\beta,p,q,t,x\in\mathbb{R}$, $p>0$ and $q>0$. $u=u(x,t)$,  $v=v(x, t)$ are
complex unknown functions. The parameter $\beta$ is a real positive constant has to be interpreted as the birefringence intensity and describes the
coupling between the two component of the electric-field envelope.
Taking $\alpha = 1$ and $q = p - 1.$ The system \eqref{main1} leads to the following  system
\begin{eqnarray}
\label{main2} \left\lbrace
\begin{array}{l}
iu_{t} + u_{xx} = (|u|^{2p} + \beta|u|^{p - 1}|v|^{p + 1})u,    \\
iv_{t} + v_{xx} = (|v|^{2p} + \beta|v|^{p - 1}|u|^{p + 1})v,  \\
u(x,\,0) = u_{0}(x),\qquad v(x,\,0) = v_{0}(x).
\end{array}
\right.
\end{eqnarray}
In general a single mode fiber can support two distinct modes of polarization which are orthogonal to each other. This phenomenon is known as birefringence. Among these two modes one corresponds to the ordinary ray (O-ray) in which the refractive index of the medium is constant along every direction of the incident ray. The other is the extraordinary ray (E-ray) whose refractive index for the medium varies with the direction of the incident ray. In an ideal fiber these two modes are degenerate, while in a real fiber due to the fiber nonlinearity this degeneracy is broken and the phenomenon is known as {\it modal birefringence} \cite{agarwal}. The focusing nonlinear terms in \eqref{main2} describes the dependence of the refraction index of material on the electric field intensity and the birefringence effects.The study of the propagation of pulses  in nonlinear optical fiber has been of great interest in the last years. I. P. Kaminov \cite{kaminov} showed that single-mode optical fibers are not really ”single-mode” but actually bimodal due to the presence of birefringence which can deeply influence the way in which an optical evolves during the propagation along the fiber. Indeed, it can occur that the linear birefringence makes a pulse split in two, while nonlinear birefringence traps them together against splitting. In 1979, T. Kato \cite{ka} established a remarkable result for the regularizing property of solutions to the initial value problem for the KdV equation. He proved that if the initial function $u_{0}(x)\in L_{b}^{2} =  H^{2}(\mathbb{R})\cap L^{2}(e^{bx}dx)$ ($b > 0$), then the solution $u(t)$ is $C^{\infty}(\mathbb{R})$ for $t > 0.$ A main ingredient in the proof was the fact that formally the semigroup $e^{-t\partial_{x}^{3}}$ in $L_{b}^{2}$ is equivalent to $U_{b} = e^{-t(\partial x - b)^{3}}$ in $L^{2}$ when $t > 0.$ A number of results concerning the gain in regularity for various nonlinear evolution equations have appeared from different points view. It is shown that $C^{\infty}$ solutions $u(x,\,t)$ are obtained for all $t > 0$ if the initial data $\varphi(x)$ decays faster than polynomially on $\mathbb{R}^{+}$ and has certain initial Sobolev regularity. Here we quantify the gain of each derivative by the degree of vanishing of the initial data at infinity.
The gain of regularity for a higher order Schr\"{o}dinger equation has been also proved in \cite{tavo}.
In 1986, N. Hayashi {\it et al.} \cite{hayashi1, hayashi2} showed that if the initial data $u_{0}(x)$ decreases rapidly enough, then the solution of the Schr\"{o}dinger equation
\begin{align}
\label{main4}
\begin{cases}
iu_{t} + u_{xx} = \lambda|u|^{p}u,\quad x\in\mathbb{R},\ t\in \mathbb{R} \\
u(x,\,0) = u_{0}(x),\quad x\in\mathbb{R}
\end{cases}
\end{align}
with $\lambda\in\mathbb{R},$ $p > 1$ and $u = u(x,\,t)$ is a complex unknown function that becomes smooth for $t\neq 0,$ provided that the initial functions in $H^{1}(\mathbb{R})$ decay rapidly enough as $|x|\rightarrow +\infty.$ In 1986, N. Hayashi {\it et al.} \cite{hayashi2} studied the more general case. Indeed, they considered the initial value problem
\begin{align}
\label{main5}
\begin{cases}
iu_{t} +  \Delta u = f\left(|u|^{2}u\right),\quad x\in\mathbb{R}^{n},\ t\in \mathbb{R}  \\
u(x,\,0) = u_{0}(x),\quad x\in\mathbb{R}^{n}
\end{cases}
\end{align}
and they showed that if the initial $u_{0}(x)$ decreases sufficiently rapidly, then the solution of \eqref{main5} becomes smooth for $t\neq 0,$ provided the nonlinear term $f\left(|u|^{2}u\right)$ is smooth.  In J. C. Ceballos {\it et. al.} \cite{ceballo} the authors studied the gain in regularity for the system (case $p = 1$)
\begin{align}
\label{main6}
\begin{cases}
iu_{t} + u_{xx} + |u|^{2}u + \beta|v|^{2}u = 0 \\
iv_{t} + v_{xx} + |v|^{2}v + \beta|u|^{2}v = 0 \\
u(x,\,0) = u_{0}(x),\qquad v(x,\,0) = v_{0}(x)
\end{cases}
\end{align}
where $(x,\,t)\in \mathbb{R}\times \mathbb{R}.$ $u = u(x,\,t)$ and $v = v(x,\,t)$ are complex unknown functions and $\beta$ is a real positive constant which depends on the anisotropy of the fiber. C. R. Menyuk \cite{menyuk1, menyuk2} showed that the evolution of two orthogonal pulse envelopes in birefringent optical fiber is governed by the coupled nonlinear Schr\"{o}dinger system (\ref{main6}). If $\beta = 0$ the equations in (\ref{main6}) are two copies of a single nonlinear Schr\"{o}dinger equation which is integrable; when $\beta = 1,$ (\ref{main6}) is known as a Manakov system \cite{manak}. In all the other cases the situation is much more complicated from different points of view. The Cauchy problem for the system \eqref{main6} has been studied by many authors and over different point views see \cite{caze, ceballo, che, ma, rober} and references therein. Thus, it is natural to ask whether the equation \eqref{main2} has a gain in regularity. It might be expected that the Schr\"{o}dinger systems have an analogous regularizing effect as that of the \eqref{main4} equation. This is our motivation for the study of gain in regularity. Our aim in this paper is to show that the Schr\"{o}dinger systems have a regularizing effect. Indeed, that all solutions of finite energy to \eqref{main2} are smooth for $t\neq 0$ provided that the initial functions in $H^{1}(\mathbb{R})$ decay rapidly enough as $|x|\rightarrow +\infty.$ For to prove this result, the main tool is to use a operator $J$ defined by
\begin{align*}
Ju = e^{ix^{2}/4t}(2it)\partial_{x}(e^{-ix^{2}/4t}u) = (x + 2it\partial_{x})u.
\end{align*}
which commute with the operator of Schr\"{o}dinger. \\
\\
This paper is organized as follows: Before describing the main results, in Section 2 we briefly outline the notation and terminology to be used later on and we present some Theorems. In Section 3 we find estimates of finite energy. In Section 4 we find estimates for the operator $J$. In Section 5 we prove our main result. Our main result is\\

$\textbf{Main Theorem.}$
\begin{enumerate}
\item Let $p>1$ odd integer number, $(u_0, v_0)\in H^1(\mathbb{R})\times H^1(\mathbb{R})$ and $(x^nu_0, x^nv_0)\in L^2(\mathbb{R})\times L^2(\mathbb{R})$
for some $n\in\mathbb{N}$. Then, there exists a unique solution $(u(x, t), v(x, t))$ of (\ref{main2})
satisfying
\begin{align}
\label{ecudoc1}&(u,v)\in C_b(\mathbb{R}:H^1(\mathbb{R}))\times C_b(\mathbb{R}:H^1(\mathbb{R}))\\
\label{ecudoc2}&(J^mu,J^mv)\in C(\mathbb{R}:L^2(\mathbb{R}))\times C(\mathbb{R}:L^2(\mathbb{R})),\ \ m=1,2,\ldots,n.
\end{align}
Moreover $(u,v)$ satisfies the integral identities:\\
Densities Conservation
\begin{align*}
\|u\|_{L^{2}(\mathbb{R})} = \|u_{0}\|_{L^{2}(\mathbb{R})}\qquad and \qquad \|v\|_{L^{2}(\mathbb{R})} = \|v_{0}\|_{L^{2}(\mathbb{R})}.
\end{align*}
Energy Conservation
\begin{align*}
&\ \|u_{x}\|_{L^{2}(\mathbb{R})}^{2} + \|v_{x}\|_{L^{2}(\mathbb{R})}^{2} + \frac{1}{p + 1}\|u\|_{L^{2(p + 1)}(\mathbb{R})}^{2(p + 1)} + \frac{1}{p + 1}\|v\|_{L^{2(p + 1)}(\mathbb{R})}^{2(p + 1)}  \\
&\qquad  + \frac{2\beta}{p + 1}\int_{\mathbb{R}}|u|^{p + 1}|v|^{p + 1}dx   \\
= \ & \|u_{0x}\|_{L^{2}(\mathbb{R})}^{2} + \|v_{0x}\|_{L^{2}(\mathbb{R})}^{2} + \frac{1}{p + 1}\|u_{0}\|_{L^{2(p + 1)}(\mathbb{R})}^{2(p + 1)} + \frac{1}{p + 1}\|v_{0}\|_{L^{2(p + 1)}(\mathbb{R})}^{2(p + 1)}   \\
&\qquad  + \frac{2\beta}{p + 1}\int_{\mathbb{R}}|u_{0}|^{p + 1}|v_{0}|^{p + 1}dx.
\end{align*}
\item Let $(u_{0},\,v_{0})\in H^{1}(\mathbb{R})\times H^{1}(\mathbb{R})$ with $(x^{n}u_{0},\,x^{n}v_{0})\in L^{2}(\mathbb{R})\times L^{2}(\mathbb{R})$ and  $p > 1$ odd integer number. Then there exists a positive constant $C_{m}$ depending on $\|u_{0}\|_{H^{1}(\mathbb{R})},$ $\|v_{0}\|_{H^{1}(\mathbb{R})}$ and $\|x^{n}u_{0}\|_{L^{2}(\mathbb{R})},$ $\|x^{n}v_{0}\|_{L^{2}(\mathbb{R})}$ but independent of $t$ such that
\begin{align}\label{ecudoc3}
\|J^{m}u\|_{L^{2}(\mathbb{R})} \leq C_{m}e^{t}\quad {\rm and}\quad \|J^{m}v\|_{L^{2}(\mathbb{R})} \leq C_{m}e^{t},
\end{align}
for $m=1,2,...,n$.
\item Let $\beta < 1$, $p>1$ odd integer number, $(u_{0},\,v_{0})\in L^{2(p + 1)}(\mathbb{R})\times L^{2(p + 1)}(\mathbb{R})$
, $(u_{0},\,v_{0})\in H^1(\mathbb{R})\times H^1(\mathbb{R})$  and $(xu_{0}(x),\,xv_{0}(x))\in L^{2}(\mathbb{R})\times L^{2}(\mathbb{R}).$ Then, for all $t\neq0$ we have
\begin{align*}
\|u\|_{L^{\infty}(\mathbb{R})} \leq \frac{C}{t^{1/4}}\quad{\rm and}\quad \|v\|_{L^{\infty}(\mathbb{R})} \leq \frac{C}{t^{1/4}}
\end{align*}
and
\begin{align*}
\|u\|_{L^{2(p + 1)}(\mathbb{R})}  \leq \frac{C}{t^{2(p + 1)}}\quad {\rm and}\quad\|v\|_{L^{2(p + 1)}(\mathbb{R})}  \leq \frac{C}{t^{2(p + 1)}}.
\end{align*}
\end{enumerate}
Throughout this paper $C$ is a generic constant, not necessarily the same. at each occasion
(it will change from line to line), which depends in an increasing way on the indicated
quantities.

\renewcommand{\theequation}{\thesection.\arabic{equation}}
\setcounter{equation}{0}

\section{Preliminares}
In this section, some notations will be introduced. The Lebesgue space
 $L^p(\mathbb{R})$, $1\le p\le+\infty$, with norm denoted by $\|\cdot\|_{L^p}$ will be used.
 Let $m$ be a positive integer, i.e., $m\in\mathbb{N}$, the usual Sobolev spaces
$W^{m,p}(\mathbb{R})=\{u\in L^p(\mathbb{R})\,:\, \|D^\alpha
u\|_{L^p}<+\infty,\ \forall |\alpha|\le m\}$, with norm denoted by
$\|\cdot\|_{W^{m,p}}$ is considered. When $p=2$,  it is stablished
$H^m(\mathbb{R}):=W^{m,2}(\mathbb{R})$  denoting the respective norm by
$\|\cdot\|_{H^m}$.  \\
For any interval $I$ of $\mathbb{R}$ and any Banach space $X$ with the norm $\|\cdot \|_{X},$ we denote $C(I:\,X)$ (respectively $C_{b}(I:\,X))$ the space of continuous (respectively bounded continuous) functions from $I$ to $X$. We denote $C^{k}(I:\,X)$ ($k\geq 1$) the space of $k$-times continuously differentiable functions from $I$ to $X$. \\
If $X$ is a Banach space, $L^p(I:\,X)$  indicates the space of valued
functions in $X$ defined on the interval $I$ that are integrable
in the Bochner sense, and its norm will be denoted by
$\|\cdot\|_{L^p(X)}$.

{\bf Remark.}
We only consider the case $t>0.$ The case $t<0$ can be treated analogously. \\
\\
The following results are going to be used several times from now on.
\begin{theorem}
\label{lemma201}
(The Gagliardo-Nirenberg inequality) Let $q,\,r$ be any real numbers satisfying $1\leq q,\,r\leq +\infty$ and let $j$ and $m$ be non-negative integers such that $j\leq m.$ Then
\begin{align}
\label{201}\|D^{j}u\|_{L^{p}(\mathbb{R})} \leq C\|D^{m}u\|_{L^{r}(\mathbb{R})}^{a}\|u\|_{L^{q}(\mathbb{R})}^{1 - a},
\end{align}
where
\begin{align}
\frac{1}{p} = j + a\left(\frac{1}{r} - m\right) + (1 - a)\frac{1}{q},
\end{align}
for all $a$ in the interval $\frac{j}{m}\leq a\leq 1,$ and $C$ is a positive constant depending only $m,$ $j,$ $q,$ $r,$ and $a$.
\end{theorem}
\proof See \cite{friedman}.\\
How consequence the theorem \eqref{lemma201} we have the following corollary.
\begin{corollary}
For all $u,\ v\in H^{1}(\mathbb{R})$ we have
\begin{equation}
\begin{aligned}
\label{lemma202}
\|u\|_{L^{\infty}(\mathbb{R})}^{2} \leq C\|u\|_{L^{2}(\mathbb{R})}\|u_{x}\|_{L^{2}(\mathbb{R})}
\end{aligned}
\end{equation}
\end{corollary}
\begin{lemma}
\label{lemma203}
Let $u$ and $v$ be the solutions of \eqref{main2}, then we have
\begin{align*}
(a)\ \frac{d}{dt}(|u|^{2}) = 2Im(u\overline{u}_{xx})\qquad {\rm and}\qquad (b)\ \frac{d}{dt}(|v|^{2}) = 2Im(v\overline{v}_{xx})
\end{align*}
\end{lemma}
\proof Caso $(a)$ Multiplying \eqref{main2}$_{1}$ by $\overline{u}$ we have
\begin{align}
\label{204}i\overline{u}\,u_{t} + \overline{u}\,u_{xx} = (|u|^{2p} + \beta|u|^{p - 1}|v|^{p + 1})|u|^{2}.
\end{align}
Applying conjugate in the above equation
\begin{align}
\label{205}-iu\,\overline{u}_{t} + u\,\overline{u}_{xx} = (|u|^{2p} + \beta|u|^{p - 1}|v|^{p + 1})|u|^{2}.
\end{align}
Subtracting \eqref{204} with \eqref{205} we obtain
\begin{align}
&i\frac{d}{dt}(|u|^{2}) + \overline{u}u_{xx} - u\overline{u}_{xx} = 0 \iff i\frac{d}{dt}(|u|^{2}) = u\overline{u}_{xx} - \overline{u}u_{xx} = 2iIm\,u\overline{u}_{xx} \nonumber  \\
\label{206}\iff & \qquad\qquad\frac{d}{dt}(|u|^{2}) = 2Im\,u\overline{u}_{xx}.
\end{align}
In a similar way we obtain $(b).$ \\
\\
The operator $J$ commutes with the operator $L$ defined by $L = (i\partial_{t} + \partial_{x}^{2}),$ that is, $LJ - JL\equiv [L,\,J] = 0.$ In general,
\begin{align*}
J^m(u) = e^{ix^{2}/4t}(2it)^m\partial_x^m(e^{-ix^2/4t}u)=(x+2it\partial_x)^m u, \hspace{0.3cm}m\in\mathbb{N},
\end{align*}
where $J^m(u)=J(J^{m-1}u)$, $m\in \mathbb{N}$. Hence, applying $J^m$ to be equations  \eqref{main2}
\begin{eqnarray}
\label{main3} \left.
\begin{array}{l}
i(J^{m}u)_{t} + (J^{m}u)_{xx} = J^{m}(|u|^{2p}u) + \beta J^{m}(|u|^{p - 1}|v|^{p + 1}u),\quad   t > 0,\quad x\in\mathbb{R}   \\
i(J^{m}u)_{t} + (J^{m}u)_{xx} = J^{m}(|v|^{2p}v) + \beta J^{m}(|v|^{p - 1}|u|^{p + 1}v),\quad   t > 0,\quad x\in\mathbb{R}   \\
J^{m}u(x,\,0) = x^{m}u_{0}(x),\quad J^{m}u(x,\,0) = x^{m}v_{0}(x),\quad x\in\mathbb{R}.
\end{array}
\right.
\end{eqnarray}
which allows us to get the estimates of the Main Theorem.

\renewcommand{\theequation}{\thesection.\arabic{equation}}
\setcounter{equation}{0}

\section{Finite energy solutions}

We begin obtaining estimates for the norm-$H^{1}(\mathbb{R}).$ To estimate in $H^{1}(\mathbb{R})$ is important because in $\mathbb{R}$ we have the following Sobolev immersion (continuous injection) $H^{1}(\mathbb{R})\hookrightarrow L^{\infty}(\mathbb{R})$.
\begin{lemma}
\label{lemma301}
(Densities conservation). Let $(u,\,v)$ be the solution to \eqref{main2}. Let $(u_{0},\,v_{0}) \in L^{2}(\mathbb{R})\times L^{2}(\mathbb{R}),$ then
\begin{align}
\label{301}\|u\|_{L^{2}(\mathbb{R})} = \|u_{0}\|_{L^{2}(\mathbb{R})}\qquad and \qquad \|v\|_{L^{2}(\mathbb{R})} = \|v_{0}\|_{L^{2}(\mathbb{R})}.
\end{align}
\end{lemma}
\proof Integrating \eqref{206} over $x\in \mathbb{R}$ and using integrating by parts
\begin{align*}
\frac{d}{dt}\|u\|_{L^{2}(\mathbb{R})}^2 = -2Im\int_{\mathbb{R}}|u_{x}|^{2}dx = 0 \iff \frac{d}{dt}\|u\|_{L^{2}(\mathbb{R})}^2 = 0
\end{align*}
and integrating over $t\in [0,\,T]$ we get the first term. Similarly for $v.$ The lemma follows.
\begin{lemma}
(Energy conservation). Let $(u,\,v)$ be the solution to \eqref{main2}. Let $(u_{0},\,v_{0})\in L^{2(p + 1)}(\mathbb{R})\times L^{2(p + 1)}(\mathbb{R})$, $(u_{0},\,v_{0}) \in H^{1}(\mathbb{R})\times H^{1}(\mathbb{R})$ and $p>1$ odd integer number, then
\label{lemma302}
\begin{align}
&\ \|u_{x}\|_{L^{2}(\mathbb{R})}^{2} + \|v_{x}\|_{L^{2}(\mathbb{R})}^{2} + \frac{1}{p + 1}\|u\|_{L^{2(p + 1)}(\mathbb{R})}^{2(p + 1)} + \frac{1}{p + 1}\|v\|_{L^{2(p + 1)}(\mathbb{R})}^{2(p + 1)} \nonumber \\
&\qquad  + \frac{2\beta}{p + 1}\int_{\mathbb{R}}|u|^{p + 1}|v|^{p + 1}dx  \nonumber \\
= \ & \|u_{0x}\|_{L^{2}(\mathbb{R})}^{2} + \|v_{0x}\|_{L^{2}(\mathbb{R})}^{2} + \frac{1}{p + 1}\|u_{0}\|_{L^{2(p + 1)}(\mathbb{R})}^{2(p + 1)} + \frac{1}{p + 1}\|v_{0}\|_{L^{2(p + 1)}(\mathbb{R})}^{2(p + 1)}  \nonumber \\
\label{302}&\qquad  + \frac{2\beta}{p + 1}\int_{\mathbb{R}}|u_{0}|^{p + 1}|v_{0}|^{p + 1}dx.
\end{align}
\end{lemma}
\proof Applying $\partial_{x}$ over \eqref{main2}$_{1}$ we have
\begin{align*}
&iu_{xt} + u_{xxx}= (|u|^{2p})_{x}u+\beta(|u|^{p - 1}|v|^{p + 1})_xu+ (|u|^{2p}+\beta|u|^{p-1}|v|^{p+1})u_x.
\end{align*}
Multiplying the above equation by $\overline{u}_{x}$ we have
\begin{align*}
&i\overline{u}_{x}u_{xt} + \overline{u}_{x}u_{xxx}=(|u|^{2p})_{x}u\overline{u}_{x}+\beta(|u|^{p-1}|v|^{p+1})_xu\overline{u}_x
 + |u|^{2p}|u_x|^2+\beta |u|^{p-1}|v|^{p+1}|u_x|^2.
\end{align*}
Integrating over $\mathbb{R}$ yields
\begin{align*}
i\int_{\mathbb{R}}\overline{u}_{x}u_{xt}dx + \int_{\mathbb{R}}\overline{u}_{x}u_{xxx}dx =& \int_{\mathbb{R}}(|u|^{2p})_{x}u\overline{u}_{x}dx + \int_{\mathbb{R}}|u|^{2p}|u_{x}|^{2}dx
 + \beta\int_{\mathbb{R}}(|u|^{p - 1})_{x}|v|^{p + 1}u\overline{u}_{x}dx\\
& + \beta\int_{\mathbb{R}}|u|^{p - 1}(|v|^{p + 1})_{x}u\overline{u}_{x}dx + \beta\int_{\mathbb{R}}|u|^{p - 1}|v|^{p + 1}|u_{x}|^{2}dx.
\end{align*}
Integrating by parts the second term in the left hand side we obtain
\begin{align*}
i\int_{\mathbb{R}}\overline{u}_{x}u_{xt}dx - \int_{\mathbb{R}}|u_{xx}|^{2}dx =& \int_{\mathbb{R}}(|u|^{2p})_{x}u\overline{u}_{x}dx + \int_{\mathbb{R}}|u|^{2p}|u_{x}|^{2}dx
 + \beta\int_{\mathbb{R}}(|u|^{p - 1})_{x}|v|^{p + 1}u\overline{u}_{x}dx\\&  + \beta\int_{\mathbb{R}}|u|^{p - 1}(|v|^{p + 1})_{x}u\overline{u}_{x}dx + \beta\int_{\mathbb{R}}|u|^{p - 1}|v|^{p + 1}|u_{x}|^{2}dx.
\end{align*}
Applying conjugate
\begin{align*}
-i\int_{\mathbb{R}}u_{x}\overline{u}_{xt}dx - \int_{\mathbb{R}}|u_{xx}|^{2}dx =& \int_{\mathbb{R}}(|u|^{2p})_{x}\overline{u\overline{u}_{x}}dx + \int_{\mathbb{R}}|u|^{2p}|u_{x}|^{2}dx
 + \beta\int_{\mathbb{R}}(|u|^{p - 1})_{x}|v|^{p + 1}\overline{u\overline{u}_{x}}dx \\& + \beta\int_{\mathbb{R}}|u|^{p - 1}(|v|^{p + 1})_{x}\overline{u\overline{u}_{x}}dx + \beta\int_{\mathbb{R}}|u|^{p - 1}|v|^{p + 1}|u_{x}|^{2}dx.
\end{align*}
Subtracting the above equations
\begin{align*}
\frac{d}{dt}\int_{\mathbb{R}}|u_{x}|^{2}dx =& 2Im\int_{\mathbb{R}}(|u|^{2p})_{x}u\overline{u}_{x}dx
+ 2\beta Im\int_{\mathbb{R}}(|u|^{p - 1})_{x}|v|^{p + 1}u\overline{u}_{x}dx \\&+ 2\beta Im\int_{\mathbb{R}}|u|^{p - 1}(|v|^{p + 1})_{x}u\overline{u}_{x}dx.
\end{align*}
Integrating by parts
\begin{align*}
\frac{d}{dt}\int_{\mathbb{R}}|u_{x}|^{2}dx =& - 2Im\int_{\mathbb{R}}|u|^{2p}u\overline{u}_{xx}dx
+ 2\beta Im\int_{\mathbb{R}}(|u|^{p - 1})_{x}|v|^{p + 1}u\overline{u}_{x}dx\\& + 2\beta Im\int_{\mathbb{R}}|u|^{p - 1}(|v|^{p + 1})_{x}u\overline{u}_{x}dx.
\end{align*}
Then
\begin{align*}
&\frac{d}{dt}\int_{\mathbb{R}}|u_{x}|^{2}dx = - \int_{\mathbb{R}}|u|^{2p}\ (2Im u\overline{u}_{xx})dx + 2\beta Im\int_{\mathbb{R}}(|u|^{p - 1}|v|^{p + 1})_{x}u\overline{u}_{x}dx.
\end{align*}
Using Lemma \ref{lemma203} we obtain
\begin{align*}
&\frac{d}{dt}\int_{\mathbb{R}}|u_{x}|^{2}dx = - \int_{\mathbb{R}}|u|^{2p}\ \frac{d}{dt}(|u|^{2})dx - 2\beta Im\int_{\mathbb{R}}|u|^{p - 1}|v|^{p + 1}u\overline{u}_{xx}dx
\end{align*}
where
\begin{align*}
&\frac{d}{dt}\int_{\mathbb{R}}|u_{x}|^{2}dx + \frac{1}{p + 1}\frac{d}{dt}\int_{\mathbb{R}}|u|^{2(p + 1)}dx + \beta \int_{\mathbb{R}}|u|^{p - 1}|v|^{p + 1}\frac{d}{dt}(|u|^{2})dx = 0.
\end{align*}
Hence
\begin{align*}
&\frac{d}{dt}\left[\int_{\mathbb{R}}|u_{x}|^{2}dx + \frac{1}{p + 1}\int_{\mathbb{R}}|u|^{2(p + 1)}dx\right] + \beta \int_{\mathbb{R}}|u|^{p - 1}|v|^{p + 1}\frac{d}{dt}(|u|^{2})dx = 0.
\end{align*}
Thus
\begin{align*}
\frac{d}{dt}\left[\int_{\mathbb{R}}|u_{x}|^{2}dx + \frac{1}{p + 1}\int_{\mathbb{R}}|u|^{2(p + 1)}dx\right]
+ \beta \int_{\mathbb{R}}\frac{2}{p+1}\frac{d}{dt}((|u|^{2})^{\frac{p+1}{2}})|v|^{p + 1}dx = 0.
\end{align*}
Thereby
\begin{align}
\label{303}\frac{d}{dt}\left[\int_{\mathbb{R}}|u_{x}|^{2}dx + \frac{1}{p + 1}\int_{\mathbb{R}}|u|^{2(p + 1)}dx\right] + \frac{2\beta}{p + 1}\int_{\mathbb{R}}\frac{d}{dt}(|u|^{p + 1})|v|^{p + 1}dx = 0.
\end{align}
Similarly we have
\begin{align}
\label{304}\frac{d}{dt}\left[\int_{\mathbb{R}}|v_{x}|^{2}dx + \frac{1}{p + 1}\int_{\mathbb{R}}|v|^{2(p + 1)}dx\right] + \frac{2\beta}{p + 1}\int_{\mathbb{R}}\frac{d}{dt}(|v|^{p + 1})|u|^{p + 1}dx = 0.
\end{align}
Adding \eqref{303} with \eqref{304} we obtain
\begin{align*}
&\frac{d}{dt}\left[\int_{\mathbb{R}}|u_{x}|^{2}dx + \int_{\mathbb{R}}|v_{x}|^{2}dx + \frac{1}{p + 1}\int_{\mathbb{R}}|u|^{2(p + 1)}dx + \frac{1}{p + 1}\int_{\mathbb{R}}|v|^{2(p + 1)}dx\right. \\
&\qquad \left. + \frac{2\beta}{p + 1}\int_{\mathbb{R}}|u|^{p + 1}|v|^{p + 1}dx\right] = 0.
\end{align*}
Therefore, we obtain
\begin{align}
&\frac{d}{dt}\left[\|u_{x}\|_{L^{2}(\mathbb{R})}^{2} + \|v_{x}\|_{L^{2}(\mathbb{R})}^{2} + \frac{1}{p + 1}\|u\|_{L^{2(p + 1)}(\mathbb{R})}^{2(p + 1)} + \frac{1}{p + 1}\|v\|_{L^{2(p + 1)}(\mathbb{R})}^{2(p + 1)}\right. \nonumber \\
\label{305}&\qquad \left. + \frac{2\beta}{p + 1}\int_{\mathbb{R}}|u|^{p + 1}|v|^{p + 1}dx\right] = 0.
\end{align}
Integrating over $t\in[0,\,T]$ the result follows. \\
\\
{\bf Remark.} Integrating \eqref{305} over $t\in [0,\,T]$ we have
\begin{align*}
&\|u_{x}\|_{L^{2}(\mathbb{R})}^{2} + \|v_{x}\|_{L^{2}(\mathbb{R})}^{2} + \frac{1}{p + 1}\|u\|_{L^{2(p + 1)}(\mathbb{R})}^{2(p + 1)} + \frac{1}{p + 1}\|v\|_{L^{2(p + 1)}(\mathbb{R})}^{2(p + 1)} \nonumber \\
&\qquad  + \frac{2\beta}{p + 1}\int_{\mathbb{R}}|u|^{p + 1}|v|^{p + 1}dx \nonumber \\
=\ & \|u_{0x}\|_{L^{2}(\mathbb{R})}^{2} + \|v_{0x}\|_{L^{2}(\mathbb{R})}^{2} + \frac{1}{p + 1}\|u_{0}\|_{L^{2(p + 1)}(\mathbb{R})}^{2(p + 1)} + \frac{1}{p + 1}\|v_{0}\|_{L^{2(p + 1)}(\mathbb{R})}^{2(p + 1)} \nonumber \\
&\qquad  + \frac{2\beta}{p + 1}\int_{\mathbb{R}}|u_{0}|^{p + 1}|v_{0}|^{p + 1}dx.
\end{align*}
Then
\begin{align*}
&\|u_{x}\|_{L^{2}(\mathbb{R})}^{2} + \|v_{x}\|_{L^{2}(\mathbb{R})}^{2} + \frac{1}{p + 1}\|u\|_{L^{2(p + 1)}(\mathbb{R})}^{2(p + 1)} + \frac{1}{p + 1}\|v\|_{L^{2(p + 1)}(\mathbb{R})}^{2(p + 1)} \nonumber \\
=\ & \|u_{0x}\|_{L^{2}(\mathbb{R})}^{2} + \|v_{0x}\|_{L^{2}(\mathbb{R})}^{2} + \frac{1}{p + 1}\|u_{0}\|_{L^{2(p + 1)}(\mathbb{R})}^{2(p + 1)} + \frac{1}{p + 1}\|v_{0}\|_{L^{2(p + 1)}(\mathbb{R})}^{2(p + 1)} \nonumber \\
&\qquad  + \frac{2\beta}{p + 1}\int_{\mathbb{R}}|u_{0}|^{p + 1}|v_{0}|^{p + 1}dx - \frac{2\beta}{p + 1}\int_{\mathbb{R}}|u|^{p + 1}|v|^{p + 1}dx \nonumber \end{align*}
Then using the Young inequality, we have
\begin{align*}
&\|u_{x}\|_{L^{2}(\mathbb{R})}^{2} + \|v_{x}\|_{L^{2}(\mathbb{R})}^{2} + \frac{1}{p + 1}\|u\|_{L^{2(p + 1)}(\mathbb{R})}^{2(p + 1)} + \frac{1}{p + 1}\|v\|_{L^{2(p + 1)}(\mathbb{R})}^{2(p + 1)} \nonumber \\
\leq \ & \|u_{0x}\|_{L^{2}(\mathbb{R})}^{2} + \|v_{0x}\|_{L^{2}(\mathbb{R})}^{2} + \frac{1}{p + 1}\|u_{0}\|_{L^{2(p + 1)}(\mathbb{R})}^{2(p + 1)} + \frac{1}{p + 1}\|v_{0}\|_{L^{2(p + 1)}(\mathbb{R})}^{2(p + 1)} \nonumber \\
&+ \frac{2\beta}{p + 1}\int_{\mathbb{R}}|u_{0}|^{p + 1}|v_{0}|^{p + 1}dx + \frac{\beta}{p + 1}\int_{\mathbb{R}}|u|^{2(p + 1)}dx + \frac{\beta}{p + 1}\int_{\mathbb{R}}|v|^{2(p + 1)}dx \nonumber \\
\leq \ & \|u_{0x}\|_{L^{2}(\mathbb{R})}^{2} + \|v_{0x}\|_{L^{2}(\mathbb{R})}^{2} + \frac{1}{p + 1}\|u_{0}\|_{L^{2(p + 1)}(\mathbb{R})}^{2(p + 1)} + \frac{1}{p + 1}\|v_{0}\|_{L^{2(p + 1)}(\mathbb{R})}^{2(p + 1)} \nonumber \\
& + \frac{\beta}{p + 1}\|u_{0}\|_{L^{2(p + 1)}(\mathbb{R})}^{2(p + 1)}+ \frac{\beta}{p + 1}\|v_{0}\|_{L^{2(p + 1)}(\mathbb{R})}^{2(p + 1)}\nonumber \\
& + \frac{\beta}{p + 1}\|u\|_{L^{2(p + 1)}(\mathbb{R})}^{2(p + 1)}+ \frac{\beta}{p + 1}\|v\|_{L^{2(p + 1)}(\mathbb{R})}^{2(p + 1)}.
\end{align*}
Thus, we obtain
\begin{align}
&\|u_{x}\|_{L^{2}(\mathbb{R})}^{2} + \|v_{x}\|_{L^{2}(\mathbb{R})}^{2} + \frac{1 - \beta}{p + 1}\|u\|_{L^{2(p + 1)}(\mathbb{R})}^{2(p + 1)} + \frac{1 - \beta}{p + 1}\|v\|_{L^{2(p + 1)}(\mathbb{R})}^{2(p + 1)} \nonumber \\
\label{306}\leq \ & \|u_{0x}\|_{L^{2}(\mathbb{R})}^{2} + \|v_{0x}\|_{L^{2}(\mathbb{R})}^{2} + \frac{1 + \beta}{p + 1}\|u_{0}\|_{L^{2(p + 1)}(\mathbb{R})}^{2(p + 1)} + \frac{1 + \beta}{p + 1}\|v_{0}\|_{L^{2(p + 1)}(\mathbb{R})}^{2(p + 1)}.
\end{align}
We follows that, if $0<\beta<1$ then, for all $1<p<+\infty$ we have $u,v\in L^{2(p+1)}(\mathbb{R})$.

\renewcommand{\theequation}{\thesection.\arabic{equation}}
\setcounter{equation}{0}

\section{A priori estimates}

In the proof stated below it is shown that $L^{\infty}(\mathbb{R})$ estimates of solutions lead to obtain a priori estimates of $Ju.$ We estimate a Gronwall's inequality type and we establish decay of perturbed solutions.
\begin{theorem}
\label{theorem401}
Let $\beta < 1$, $p>1$ odd integer number $(u_{0},\,v_{0})\in L^{2(p + 1)}(\mathbb{R})\times L^{2(p + 1)}(\mathbb{R})$, $(u_0,\, v_0)\in H^1(\mathbb{R})\times H^1(\mathbb{R})$ and $(xu_{0}(x),\,xv_{0}(x))\in L^{2}(\mathbb{R})\times L^{2}(\mathbb{R}).$ Then
\end{theorem}
\begin{align}
\label{401}
\|u\|_{L^{\infty}(\mathbb{R})} \leq \frac{C}{t^{1/4}},\quad{\rm }\quad \|v\|_{L^{\infty}(\mathbb{R})} \leq \frac{C}{t^{1/4}}
\end{align}
and
\begin{align}
\label{4011}\|u\|_{L^{2(p + 1)}(\mathbb{R})}  \leq \frac{C}{t^{2(p + 1)}},\quad {\rm }\quad\|v\|_{L^{2(p + 1)}(\mathbb{R})}  \leq \frac{C}{t^{2(p + 1)}}.
\end{align}
\\
\proof We rewrite the equation \eqref{main2}$_{1}$ as
\begin{align*}
Lu = |u|^{2p}u + \beta|u|^{p - 1}|v|^{p + 1}u.
\end{align*}
Then we consider the operator $J$ such that $LJ = JL,$ thus
\begin{align*}
L(Ju) = J(|u|^{2p}u) + \beta J(|u|^{p - 1}|v|^{p + 1}u).
\end{align*}
Thus
\begin{align*}
i(Ju)_{t} + (Ju)_{xx} = J(|u|^{2p}u) + \beta J(|u|^{p - 1}|v|^{p + 1}u).
\end{align*}
Multiplying the above equation by $(\overline{Ju})$ we have
\begin{align}
\label{402}i(\overline{Ju})(Ju)_{t} + (\overline{Ju})(Ju)_{xx} = J(|u|^{2p}u)(\overline{Ju}) + \beta J(|u|^{p - 1}|v|^{p + 1}u)(\overline{Ju}).
\end{align}
Applying conjugate
\begin{align}
\label{403}-i(Ju)(\overline{Ju})_{t} + (Ju)(\overline{Ju})_{xx} = \overline{J(|u|^{2p}u)(\overline{Ju})} + \beta \overline{J(|u|^{p - 1}|v|^{p + 1}u)(\overline{Ju})}.
\end{align}
Integrating over $x\in\mathbb{R},$ subtracting and dividing by $i,$ we have
\begin{align}
\label{404}\frac{d}{dt}\|Ju\|_{L^{2}(\mathbb{R})}^{2} = 2Im\int_{\mathbb{R}}J(|u|^{2p}u)(\overline{Ju})dx + 2\beta Im\int_{\mathbb{R}}J(|u|^{p - 1}|v|^{p + 1}u)(\overline{Ju})dx.
\end{align}
We estimate the two terms on the right hand side of \eqref{404}
\begin{align*}
J(|u|^{2p}u) = x|u|^{2p}u + 2it(|u|^{2p})_{x}u + 2it|u|^{2p}u_{x},\qquad {\rm and}\qquad \overline{Ju} = x\overline{u} - 2it\overline{u}_{x}.
\end{align*}
Then
\begin{align*}
J(|u|^{2p}u)(\overline{Ju}) = &\ x^{2}|u|^{2p + 2} - 2itx|u|^{2p}u\overline{u}_{x} + 2itx(|u|^{2p})_{x}|u|^{2} + 4t^{2}(|u|^{2p})_{x}u\overline{u}_{x}  \\
& \ + 2itx|u|^{2p}\overline{u}u_{x} + 4t^{2}|u|^{2p}|u_{x}|^{2}.
\end{align*}
Hence
\begin{align*}
J(|u|^{2p}u)(\overline{Ju}) = &\ x^{2}|u|^{2p + 2} + 4t^{2}|u|^{2p}|u_{x}|^{2} + 4t^{2}(|u|^{2p})_{x}u\overline{u}_{x} \\
&\ + 2itx(|u|^{2p})_{x}|u|^{2} + 2itx|u|^{2p}(u\overline{u}_{x} - \overline{u\overline{u}_{x}}).
\end{align*}
Taking Imaginary part
\begin{align*}
Im[J(|u|^{2p}u)(\overline{Ju})] = 2tx(|u|^{2p})_{x}|u|^{2} + 4t^{2}Im[(|u|^{2p})_{x}u\overline{u}_{x}].
\end{align*}
Integrating over $x\in\mathbb{R}$
\begin{align*}
2Im\int_{\mathbb{R}}J(|u|^{2p}u)(\overline{Ju})dx = &\ 4t\int_{\mathbb{R}}x\ \frac{d}{dx}(|u|^{2p})\ |u|^{2}dx + 8t^{2}Im\int_{\mathbb{R}}\frac{d}{dx}(|u|^{2p})\ u\overline{u}_{x}dx \\
=&\ \frac{4pt}{(p + 1)}\int_{\mathbb{R}}x\ \frac{d}{dx}\left[(|u|^{2p})^{\frac{(p + 1)}{p}}\right]dx - 8t^{2}Im\int_{\mathbb{R}}|u|^{2p}u\overline{u}_{xx}dx \\
=&\ -\frac{4pt}{(p + 1)}\int_{\mathbb{R}}|u|^{2(p + 1)}dx - 4t^{2}\int_{\mathbb{R}}|u|^{2p}\ 2Im(u\overline{u}_{xx})dx.
\end{align*}
Using Lemma \ref{lemma203}
\begin{align*}
&\ 2Im\int_{\mathbb{R}}J(|u|^{2p}u)(\overline{Ju})dx = -\frac{4pt}{(p + 1)}\int_{\mathbb{R}}|u|^{2(p + 1)}dx - 4t^{2}\int_{\mathbb{R}}|u|^{2p}\ \frac{d}{dt}(|u|^{2})dx \\
= &\ -\frac{4pt}{(p + 1)}\int_{\mathbb{R}}|u|^{2(p + 1)}dx - \frac{4t^{2}}{(p + 1)}\frac{d}{dt}\int_{\mathbb{R}}|u|^{2(p + 1)}dx.
\end{align*}
On the other hand, using the identity
\begin{align*}
- \frac{4t^{2}}{(p + 1)}\frac{d}{dt}\int_{\mathbb{R}}|u|^{2(p + 1)}dx = \frac{8t}{(p + 1)}\int_{\mathbb{R}}|u|^{2(p + 1)}dx - \frac{d}{dt}\left[\frac{4t^{2}}{(p + 1)}\int_{\mathbb{R}}|u|^{2(p + 1)}dx\right]
\end{align*}
we have
\begin{align*}
\lefteqn{2Im\int_{\mathbb{R}}J(|u|^{2p}u)(\overline{Ju})dx } \\
=\ &  -\frac{4pt}{(p + 1)}\|u\|_{L^{2(p + 1)}(\mathbb{R})}^{2(p + 1)} + \frac{8t}{(p + 1)}\|u\|_{L^{2(p + 1)}(\mathbb{R})}^{2(p + 1)} - \frac{d}{dt}\left[\frac{4t^{2}}{(p + 1)}\|u\|_{L^{2(p + 1)}(\mathbb{R})}^{2(p + 1)}\right].
\end{align*}
Thereby
\begin{align}
\label{405}2Im\int_{\mathbb{R}}J(|u|^{2p}u)(\overline{Ju})dx = \frac{4(2 - p)t}{(p + 1)}\|u\|_{L^{2(p + 1)}(\mathbb{R})}^{2(p + 1)} - \frac{d}{dt}\left[\frac{4t^{2}}{(p + 1)}\|u\|_{L^{2(p + 1)}(\mathbb{R})}^{2(p + 1)}\right].
\end{align}
Now we estimate the second term in the right hand side of \eqref{404}, that is,
\begin{align*}
2\beta Im\int_{\mathbb{R}}J(|u|^{p - 1}|v|^{p + 1}u)(\overline{Ju})dx
\end{align*}
\begin{align*}
J(|u|^{p - 1}|v|^{p + 1}u) = \ & x|u|^{p - 1}|v|^{p + 1}u + 2it(|u|^{p - 1})_{x}|v|^{p + 1}u + 2it|u|^{p - 1}(|v|^{p + 1})_{x}u \\
&+ 2it|u|^{p - 1}|v|^{p + 1}u_{x}.
\end{align*}
Then
\begin{align*}
J(|u|^{p - 1}|v|^{p + 1}u)\overline{Ju}
= \ & x^{2}|u|^{p - 1}|v|^{p + 1}|u|^{2} + 4t^{2}|u|^{p - 1}|v|^{p + 1}|u_{x}|^{2} \\
\ & + 2itx(|u|^{p - 1})_{x}|v|^{p + 1}|u|^{2} + 2itx|u|^{p - 1}(|v|^{p + 1})_{x}|u|^{2} \\
\ & - 2itx|u|^{p - 1}|v|^{p + 1}(u\overline{u}_{x} - \overline{u\overline{u}_{x}}) + 4t^{2}(|u|^{p - 1}|v|^{p + 1})_{x}u\overline{u}_{x} \\
= \ & x^{2}|u|^{p - 1}|v|^{p + 1}|u|^{2} + 4t^{2}|u|^{p - 1}|v|^{p + 1}|u_{x}|^{2} \\
\ & + 2itx(|u|^{p - 1}|v|^{p + 1})_{x}|u|^{2} + 4tx|u|^{p - 1}|v|^{p + 1}Im(u\overline{u}_{x})  \\
\ & + 4t^{2}(|u|^{p - 1}|v|^{p + 1})_{x}u\overline{u}_{x}.
\end{align*}
Taking the imaginary part
\begin{align*}
Im[J(|u|^{p - 1}|v|^{p + 1}u)\overline{Ju}] = 2tx(|u|^{p - 1}|v|^{p + 1})_{x}|u|^{2} + 4t^{2}(|u|^{p - 1}|v|^{p + 1})_{x}Im(u\overline{u}_{x}).
\end{align*}
Integrating over $x\in\mathbb{R}$
\begin{align*}
\lefteqn{2\beta Im\int_{\mathbb{R}}J(|u|^{p - 1}|v|^{p + 1}u)\overline{Ju}\,dx } \\
= \  &  4\beta t\int_{\mathbb{R}}x\ \frac{d}{dx}(|u|^{p - 1}|v|^{p + 1})\ |u|^{2}\,dx + 8\beta t^{2}\int_{\mathbb{R}}\frac{d}{dx}(|u|^{p - 1}|v|^{p + 1})\ Im(u\overline{u}_{x})\,dx \\
= \ & - 4\beta t\int_{\mathbb{R}}|u|^{p - 1}|v|^{p + 1}|u|^{2}\,dx - 4\beta t\int_{\mathbb{R}}x|u|^{p - 1}|v|^{p + 1}\ \frac{d}{dx}(|u|^{2})dx \\
& \ - 4\beta t^{2}\int_{\mathbb{R}}|u|^{p - 1}|v|^{p + 1}\ 2Im(u\overline{u}_{xx})\,dx.
\end{align*}
Hence using Lemma \ref{lemma203} we have
\begin{align}
\lefteqn{2\beta Im\int_{\mathbb{R}}J(|u|^{p - 1}|v|^{p + 1}u)\overline{Ju}dx } \nonumber \\
= \ & - 4\beta t\int_{\mathbb{R}}|u|^{p + 1}|v|^{p + 1}dx - 4\beta t\int_{\mathbb{R}}x|u|^{p - 1}\frac{d}{dx}(|u|^{2})\ |v|^{p + 1}dx \nonumber \\
& \ - 4\beta t^{2}\int_{\mathbb{R}}|u|^{p - 1}|v|^{p + 1}\frac{d}{dt}(|u|^{2})dx \nonumber \\
= \ & - 4\beta t\int_{\mathbb{R}}|u|^{p + 1}|v|^{p + 1}dx - \frac{8\beta t}{(p + 1)} \int_{\mathbb{R}}x\,\frac{d}{dx}(|u|^{p + 1})\ |v|^{p + 1}dx \nonumber \\
\label{406}& \ - \frac{8\beta t^{2}}{(p + 1)}\int_{\mathbb{R}}\frac{d}{dt}(|u|^{p + 1})\ |v|^{p + 1}dx.
\end{align}
Replacing \eqref{405} and \eqref{406} into \eqref{404} it follows that
\begin{align*}
& \frac{d}{dt}\|Ju\|_{L^{2}(\mathbb{R})}^{2} 
=\  \frac{4(2 - p)t}{(p + 1)}\|u\|_{L^{2(p + 1)}(\mathbb{R})}^{2(p + 1)} - \frac{d}{dt}\left[\frac{4t^{2}}{(p + 1)}\|u\|_{L^{2(p + 1)}(\mathbb{R})}^{2(p + 1)}\right] - 4\beta t\int_{\mathbb{R}}|u|^{p + 1}|v|^{p + 1}dx \nonumber \\
\ & - \frac{8\beta t}{(p + 1)} \int_{\mathbb{R}}x\,\frac{d}{dx}(|u|^{p + 1})\ |v|^{p + 1}dx - \frac{8\beta t^{2}}{(p + 1)}\int_{\mathbb{R}}\frac{d}{dt}(|u|^{p + 1})\ |v|^{p + 1}dx \nonumber
\end{align*}
or
\begin{align}
\lefteqn{\frac{d}{dt}\left[\|Ju\|_{L^{2}(\mathbb{R})}^{2} + \frac{4t^{2}}{(p + 1)}\|u\|_{L^{2(p + 1)}(\mathbb{R})}^{2(p + 1)}\right] } \nonumber \\
=\ & \frac{4(2 - p)t}{(p + 1)}\|u\|_{L^{2(p + 1)}(\mathbb{R})}^{2(p + 1)} - 4\beta t\int_{\mathbb{R}}|u|^{p + 1}|v|^{p + 1}dx - \frac{8\beta t}{(p + 1)}\int_{\mathbb{R}}x\ \frac{d}{dx}(|u|^{p + 1})\ |v|^{p + 1}dx \nonumber \\
\label{407}& \ - \frac{8\beta t^{2}}{(p + 1)}\int_{\mathbb{R}}\frac{d}{dt}(|u|^{p + 1})|v|^{p + 1}dx.
\end{align}
Performing similar calculations for the equation \eqref{main2}$_{2}$ we obtain
\begin{align}
\lefteqn{\frac{d}{dt}\left[\|Jv\|_{L^{2}(\mathbb{R})}^{2} + \frac{4t^{2}}{(p + 1)}\|v\|_{L^{2(p + 1)}(\mathbb{R})}^{2(p + 1)}\right] } \nonumber \\
=\ & \frac{4(2 - p)t}{(p + 1)}\|v\|_{L^{2(p + 1)}(\mathbb{R})}^{2(p + 1)} - 4\beta t\int_{\mathbb{R}}|u|^{p + 1}|v|^{p + 1}dx - \frac{8\beta t}{(p + 1)}\int_{\mathbb{R}}x|u|^{p + 1}\ \frac{d}{dx}(|v|^{p + 1})dx \nonumber \\
\label{408}& \ - \frac{8\beta t^{2}}{(p + 1)}\int_{\mathbb{R}}\frac{d}{dt}(|v|^{p + 1})|u|^{p + 1}dx.
\end{align}
Adding \eqref{407} with \eqref{408} it follows that
\begin{align*}
&\frac{d}{dt}\left[\|Ju\|_{L^{2}(\mathbb{R})}^{2} + \|Jv\|_{L^{2}(\mathbb{R})}^{2} + \frac{4t^{2}}{(p + 1)}\|u\|_{L^{2(p + 1)}(\mathbb{R})}^{2(p + 1)} + \frac{4t^{2}}{(p + 1)}\|v\|_{L^{2(p + 1)}(\mathbb{R})}^{2(p + 1)}\right]  \nonumber \\
=\ & \frac{4(2 - p)t}{(p + 1)}\|u\|_{L^{2(p + 1)}(\mathbb{R})}^{2(p + 1)} + \frac{4(2 - p)t}{(p + 1)}\|v\|_{L^{2(p + 1)}(\mathbb{R})}^{2(p + 1)} - 8\beta t\int_{\mathbb{R}}|u|^{p + 1}|v|^{p + 1}dx \nonumber \\
\ & - \frac{8\beta t}{(p + 1)}\int_{\mathbb{R}}x\ \frac{d}{dx}(|u|^{p + 1})\ |v|^{p + 1}dx - \frac{8\beta t}{(p + 1)}\int_{\mathbb{R}}x\ |u|^{p + 1}\ \frac{d}{dx}(|v|^{p + 1})dx \nonumber \\
& \ - \frac{8\beta t^{2}}{(p + 1)}\int_{\mathbb{R}}\frac{d}{dt}(|u|^{p + 1})\ |v|^{p + 1}dx - \frac{8\beta t^{2}}{(p + 1)}\int_{\mathbb{R}}|u|^{p + 1}\ \frac{d}{dt}(|v|^{p + 1})dx.
\end{align*}
Then using integrating by part
\begin{align}
&\frac{d}{dt}\left[\|Ju\|_{L^{2}(\mathbb{R})}^{2} + \|Jv\|_{L^{2}(\mathbb{R})}^{2} + \frac{4t^{2}}{(p + 1)}|u\|_{L^{2(p + 1)}(\mathbb{R})}^{2(p + 1)} + \frac{4t^{2}}{(p + 1)}\|v\|_{L^{2(p + 1)}(\mathbb{R})}^{2(p + 1)}\right]  \nonumber \\
=\ & \frac{4(2 - p)t}{(p + 1)}\|u\|_{L^{2(p + 1)}(\mathbb{R})}^{2(p + 1)}
+ \frac{4(2 - p)t}{(p + 1)}\|v\|_{L^{2(p + 1)}(\mathbb{R})}^{2(p + 1)}- \frac{8\beta pt}{(p + 1)}\int_{\mathbb{R}}|u|^{p + 1}|v|^{p + 1}dx  \nonumber \\
\label{409}\ &  - \frac{8\beta t^{2}}{(p + 1)}\frac{d}{dt}\int_{\mathbb{R}}|u|^{p + 1}|v|^{p + 1}dx.
\end{align}
On the other hand, using the identity
\begin{align*}
\frac{-8\beta t^{2}}{(p + 1)}\frac{d}{dt}\int_{\mathbb{R}}|u|^{p + 1}|v|^{p + 1}dx =  \frac{16\beta t}{(p + 1)}\int_{\mathbb{R}}|u|^{p + 1}|v|^{p + 1}dx - \frac{d}{dt}\left[\frac{8\beta t^{2}}{(p + 1)}\int_{\mathbb{R}}|u|^{p + 1}|v|^{p + 1}dx\right].
\end{align*}
Replacing the above equality into \eqref{409} we obtain
\begin{align*}
&\frac{d}{dt}\left[\|Ju\|_{L^{2}(\mathbb{R})}^{2} + \|Jv\|_{L^{2}(\mathbb{R})}^{2} + \frac{4t^{2}}{(p + 1)}\|u\|_{L^{2(p + 1)}(\mathbb{R})}^{2(p + 1)} + \frac{4t^{2}}{(p + 1)}\|v\|_{L^{2(p + 1)}(\mathbb{R})}^{2(p + 1)}\right]  \nonumber \\
=\ & \frac{4(2 - p)t}{(p + 1)}\|u\|_{L^{2(p + 1)}(\mathbb{R})}^{2(p + 1)} + \frac{4(2 - p)t}{(p + 1)}\|v\|_{L^{2(p + 1)}(\mathbb{R})}^{2(p + 1)} - \frac{8\beta pt}{(p + 1)}\int_{\mathbb{R}}|u|^{p + 1}|v|^{p + 1}dx \nonumber \\
\ & + \frac{16\beta t}{(p + 1)}\int_{\mathbb{R}}|u|^{p + 1}|v|^{p + 1}dx - \frac{d}{dt}\left[\frac{8\beta t^{2}}{(p + 1)}\int_{\mathbb{R}}|u|^{p + 1}|v|^{p + 1}dx\right].
\end{align*}
Thus
\begin{align*}
&\frac{d}{dt}\left[\ \|Ju\|_{L^{2}(\mathbb{R})}^{2} + \|Jv\|_{L^{2}(\mathbb{R})}^{2} \right. \nonumber \\
& \left. \quad + \frac{4t^{2}}{(p + 1)}\|u\|_{L^{2(p + 1)}(\mathbb{R})}^{2(p + 1)} + \frac{4t^{2}}{(p + 1)}\|v\|_{L^{2(p + 1)}(\mathbb{R})}^{2(p + 1)} + \frac{8\beta t^{2}}{(p + 1)}\int_{\mathbb{R}}|u|^{p + 1}|v|^{p + 1}dx\right]  \nonumber \\
=\ & \frac{4(2 - p)t}{(p + 1)}\|u\|_{L^{2(p + 1)}(\mathbb{R})}^{2(p + 1)} + \frac{4(2 - p)t}{(p + 1)}\|v\|_{L^{2(p + 1)}(\mathbb{R})}^{2(p + 1)} + \frac{8\beta (2 - p)t}{(p + 1)}\int_{\mathbb{R}}|u|^{p + 1}|v|^{p + 1}dx.
\end{align*}
But, $2 - p < 1,$ because $p>1.$ Then performing straightforward estimates in the above equation we have
\begin{align}
&\frac{d}{dt}\left[\ \|Ju\|_{L^{2}(\mathbb{R})}^{2} + \|Jv\|_{L^{2}(\mathbb{R})}^{2} \right. \nonumber \\
& \left. \quad + \frac{4t^{2}}{(p + 1)}\|u\|_{L^{2(p + 1)}(\mathbb{R})}^{2(p + 1)} + \frac{4t^{2}}{(p + 1)}\|v\|_{L^{2(p + 1)}(\mathbb{R})}^{2(p + 1)} + \frac{8\beta t^{2}}{(p + 1)}\int_{\mathbb{R}}|u|^{p + 1}|v|^{p + 1}dx\right]  \nonumber \\
\label{410}\leq \ & \frac{4t}{(p + 1)}\|u\|_{L^{2(p + 1)}(\mathbb{R})}^{2(p + 1)} + \frac{4t}{(p + 1)}\|v\|_{L^{2(p + 1)}(\mathbb{R})}^{2(p + 1)} + \frac{8\beta t}{(p + 1)}\int_{\mathbb{R}}|u|^{p + 1}|v|^{p + 1}dx.
\end{align}
Therefore
\begin{align}
\label{411}\frac{d}{dt}\left[\ \|Ju\|_{L^{2}(\mathbb{R})}^{2} + \|Jv\|_{L^{2}(\mathbb{R})}^{2} + t^{2}f(t)\right] \leq tf(t),
\end{align}
where
\begin{align*}
f(t) = \frac{4}{(p + 1)}\left[\|u\|_{L^{2(p + 1)}(\mathbb{R})}^{2(p + 1)} + \|v\|_{L^{2(p + 1)}(\mathbb{R})}^{2(p + 1)} + 2\beta\int_{\mathbb{R}}|u|^{p + 1}|v|^{p + 1}dx\right].
\end{align*}
Integrating \eqref{411} over $t\in [0,\,T],$ with $T$ arbitrary, using the positivity of the two first terms and performing straightforward calculations we have
\begin{align}\label{ecudoc4}
t^{2}f(t)\leq\|Ju\|_{L^{2}(\mathbb{R})}^{2} + \|Jv\|_{L^{2}(\mathbb{R})}^{2} + t^{2}f(t) \leq  \|xu_{0}\|_{L^{2}(\mathbb{R})}^{2} + \|xv_{0}\|_{L^{2}(\mathbb{R})}^{2} + \int_{0}^{t}s\,f(s)ds.
\end{align}
Which we can rewrite as
\begin{align}
\label{412}t^{2}f(t) \leq  \Upsilon + \int_{1}^{t}sf(s)ds,
\end{align}
where
\begin{align*}
\Upsilon = \left[\|xu_{0}\|_{L^{2}(\mathbb{R})}^{2} + \|xv_{0}\|_{L^{2}(\mathbb{R})}^{2} + \int_{0}^{1}s\,f(s)ds\right].
\end{align*}
Using \eqref{306} of the Lemma \ref{lemma302} it is easy to see that
\begin{align*}
\Upsilon \leq \Theta(\beta,\,p,\,\|u_{0}\|_{L^{2(p + 1)}(\mathbb{R})},\,\|v_{0}\|_{L^{2(p + 1)}(\mathbb{R})},\,\|xu_{0}\|_{L^{2}(\mathbb{R})},\,\|xv_{0}\|_{L^{2}(\mathbb{R})},\,\|u_{0x}\|_{L^{2}(\mathbb{R})},\,\|v_{0x}\|_{L^{2}(\mathbb{R})}).
\end{align*}
Thereby, $\Upsilon \leq \Theta$ which is finite by hypothesis. Hence \eqref{412} can be written in the following way
\begin{align*}
F(t) \leq \Theta + \int_{1}^{t}G(s)F(s)ds,
\end{align*}
for $F(t)\equiv t^{2}f(t)$ and
\begin{align*}
G(t) \equiv \frac{1}{t},\quad {\rm for}\quad t \geq 1.
\end{align*}
Using the Gronwall's inequality, we have $F(t)\leq \Theta t,$ $\forall\;t\geq1.$ Using the hypothesis on $u_{0},\,v_{0}$ and integration of \eqref{305} together imply that $f(t)$ is uniformly bounded for all $t,$ in particular for $0\leq t \leq 1.$ This way, there exists a constant
\begin{align*}
C = C(\beta,\,\|u_{0}\|_{L^{2(p + 1)}(\mathbb{R})},\,\|v_{0}\|_{L^{2(p + 1)}(\mathbb{R})},\,\|xu_{0}\|_{L^{2}(\mathbb{R})},\,\|xv_{0}\|_{L^{2}(\mathbb{R})},\,\|u_{0x}\|_{L^{2}(\mathbb{R})},\,\|v_{0x}\|_{L^{2}(\mathbb{R})})
\end{align*}
such that
\begin{align}
\label{413}t^{2}\,f(t) \leq  Ct,\quad {\rm for\ any}\quad t>0
\end{align}
From \eqref{ecudoc4} we have
\begin{align*}
&\|Ju\|_{L^2(\mathbb{R})}^2+\|Jv\|_{L^2(\mathbb{R})}^2+ \frac{4t^{2}}{(p + 1)}\|u\|_{L^{2(p + 1)}(\mathbb{R})}^{2(p + 1)} \\
&+ \frac{4t^{2}}{(p + 1)}\|v\|_{L^{2(p + 1)}(\mathbb{R})}^{2(p + 1)} + \frac{8\beta t^{2}}{(p + 1)}\int_{\mathbb{R}}|u|^{p + 1}|v|^{p + 1}dx
\leq Ct, \ {\rm for\ any}\  t>0.
\end{align*}
Then using that $Ju = e^{ix^{2}/4t}(2it)\partial_{x}(e^{-ix^{2}/4t}u)$ and $Jv = e^{ix^{2}/4t}(2it)\partial_{x}(e^{-ix^{2}/4t}v)$ we obtain
\begin{align*}
&4t^{2}\|\partial_{x}(e^{-ix^{2}/4t}u)\|_{L^{2}(\mathbb{R})}^{2} + 4t^{2}\|\partial_{x}(e^{-ix^{2}/4t}v)\|_{L^{2}(\mathbb{R})}^{2}  \nonumber \\
&   + \frac{4t^{2}}{(p + 1)}\|u\|_{L^{2(p + 1)}(\mathbb{R})}^{2(p + 1)} + \frac{4t^{2}}{(p + 1)}\|v\|_{L^{2(p + 1)}(\mathbb{R})}^{2(p + 1)} + \frac{8\beta t^{2}}{(p + 1)}\int_{\mathbb{R}}|u|^{p + 1}|v|^{p + 1}dx  \leq Ct,
\end{align*}
for  any  $t>0.$ Then
\begin{align*}
&\|\partial_{x}(e^{-ix^{2}/4t}u)\|_{L^{2}(\mathbb{R})}^{2} + \|\partial_{x}(e^{-ix^{2}/4t}v)\|_{L^{2}(\mathbb{R})}^{2}  \nonumber \\
&   + \frac{1}{(p + 1)}\|u\|_{L^{2(p + 1)}(\mathbb{R})}^{2(p + 1)} + \frac{1}{(p + 1)}\|v\|_{L^{2(p + 1)}(\mathbb{R})}^{2(p + 1)} + \frac{2\beta}{(p + 1)}\int_{\mathbb{R}}|u|^{p + 1}|v|^{p + 1}dx  \leq Ct^{-1}.
\end{align*}
Using the Young inequality in the last term on the left hand side and performing Straightforward estimates we have for $\beta < 1$
\begin{align}
&\|\partial_{x}(e^{-ix^{2}/4t}u)\|_{L^{2}(\mathbb{R})}^{2} + \|\partial_{x}(e^{-ix^{2}/4t}v)\|_{L^{2}(\mathbb{R})}^{2}  \nonumber \\
\label{4133}&   + \frac{(1 - \beta)}{(p + 1)}\|u\|_{L^{2(p + 1)}(\mathbb{R})}^{2(p + 1)} + \frac{(1 - \beta)}{(p + 1)}\|v\|_{L^{2(p + 1)}(\mathbb{R})}^{2(p + 1)}  \leq Ct^{-1}.
\end{align}
Hence
\begin{align*}
\|\partial_{x}(e^{-ix^{2}/4t}u)\|_{L^{2}(\mathbb{R})}^{2} + \|\partial_{x}(e^{-ix^{2}/4t}v)\|_{L^{2}(\mathbb{R})}^{2} \leq \frac{C}{t}
\end{align*}
for any $t>0.$ Moreover, using the Gagliardo-Nirenberg inequality, we have
\begin{align*}
\|e^{-ix^{2}/4t}u\|_{L^{\infty}(\mathbb{R})} \leq C\|\partial_{x}(e^{-ix^{2}/4t}u)\|_{L^{2}(\mathbb{R})}^{1/2}\|u\|_{L^{2}(\mathbb{R})}^{1/2}
\end{align*}
and
\begin{align*}
\|e^{-ix^{2}/4t}v\|_{L^{\infty}(\mathbb{R})} \leq C\|\partial_{x}(e^{-ix^{2}/4t}v)\|_{L^{2}(\mathbb{R})}^{1/2}\|v\|_{L^{2}(\mathbb{R})}^{1/2}.
\end{align*}
From the above estimates, we deduce
\begin{align*}
\|e^{-ix^{2}/4t}u\|_{L^{\infty}(\mathbb{R})}^{4} + \|e^{-ix^{2}/4t}v\|_{L^{\infty}(\mathbb{R})}^{4} \leq \frac{C}{t}
\end{align*}
for any $t>0.$ Therefore,
\begin{align*}
\|u\|_{L^{\infty}(\mathbb{R})} \leq \frac{C}{t^{1/4}}\quad {\rm and}\quad \|v\|_{L^{\infty}(\mathbb{R})} \leq \frac{C}{t^{1/4}}.
\end{align*}
Moreover,
\begin{align*}
\frac{(1 - \beta)}{(p + 1)}\|u\|_{L^{2(p + 1)}(\mathbb{R})}^{2(p + 1)} + \frac{(1 - \beta)}{(p + 1)}\|v\|_{L^{2(p + 1)}(\mathbb{R})}^{2(p + 1)}  \leq \frac{C}{t}.
\end{align*}
Then
\begin{align*}
\|u\|_{L^{2(p + 1)}(\mathbb{R})}  \leq \frac{C}{t^{2(p + 1)}}\quad {\rm and}\quad\|v\|_{L^{2(p + 1)}(\mathbb{R})}  \leq \frac{C}{t^{2(p + 1)}}.
\end{align*}
The Theorem \ref{theorem401} follows.

{\bf Remark} Using the Gagliardo-Nirenberg inequality we have
\begin{align*}
\|u\|_{L^{p}(\mathbb{R})} \leq C\|u\|_{L^{2}(\mathbb{R})}^{2/p}\|u\|_{L^{\infty}(\mathbb{R})}^{(p - 2)/p}.
\end{align*}
Using Lemma \ref{lemma301} and \eqref{401} we deduce the following $L^{p}$ estimate
\begin{align*}
\|u\|_{L^{p}(\mathbb{R})} \leq \frac{C}{t^{(p - 2)/4p}} \quad {\rm and}\quad \|v\|_{L^{p}(\mathbb{R})} \leq \frac{C}{t^{(p - 2)/4p}}.
\end{align*}
for $2 < p \leq +\infty.$ \\

\begin{lemma}\label{lema1}
Let $k\in\mathbb{N}$, $w\in L^\infty(\mathbb{R})$ and its derivatives of order $m$, $\partial_x^m w\in L^2(\mathbb{R})$, then
\begin{align}
\|\partial_x^m(|w|^{2k}w)\|_{L^2(\mathbb{R})}\leq C_m\|\partial_x^m w\|_{L^2(\mathbb{R})}\|w\|_{L^\infty(\mathbb{R})}^{2k},
\end{align}
where $C_m$ is a constant that depends of $m$.
\end{lemma}
\proof By induction over $k$. If $k=1$, we have
\begin{align*}
&\|\partial_x^m(|w|^2w)\|_{L^2(\mathbb{R})}=\|\partial_x^m(w\overline{w}w)\|_{L^2(\mathbb{R})}
\leq\sum_{\zeta_1+\zeta_2+\zeta_3=m}\|\partial_x^{\zeta_1}w\cdot\partial_x^{\zeta_2}\overline{w}\cdot\partial_x^{\zeta_3}w\|_{L^2(\mathbb{R})}\\
&\leq\sum_{\zeta_1+\zeta_2+\zeta_3=m}\|\partial_x^{\zeta_1}w\|_{L^{2m/\zeta_1}(\mathbb{R})}\|  \partial_x^{\zeta_2}\overline{w}\|_{L^{2m/\zeta_2}(\mathbb{R})}\|\partial_x^{\zeta_3}w\|_{L^{2m/\zeta_3}(\mathbb{R})}\\
&\leq C_m\sum_{\zeta_1+\zeta_2+\zeta_3=m}\|\partial_x^mw\|_{L^2(\mathbb{R})}^{\zeta_1/m}\|w\|_{L^\infty(\mathbb{R})}^{(m-\zeta_1)/m}
\|\partial_x^mw\|_{L^2(\mathbb{R})}^{\zeta_2/m}\|w\|_{L^\infty(\mathbb{R})}^{(m-\zeta_2)/m}
\|\partial_x^mw\|_{L^2(\mathbb{R})}^{\zeta_3/m}\|w\|_{L^\infty(\mathbb{R})}^{(m-\zeta_3)/m}\\
&\leq C_m\|\partial_x^mw\|_{L^2(\mathbb{R})}\|w\|_{L^\infty(\mathbb{R})}^2.
\end{align*}
Suppose it is valid for $k$, show that it is valid for $k+1$.
\begin{align*}
&\|\partial_x^m(|w|^{2(k+1)}w)\|_{L^2(\mathbb{R})}=\|\partial_x^m((|w|^{2k}w)\overline{w}w)\|_{L^2(\mathbb{R})}
\leq\sum_{\zeta_1+\zeta_2+\zeta_3=m}\|\partial_x^{\zeta_1}(|w|^{2k}w)\cdot\partial_x^{\zeta_2}\overline{w}\cdot\partial_x^{\zeta_3}w\|_{L^2(\mathbb{R})}\\
&\leq\sum_{\zeta_1+\zeta_2+\zeta_3=m}\|\partial_x^{\zeta_1}(|w|^{2k}w)\|_{L^{2m/\zeta_1}(\mathbb{R})}\|  \partial_x^{\zeta_2}\overline{w}\|_{L^{2m/\zeta_2}(\mathbb{R})}\|\partial_x^{\zeta_3}w\|_{L^{2m/\zeta_3}(\mathbb{R})}\\
&\leq C_m\sum_{\zeta_1+\zeta_2+\zeta_3=m}\|\partial_x^m(|w|^{2k}w)\|_{L^2(\mathbb{R})}^{\zeta_1/m}\||w|^{2k}w\|_{L^\infty(\mathbb{R})}^{(m-\zeta_1)/m}
\|\partial_x^mw\|_{L^2(\mathbb{R})}^{\zeta_2/m}\|w\|_{L^\infty(\mathbb{R})}^{(m-\zeta_2)/m}\\
&\hspace{3cm}\times
\|\partial_x^mw\|_{L^2(\mathbb{R})}^{\zeta_3/m}\|w\|_{L^\infty(\mathbb{R})}^{(m-\zeta_3)/m},
\end{align*}
by induction hypothesis, we get
\begin{align*}
\|\partial_x^m(|w|^{2(k+1)}w)\|_{L^2(\mathbb{R})}
&\leq C_m\sum_{\zeta_1+\zeta_2+\zeta_3=m}\|\partial_x^mw\|_{L^2(\mathbb{R})}^{\zeta_1/m}\|w\|_{L^\infty(\mathbb{R})}^{2k\zeta_1/m}
\|w\|_{L^\infty(\mathbb{R})}^{2k(m-\zeta_1)/m}
\|w\|_{L^\infty(\mathbb{R})}^{(m-\zeta_1)/m}\\
&\hspace{3cm}\times\|\partial_x^mw\|_{L^2(\mathbb{R})}^{\zeta_2/m}\|w\|_{L^\infty(\mathbb{R})}^{(m-\zeta_2)/m}
\|\partial_x^mw\|_{L^2(\mathbb{R})}^{\zeta_3/m}\|w\|_{L^\infty(\mathbb{R})}^{(m-\zeta_3)/m}\\
&\leq C_m\|\partial_x^mw\|_{L^2(\mathbb{R})}\|w\|_{L^\infty(\mathbb{R})}^{2(k+1)}.
\end{align*}
The Lemma follows.
\begin{lemma}\label{lema2}
Let $p>1$ an odd integer number, $w,\ z\in L^\infty(\mathbb{R})$ and its derivatives of order $m$, $\partial_x^m w,\ \partial_x^m z\in L^2(\mathbb{R})$, then
\begin{align}
\|\partial_x^m(|w|^{p-1}|z|^{p+1}w)\|_{L^2(\mathbb{R})}\leq C_m(\|\partial_x^m w\|_{L^2(\mathbb{R})}+\|\partial_x^m z\|_{L^2(\mathbb{R})})
\end{align}
where $C_m$ is a constant that depends of $m$.
\end{lemma}
\proof Denote by $\bigoplus_{i=1}^n \zeta_i$ the sum $\zeta_1+...+\zeta_n$, then we obtain
\begin{align*}
&\|\partial_x^m(|w|^{p-1}|z|^{p+1}w)\|_{L^2(\mathbb{R})}
=\|\partial_x^m(\underbrace{(w\overline{w})\cdots(w\overline{w})}_{(p-1)-times}w\underbrace{(z\overline{z})\cdots(z\overline{z})}_{(p+1)-times} )\|_{L^2(\mathbb{R})}\|\\
&=\|\partial_x^m(\underbrace{w\cdots w}_{p-times}
\underbrace{\overline{w}\cdots\overline{w}}_{(p-1)-times}
\underbrace{z\cdots z}_{(p+1)-times}\underbrace{\overline{z}\cdots\overline{z}}_{(p+1)-times}\|_{L^2(\mathbb{R})}\\
&=\|\sum_{\bigoplus_{i=1}^{4p+1}\zeta_i=m}\partial_x^{\zeta_1}w\cdots\partial_x^{\zeta_p}w
\cdot\partial_x^{\zeta_{p+1}}\overline{w}\cdots\partial_x^{\zeta_{p+(p-1)}}\overline{w}
\cdot\partial_x^{\zeta_{p+(p-1)+1}}z\cdots\partial_x^{\zeta_{p+(p-1)+(p+1)}}z\\
&\hspace{3cm}\times
\partial_x^{\zeta_{p+(p-1)+(p+1)+1}}\overline{z}\cdots\partial_x^{\zeta_{p+(p-1)+(p+1)+(p+1)}}\overline{z}
\|_{L^2(\mathbb{R})}
\end{align*}
\begin{align*}
&\leq\sum_{\bigoplus_{i=1}^{4p+1}\zeta_i=m}\|\partial_x^{\zeta_1}w\cdots\partial_x^{\zeta_p}w
\cdot\partial_x^{\zeta_{p+1}}\overline{w}\cdots\partial_x^{\zeta_{2p-1}}\overline{w}
\cdot\partial_x^{\zeta_{2p}}z\cdots\partial_x^{\zeta_{3p}}z
\partial_x^{\zeta_{3p+1}}\overline{z}\cdots\partial_x^{\zeta_{4p+1}}\overline{z}
\|_{L^2(\mathbb{R})}\\
&\leq \sum_{\bigoplus_{i=1}^{4p+1}\zeta_i=m}\|\partial_x^{\zeta_1}w\|_{L^{2m/\zeta_1}(\mathbb{R})}
\cdots\|\partial_x^{\zeta_p}w\|_{L^{2m/\zeta_p}(\mathbb{R})}\\
&\hspace{1cm}\times\|\partial_x^{\zeta_{p+1}}\overline{w}\|_{L^{2m/\zeta_{p+1}}(\mathbb{R})}\cdots\|\partial_x^{\zeta_{2p-1}}\overline{w}\|_{L^{2m/\zeta_{2p-1}}(\mathbb{R})}
\cdot\|\partial_x^{\zeta_{2p}}z\|_{L^{2m/\zeta_{2p}}(\mathbb{R})}\cdots\|\partial_x^{\zeta_{3p}}z\|_{L^{2m/\zeta_{3p}}(\mathbb{R})}\\
&\hspace{1cm}\times
\|\partial_x^{\zeta_{3p+1}}\overline{z}\|_{L^{2m/\zeta_{3p+1}}(\mathbb{R})}\cdots\|\partial_x^{\zeta_{4p+1}}\overline{z}\|_{L^{2m/\zeta_{4p+1}}(\mathbb{R})}.
\end{align*}
Note that in the inequality last we used the H\"{o}lder inequality, by the Gagliardo-Nirenberg inequality, we have
\begin{align*}
\|\partial_x^m(|w|^{p-1}|z|^{p+1}w)\|_{L^2(\mathbb{R})}
&\leq C\sum_{\bigoplus_{i=1}^{4p+1}\zeta_i=m}\|\partial_x^mw\|_{L^2(\mathbb{R})}^{\zeta_1/m}\|w\|_{L^\infty(\mathbb{R})}^{(m-\zeta_1)/m}
\cdots\|\partial_x^mw\|_{L^2(\mathbb{R})}^{\zeta_p/m}\|w\|_{L^\infty(\mathbb{R})}^{(m-\zeta_p)/m}\\
&\hspace{0.5cm}\times\|\partial_x^mw\|_{L^2(\mathbb{R})}^{\zeta_{p+1}/m}\|w\|_{L^\infty(\mathbb{R})}^{(m-\zeta_{p+1})/m}
\cdots\|\partial_x^mw\|_{L^2(\mathbb{R})}^{\zeta_{2p-1}/m}\|w\|_{L^\infty(\mathbb{R})}^{(m-\zeta_{2p-1})/m}\\
&\hspace{0.5cm}\times\|\partial_x^mz\|_{L^2(\mathbb{R})}^{\zeta_{2p}/m}\|w\|_{L^\infty(\mathbb{R})}^{(m-\zeta_{2p})/m}
\cdots\|\partial_x^mw\|_{L^2(\mathbb{R})}^{\zeta_{3p}/m}\|w\|_{L^\infty(\mathbb{R})}^{(m-\zeta_{3p})/m}\\
&\hspace{0.5cm}\times\|\partial_x^mz\|_{L^2(\mathbb{R})}^{\zeta_{3p+1}/m}\|w\|_{L^\infty(\mathbb{R})}^{(m-\zeta_{3p+1})/m}
\cdots\|\partial_x^mw\|_{L^2(\mathbb{R})}^{\zeta_{4p+1}/m}\|w\|_{L^\infty(\mathbb{R})}^{(m-\zeta_{4p+1})/m}\\
&=C\sum_{\bigoplus_{i=1}^{4p+1}\zeta_i=m}\|\partial_x^mw\|_{L^2(\mathbb{R})}^{(\zeta_1+...+\zeta_{2p-1})/m}
\|w\|_{L^\infty(\mathbb{R})}^{(m-\zeta_1)/m+...+(m-\zeta_{2p-1})/m}\\
&\hspace{0.5cm}\times
\|\partial_x^mz\|_{L^2(\mathbb{R})}^{(\zeta_{2p}+\ldots+\zeta_{4p-1})/m}
\|z\|_{L^\infty(\mathbb{R})}^{(m-\zeta_{2p})/m+\ldots+(m-\zeta_{4p-1})/m}.
\end{align*}
By the Young inequality, we have
\begin{align*}
\|\partial_x^m(|w|^{p-1}|z|^{p+1}w)\|_{L^2(\mathbb{R})}
&\leq C\sum_{\bigoplus_{i=1}^{4p+1}\zeta_i=m}\frac{\zeta_1+\ldots+\zeta_{2p-1}}{m}
\|\partial_x^mw\|_{L^2(\mathbb{R})}\\
&\hspace{2.5cm}+\frac{\zeta_{2p}+\ldots+\zeta_{4p-1}}{m}\|\partial_x^mz\|_{L^2(\mathbb{R})}\\
&\leq C_m(\|\partial_x^mw\|_{L^2(\mathbb{R})}+\|\partial_x^mz\|_{L^2(\mathbb{R})}).
\end{align*}
The Lemma follows.
\begin{theorem}
\label{theorem402} Let $(u_{0},\,v_{0})\in H^{1}(\mathbb{R})\times H^{1}(\mathbb{R})$ with $(x^{n}u_{0},\,x^{n}v_{0})\in L^{2}(\mathbb{R})\times L^{2}(\mathbb{R})$ and  $p > 1$ odd integer number. Then there exists a positive constant $C_{m}$ depending on $\|u_{0}\|_{H^{1}(\mathbb{R})},$ $\|v_{0}\|_{H^{1}(\mathbb{R})}$ and $\|x^{n}u_{0}\|_{L^{2}(\mathbb{R})},$ $\|x^{n}v_{0}\|_{L^{2}(\mathbb{R})}$ but independent of $t$ such that
\begin{align}
\label{414}\|J^{m}u\|_{L^{2}(\mathbb{R})} \leq C_{m}e^{t}\quad {\rm and}\quad \|J^{m}v\|_{L^{2}(\mathbb{R})} \leq C_{m}e^{t},
\end{align}
for $m=1,2,...,n$.
\end{theorem}
\proof From \eqref{main3}$_{1}$ we have
\begin{align*}
i(J^{m}u)_{t} + (J^{m}u)_{xx} = J^{m}(|u|^{2p}u) + \beta J^{m}(|u|^{p - 1}|v|^{p + 1}u).
\end{align*}
Multiplying the above equation by $(\overline{J^{m}u})$ we have
\begin{align}
\label{415}i(\overline{J^{m}u})(J^{m}u)_{t} + (\overline{J^{m}u})(J^{m}u)_{xx} = J^{m}(|u|^{2p}u)(\overline{J^{m}u}) + \beta J^{m}(|u|^{p - 1}|v|^{p + 1}u)(\overline{J^{m}u}).
\end{align}
Applying conjugate in \eqref{415} we obtain
\begin{align}
\label{416}-i(J^{m}u)(\overline{J^{m}u})_{t} + (J^{m}u)(\overline{J^{m}u})_{xx} = \overline{J^{m}(|u|^{2p}u)(\overline{J^{m}u})} + \beta \overline{J^{m}(|u|^{p - 1}|v|^{p + 1}u)(\overline{J^{m}u})}.
\end{align}
Subtracting \eqref{415} with \eqref{416}, integrating over $x\in\mathbb{R}$ and performing straightforward calculations we have
\begin{align}
\label{417}\frac{d}{dt}\|J^{m}u\|_{L^{2}(\mathbb{R})}^{2} = 2Im\int_{\mathbb{R}}J^{m}(|u|^{2p}u)(\overline{J^{m}u})dx + 2\beta Im\int_{\mathbb{R}}J^{m}(|u|^{p - 1}|v|^{p + 1}u)(\overline{J^{m}u})dx.
\end{align}
In a similar way, using \eqref{main3}$_{2}$ we have
\begin{align}
\label{418}\frac{d}{dt}\|J^{m}v\|_{L^{2}(\mathbb{R})}^{2} = 2Im\int_{\mathbb{R}}J^{m}(|v|^{2p}v)(\overline{J^{m}v})dx + 2\beta Im\int_{\mathbb{R}}J^{m}(|v|^{p - 1}|u|^{p + 1}v)(\overline{J^{m}v})dx.
\end{align}
Adding \eqref{417} and \eqref{418} it follows that
\begin{align}
\lefteqn{\frac{d}{dt}\left[\|J^{m}u\|_{L^{2}(\mathbb{R})}^{2} + \|J^{m}v\|_{L^{2}(\mathbb{R})}^{2}\right] } \nonumber \\
=\ & 2Im\int_{\mathbb{R}}J^{m}(|u|^{2p}u)(\overline{J^{m}u})dx + 2Im\int_{\mathbb{R}}J^{m}(|v|^{2p}v)(\overline{J^{m}v})dx \nonumber \\
\label{419}&\ + 2\beta Im\int_{\mathbb{R}}J^{m}(|u|^{p - 1}|v|^{p + 1}u)(\overline{J^{m}u})dx + 2\beta Im\int_{\mathbb{R}}J^{m}(|v|^{p - 1}|u|^{p + 1}v)(\overline{J^{m}v})dx.
\end{align}
We estimate the first two terms in the right hand side. Indeed
\begin{align*}
J^{m}(|u|^{2p}u) = \ & e^{ix^{2}/4t}(2it)^{m}\partial_{x}^{m}(e^{-ix^{2}/4t}|u|^{2p}u) \\
 = \ & e^{ix^{2}/4t}(2it)^{m}\partial_{x}^{m}(|e^{-ix^{2}/4t}u|^{2p}e^{-ix^{2}/4t}u) \\
  = \ & e^{ix^{2}/4t}(2it)^{m}\partial_{x}^{m}(|w|^{2p}w)
\end{align*}
where $w = e^{-ix^{2}/4t}u.$ By Lemma \ref{lema1}, we have
\begin{align}
\label{423}\|J^{m}(|u|^{2p}\ u)\|_{L^{2}(\mathbb{R})} \leq C_m(2t)^{m}\left\|\partial_{x}^{m}w\right\|_{L^{2}(\mathbb{R})}\|w\|_{L^{\infty}(\mathbb{R})}^{2p}
= C_m\left\|J^{m}u\right\|_{L^{2}(\mathbb{R})}\|u\|_{L^{\infty}(\mathbb{R})}^{2p}.
\end{align}
In a similar way we obtain
\begin{align}
\label{424}\|J^{m}(|v|^{2p}\ v)\|_{L^{2}(\mathbb{R})} \leq C_m\left\|J^{m}v\right\|_{L^{2}(\mathbb{R})}\ \|v\|_{L^{\infty}(\mathbb{R})}^{2p}.
\end{align}
Using H\"{o}lder's  inequality in the first term of the right hand side in \eqref{419} we have
\begin{align}
&2Im\int_{\mathbb{R}}J^{m}(|u|^{2p}\ u)(\overline{J^{m}u})dx \leq \  2\|J^{m}(|u|^{2p}\ u)\|_{L^{2}(\mathbb{R})}\ \|J^{m}u\|_{L^{2}(\mathbb{R})} \nonumber  \\
\label{425}\leq \ & C_m\left\|J^{m}u\right\|_{L^{2}(\mathbb{R})}\ \|u\|_{L^{\infty}(\mathbb{R})}^{2p}\ \|J^{m}u\|_{L^{2}(\mathbb{R})} \leq  C_m\left\|J^{m}u\right\|_{L^{2}(\mathbb{R})}^{2}\ \|u\|_{L^{\infty}(\mathbb{R})}^{2p}.
\end{align}
In a similar way
\begin{align}
\label{426}2Im\int_{\mathbb{R}}J^{m}(|v|^{2p}\ v)(\overline{J^{m}v})dx \leq  C_m\left\|J^{m}v\right\|_{L^{2}(\mathbb{R})}^{2}\ \|v\|_{L^{\infty}(\mathbb{R})}^{2p}.
\end{align}
On the other hand, now we estimate the third term in \eqref{419}. Using the H\"{o}lder inequality we have
\begin{align}
\label{427}2\beta Im\int_{\mathbb{R}}J^{m}(|u|^{p - 1}|v|^{p + 1}u)(\overline{J^{m}u})dx \leq 2\beta\|J^{m}(|u|^{p - 1}|v|^{p + 1}u)\|_{L^{2}(\mathbb{R})}\ \|J^{m}u\|_{L^{2}(\mathbb{R})}.
\end{align}
We estimate the $\|J^{m}(|u|^{p - 1}|v|^{p + 1}u)\|_{L^{2}(\mathbb{R})}$ term. From the definition of $J^{m}$ we obtain
\begin{align}
J^{m}(|u|^{p - 1}|v|^{p + 1}u) = \ & e^{ix^{2}/4t}(2it)^{m}\partial_{x}^{m}(e^{-ix^{2}/4t}|u|^{p - 1}|v|^{p + 1}u) \nonumber \\
 = \ & e^{ix^{2}/4t}(2it)^{m}\partial_{x}^{m}(|e^{-ix^{2}/4t}u|^{p - 1}|e^{-ix^{2}/4t}v|^{p + 1}e^{-ix^{2}/4t}u) \nonumber  \\
\label{428}  = \ & e^{ix^{2}/4t}(2it)^{m}\partial_{x}^{m}(|w|^{p - 1}|z|^{p + 1}w)
\end{align}
where $w = e^{-ix^{2}/4t}u$ and $z = e^{-ix^{2}/4t}v.$\\
\\
{\bf Remark.} We observe that the power  nonlinearity $|u|^{p - 1}u$ is not smooth if $p$ is not an odd integer. \\
\\
By the Lemma \ref{lema2}, we get
\begin{align*}
\|J^{m}(|u|^{p-1}\ |v|^{p+1}\ u)\|_{L^{2}(\mathbb{R})}&\leq C_m(2t)^m\|\partial_x^m(|w|^{p-1}\ |z|^{p+1}\ w)\|_{L^{2}(\mathbb{R})}\\
&\leq C_m(2t)^m(\|\partial_x^mw\|_{L^{2}(\mathbb{R})}+\|\partial_x^mz\|_{L^{2}(\mathbb{R})})\\
&= C_m(\|J^mu\|_{L^{2}(\mathbb{R})}+\|J^mv\|_{L^{2}(\mathbb{R})}),
\end{align*}
then by the Young inequality, we have
\begin{align}\label{ecu1}
\|J^{m}(|u|^{p-1}\ |v|^{p+1}\ u)\|_{L^{2}(\mathbb{R})}\|J^{m}u\|_{L^{2}(\mathbb{R})}&\leq
C_m(\|J^mu\|_{L^{2}(\mathbb{R})}^2+\|J^mv\|_{L^{2}(\mathbb{R})}\|J^mu\|_{L^{2}(\mathbb{R})})\nonumber\\
&\leq C_m(\|J^mu\|_{L^{2}(\mathbb{R})}^2+\|J^mv\|_{L^{2}(\mathbb{R})}^2).
\end{align}
Now of \eqref{ecu1} and \eqref{417} we have
\begin{align}\label{ecu2}
2\beta Im[\int_\mathbb{R}J^m(|u|^{p-1}|v|^{p+1}u)\overline{J^mu}dx]
&\leq C_m(\|J^mu\|_{L^{2}(\mathbb{R})}^2+\|J^mv\|_{L^{2}(\mathbb{R})}^2).
\end{align}
In a similar way
\begin{align}\label{ecu3}
2\beta Im[\int_\mathbb{R}J^m(|v|^{p-1}|u|^{p+1}v)\overline{J^mu}dx]
&\leq C_m(\|J^mv\|_{L^{2}(\mathbb{R})}^2+\|J^mu\|_{L^{2}(\mathbb{R})}^2).
\end{align}
Replacing \eqref{425}, \eqref{426} together with \eqref{ecu2}, \eqref{ecu3} over \eqref{419}
\begin{align}
\frac{d}{dt}\left[\|J^{m}u\|_{L^{2}(\mathbb{R})}^{2} + \|J^{m}v\|_{L^{2}(\mathbb{R})}^{2}\right] \leq C_{m}\left[\|J^{m}u\|_{L^{2}(\mathbb{R})}^{2} + \|J^{m}v\|_{L^{2}(\mathbb{R})}^{2}\right].
\end{align}
Integrating over $t\in [0,\,T]$ with $T$ arbitrary and using Gronwall's inequality it follows that
\begin{align}
\|J^{m}u\|_{L^{2}(\mathbb{R})}^{2} + \|J^{m}v\|_{L^{2}(\mathbb{R})}^{2} \leq & \ \left[\|x^{m}u_{0}\|_{L^{2}(\mathbb{R})}^{2} + \|x^{m}v_{0}\|_{L^{2}(\mathbb{R})}^{2}\right]e^{C_{m}t} \nonumber.
\end{align}
The result follows. \\

\renewcommand{\theequation}{\thesection.\arabic{equation}}
\setcounter{equation}{0}

\section{The Main Result}

In this section we state and prove our theorem, which states that all solutions of finite energy to \eqref{main2} are smooth for $t \neq 0$ provided that the initial functions in $H^{1}(\mathbb{R})$ decay rapidly enough as $|x| \rightarrow +\infty.$

\begin{lemma}
\label{lemma501}
Let $p>1$ odd integer number, we have $(u^{k},\,v^{k})$ is a Cauchy sequence in $C([0,\,T]:\,H^{1}(\mathbb{R})) \times C([0, T]:\,H^{1}(\mathbb{R}))$ for any $T > 0.$ Moreover
\begin{align*}
\|u^{k} - u^{j}\|_{H^{1}(\mathbb{R})}^{2} + \|v^{k} - v^{j}\|_{H^{1}(\mathbb{R})}^{2} \leq C(T)\left[\|u_{0}^{k} - u_{0}^{j}\|_{H^{1}(\mathbb{R})}^{2} + \|v_{0}^{k} - v_{0}^{j}\|_{H^{1}(\mathbb{R})}^{2}\right]
\end{align*}
where $C(T)$ is a positive constant independent of $k$ and $j.$
\end{lemma}
\proof Let $(u^{k},\,v^{k})$ be the solution of \eqref{main2}, then
\begin{eqnarray}
\label{main11} \left\lbrace
\begin{array}{l}
iu_{t}^{k} + u_{xx}^{k} = |u^{k}|^{2p}u^{k} + \beta|u^{k}|^{p - 1}|v^{k}|^{p + 1}u^{k},    \\
iv_{t}^{k} + v_{xx}^{k} = |v^{k}|^{2p}v^{k} + \beta|v^{k}|^{p - 1}|u^{k}|^{p + 1}v^{k},
\end{array}
\right.
\end{eqnarray}
and
\begin{eqnarray}
\label{main22} \left\lbrace
\begin{array}{l}
iu_{t}^{j} + u_{xx}^{j} = |u^{j}|^{2p}u^{j} + \beta|u^{j}|^{p - 1}|v^{j}|^{p + 1}u^{j},    \\
iv_{t}^{j} + v_{xx}^{j} = |v^{j}|^{2p}v^{j} + \beta|v^{j}|^{p - 1}|u^{j}|^{p + 1}v^{j}.
\end{array}
\right.
\end{eqnarray}
Subtracting \eqref{main11}$_{1}$ with \eqref{main22}$_{1}$ we have
\begin{eqnarray}
\label{503}
\lefteqn{i(u^{k} - u^{j})_{t} + (u^{k} - u^{j})_{xx} } \nonumber \\
& = & |u^{k}|^{2p}u^{k} -  |u^{j}|^{2p}u^{j} + \beta(|u^{k}|^{p - 1}|v^{k}|^{p + 1}u^{k} - |u^{j}|^{p - 1}|v^{j}|^{p + 1}u^{j}).
\end{eqnarray}
Multiplying \eqref{503} by $(\overline{u^{k} - u^{j}})$ we obtain the following equation and the conjugate respectively:
\begin{eqnarray}
\lefteqn{i(\overline{u^{k} - u^{j}})(u^{k} - u^{j})_{t} + (\overline{u^{k} - u^{j}})(u^{k} - u^{j})_{xx} } \nonumber \\
& = & (|u^{k}|^{2p}u^{k} -  |u^{j}|^{2p}u^{j})(\overline{u^{k} - u^{j}}) \nonumber \\
\label{504}& & + \beta(|u^{k}|^{p - 1}|v^{k}|^{p + 1}u^{k} - |u^{j}|^{p - 1}|v^{j}|^{p + 1}u^{j})(\overline{u^{k} - u^{j}}).
\end{eqnarray}
\begin{eqnarray}
\lefteqn{-i(u^{k} - u^{j})(\overline{u^{k} - u^{j}})_{t} + (u^{k} - u^{j})(\overline{u^{k} - u^{j}})_{xx} } \nonumber \\
& = & \overline{(|u^{k}|^{2p}u^{k} -  |u^{j}|^{2p}u^{j})(\overline{u^{k} - u^{j}})} \nonumber \\
\label{505}& & + \beta\,\overline{(|u^{k}|^{p - 1}|v^{k}|^{p + 1}u^{k} - |u^{j}|^{p - 1}|v^{j}|^{p + 1}u^{j})(\overline{u^{k} - u^{j}})}.
\end{eqnarray}
Subtracting and integrating over $x\in\mathbb{R}$ we have
\begin{eqnarray}
\lefteqn{\frac{d}{dt}\|u^{k} - u^{j}\|_{L^{2}(\mathbb{R})}^{2} } \nonumber \\
& = & 2Im\int_{\mathbb{R}}(|u^{k}|^{2p}u^{k} -  |u^{j}|^{2p}u^{j})(\overline{u^{k} - u^{j}})dx \nonumber \\
\label{506}& &+ 2\beta Im\int_{\mathbb{R}}(|u^{k}|^{p - 1}|v^{k}|^{p + 1}u^{k} - |u^{j}|^{p - 1}|v^{j}|^{p + 1}u^{j})(\overline{u^{k} - u^{j}})dx.
\end{eqnarray}
In a similar way and performing straightforward calculations we obtain that
\begin{eqnarray}
\lefteqn{\frac{d}{dt}\|v^{k} - v^{j}\|_{L^{2}(\mathbb{R})}^{2} } \nonumber \\
& = & 2Im\int_{\mathbb{R}}(|v^{k}|^{2p}v^{k} -  |v^{j}|^{2p}v^{j})(\overline{v^{k} - v^{j}})dx \nonumber \\
\label{507}& &+ 2\beta Im\int_{\mathbb{R}}(|v^{k}|^{p - 1}|u^{k}|^{p + 1}v^{k} - |v^{j}|^{p - 1}|u^{j}|^{p + 1}v^{j})(\overline{v^{k} - v^{j}})dx.
\end{eqnarray}
We estimate the first term on the right hand side in \eqref{506}
\begin{eqnarray*}
\lefteqn{(|u^{k}|^{2p}u^{k} -  |u^{j}|^{2p}u^{j})(\overline{u^{k} - u^{j}}) } \nonumber \\
& = &  (|u^{k}|^{p} - |u^{j}|^{p})(|u^{k}|^{p} + |u^{j}|^{p})u^{k}(\overline{u^{k} - u^{j}}) + |u^{j}|^{2p}\left|u^{k} - u^{j}\right|^{2}. \nonumber \\
\end{eqnarray*}
Then
\begin{eqnarray}
\lefteqn{2Im\int_{\mathbb{R}}(|u^{k}|^{2p}u^{k} -  |u^{j}|^{2p}u^{j})(\overline{u^{k} - u^{j}})dx } \nonumber \\
& = &  2Im\int_{\mathbb{R}}(|u^{k}|^{p} - |u^{j}|^{p})(|u^{k}|^{p} + |u^{j}|^{p})u^{k}(\overline{u^{k} - u^{j}})dx \nonumber \\
& \leq &  2\|u^{k}\|_{L^{\infty}(\mathbb{R})}\left(\|u^{k}\|_{L^{\infty}(\mathbb{R})}^{p} + \|u^{j}\|_{L^{\infty}(\mathbb{R})}^{p}\right)\int_{\mathbb{R}}|\,|u^{k}|^{p} - |u^{j}|^{p}\,|\ |u^{k} - u^{j}|dx \nonumber \\
\label{508}& \leq & M_{1}\int_{\mathbb{R}}|\,|u^{k}|^{p} - |u^{j}|^{p}\,|\ |u^{k} - u^{j}|dx
\end{eqnarray}
where $M_{1} = 2\|u^{k}\|_{L^{\infty}(\mathbb{R})}\left(\|u^{k}\|_{L^{\infty}(\mathbb{R})}^{p} + \|u^{j}\|_{L^{\infty}(\mathbb{R})}^{p}\right).$ Moreover,
\begin{eqnarray*}
\lefteqn{|\,|u^{k}|^{p} - |u^{j}|^{p}\,| } \\
& = & |\,|u^{k}| - |u^{j}|\,|\ |\,|u^{k}|^{p - 1} + |u^{k}|^{p - 2}|u^{j}| + |u^{k}|^{p - 3}|u^{j}|^{2} + \cdots + |u^{j}|^{p - 1}\,|.
\end{eqnarray*}
Using $\left| |a| - |b| \right| \leq |a - b|$ we obtain
\begin{eqnarray*}
\lefteqn{|\,|u^{k}|^{p} - |u^{j}|^{p}\,|\ |u^{k} - u^{j}| } \\
& \leq & |u^{k} - u^{j}|^{2}\ |\,|u^{k}|^{p - 1} + |u^{k}|^{p - 2}|u^{j}| + |u^{k}|^{p - 3}|u^{j}|^{2} + \cdots + |u^{j}|^{p - 1}\,|.
\end{eqnarray*}
Then
\begin{eqnarray*}
\int_{\mathbb{R}}|\,|u^{k}|^{p} - |u^{j}|^{p}\,|\ |u^{k} - u^{j}|dx
 \leq M_{2}\int_{\mathbb{R}}|u^{k} - u^{j}|^{2}dx
\end{eqnarray*}
where $M_{2} = \left(\|u^{k}\|_{L^{\infty}(\mathbb{R})}^{p - 1} + \|u^{k}\|_{L^{\infty}(\mathbb{R})}^{p - 2}\|u^{j}\|_{L^{\infty}(\mathbb{R})} + \|u^{k}\|_{L^{\infty}(\mathbb{R})}^{p - 3}\|u^{j}\|_{L^{\infty}(\mathbb{R})}^{2} + \cdots + \|u^{j}\|_{L^{\infty}(\mathbb{R})}^{p - 1}\right).$ Then replacing into \eqref{508} we obtain
\begin{eqnarray}
\label{509}{2Im\int_{\mathbb{R}}(|u^{k}|^{2p}u^{k} -  |u^{j}|^{2p}u^{j})(\overline{u^{k} - u^{j}})dx } \leq M_{3}\|u^{k} - u^{j}\|_{L^{2}(\mathbb{R})}^{2},
\end{eqnarray}
where $M_{3} = M_{1}\cdot M_{2}.$ In a similar way we estimate the first term in \eqref{507}. That is,
\begin{eqnarray}
\label{510} {2Im\int_{\mathbb{R}}(|v^{k}|^{2p}v^{k} -  |v^{j}|^{2p}v^{j})(\overline{v^{k} - v^{j}})dx } \leq \widetilde{M}_{3}\|v^{k} - v^{j}\|_{L^{2}(\mathbb{R})}^{2},
\end{eqnarray}
where $\widetilde{M}_{3} = \widetilde{M}_{1}\cdot \widetilde{M}_{2}$ with $\widetilde{M}_{1} = 2\|v^{k}\|_{L^{\infty}(\mathbb{R})}\left(\|v^{k}\|_{L^{\infty}(\mathbb{R})}^{p} + \|v^{j}\|_{L^{\infty}(\mathbb{R})}^{p}\right)$ and \\
$\widetilde{M}_{2} = \left(\|v^{k}\|_{L^{\infty}(\mathbb{R})}^{p - 1} + \|v^{k}\|_{L^{\infty}(\mathbb{R})}^{p - 2}\|v^{j}| + \|v^{k}\|_{L^{\infty}(\mathbb{R})}^{p - 3}\|v^{j}\|_{L^{\infty}(\mathbb{R})}^{2} + \cdots + \|v^{j}\|_{L^{\infty}(\mathbb{R})}^{p - 1}\right).$ \\
\\
Now we estimate the second term in \eqref{506}. In fact,
\begin{align}
& 2\beta Im\int_{\mathbb{R}}(|u^{k}|^{p - 1}|v^{k}|^{p + 1}u^{k} - |u^{j}|^{p - 1}|v^{j}|^{p + 1}u^{j})(\overline{u^{k} - u^{j}})dx \nonumber \\
=\ & 2\beta Im\int_{\mathbb{R}}\left[|u^{k}|^{p - 1}|v^{k}|^{p + 1}(u^{k} - u^{j}) + (|u^{k}|^{p - 1}|v^{k}|^{p + 1} - |u^{j}|^{p - 1}|v^{j}|^{p + 1})u^{j}\right](\overline{u^{k} - u^{j}})dx \nonumber \\
\label{511}=\ & 2\beta Im\int_{\mathbb{R}}\left[(|u^{k}|^{p - 1}|v^{k}|^{p + 1} - |u^{j}|^{p - 1}|v^{j}|^{p + 1})u^{j}\right](\overline{u^{k} - u^{j}})dx.
\end{align}
Similarly we estimate the second term in \eqref{507}. That is,
\begin{align}
\label{512}2&\beta Im\int_{\mathbb{R}}\left[|v^{k}|^{p - 1}|u^{k}|^{p + 1}v^{k} - |v^{j}|^{p - 1}|u^{j}|^{p + 1}v^{j}\right](\overline{v^{k} - v^{j}})dx\nonumber\\
&=2\beta Im\int_{\mathbb{R}}\left[(|v^{k}|^{p - 1}|u^{k}|^{p + 1} - |v^{j}|^{p - 1}|u^{j}|^{p + 1})v^{j}\right](\overline{v^{k} - v^{j}})dx.
\end{align}
In \eqref{511} we have
\begin{align}
&2\beta Im\int_{\mathbb{R}}\left[(|u^{k}|^{p-1}|v^{k}|^{p+1} - |u^{j}|^{p-1}|v^{j}|^{p+1})u^{j}\right](\overline{u^{k} - u^{j}})dx \nonumber \\
=\ & 2\beta Im\int_{\mathbb{R}}\left[|u^{k}|^{p-1}(|v^{k}|^{p+1}-|v^{j}|^{p+1}) + |v^{j}|^{p+1}(|u^k|^{p-1}-|u^{j}|^{p-1})\right]u^{j}(\overline{u^{k} - u^{j}})dx \nonumber \\
\label{ecu4}\leq\ & 2\beta\int_{\mathbb{R}}\left[|u^{k}|^{p-1}\left||v^{k}|^{p+1}-|v^{j}|^{p+1}\right| + |v^{j}|^{p+1}\left||u^k|^{p-1}-|u^{j}|^{p-1}\right|\right]|u^{j}|\ |\overline{u^{k} - u^{j}}|dx
\end{align}
Using the identity
\begin{align*}
|v^k|^{p+1}-|v^j|^{p+1}&=(|v^k|-|v^j|)(|v^k|^p+|v^k|^{p-1}|v^j|+...+|v^k||v^j|^{p-1}+|v^j|^p),\\
|u^k|^{p-1}-|u^j|^{p-1}&=(|u^k|-|u^j|)(|u^k|^{p-2}+|u^k|^{p-3}|u^j|+...+|u^k||u^j|^{p-3}+|u^j|^{p-2}).
\end{align*}
Using $||a|-|b||\leq|a-b|$ and the identity, of \eqref{ecu4}, we have
\begin{align}
&2\beta Im\int_{\mathbb{R}}\left[(|u^{k}|^{p-1}|v^{k}|^{p+1} - |u^{j}|^{p-1}|v^{j}|^{p+1})u^{j}\right](\overline{u^{k} - u^{j}})dx \nonumber \\
\label{ecu5}\leq\ & 2\beta M_4\int_{\mathbb{R}}|v^k-v^j|\ |u^k-u^j|dx+2\beta M_5\int_\mathbb{R}|u^k-u^j|^2dx,
\end{align}
where
\begin{align*}
M_4=&\|u^k\|_{L^\infty(\mathbb{R})}^{p-1}\|u^j\|_{L^\infty(\mathbb{R})}(\|v^k\|_{L^\infty(\mathbb{R})}^p+
\|v^k\|_{L^\infty(\mathbb{R})}^{p-1}\|v^j\|_{L^\infty(\mathbb{R})}+...+\|v^k\|_{L^\infty(\mathbb{R})}\|v^j\|_{L^\infty(\mathbb{R})}^{p-1}\\
&+\|v^j\|_{L^\infty(\mathbb{R})}^p)\\
M_5=&\|v^j\|_{L^\infty(\mathbb{R})}^{p+1}\|u^j\|_{L^\infty(\mathbb{R})}(\|u^k\|_{L^\infty(\mathbb{R})}^{p-2}+
\|u^k\|_{L^\infty(\mathbb{R})}^{p-3}\|u^j\|_{L^\infty(\mathbb{R})}+...+\|u^k\|_{L^\infty(\mathbb{R})}\|u^j\|_{L^\infty(\mathbb{R})}^{p-3}\\
&+\|u^j\|_{L^\infty(\mathbb{R})}^{p-2}).
\end{align*}
By Young  inequality, of \eqref{ecu5}, we obtain
\begin{align}
&2\beta Im\int_{\mathbb{R}}\left[(|u^{k}|^{p-1}|v^{k}|^{p+1} - |u^{j}|^{p-1}|v^{j}|^{p+1})u^{j}\right](\overline{u^{k} - u^{j}})dx \nonumber \\
\leq\ & \beta M_4\int_{\mathbb{R}}|v^k-v^j|^2dx+(\beta M_4+2\beta M_5)\int_{\mathbb{R}} |u^k-u^j|^2dx \nonumber\\
\label{ecu6}\leq\ &  M_6\left(\int_{\mathbb{R}}|v^k-v^j|^2dx+\int_{\mathbb{R}} |u^k-u^j|^2dx\right)
\end{align}
where $M_6=max\left\{\beta M_4,\beta M_4+2\beta M_5\right\}$. Similarly we obtain
\begin{align}
&2\beta Im\int_{\mathbb{R}}\left[(|v^{k}|^{p-1}|u^{k}|^{p+1} - |v^{j}|^{p-1}|u^{j}|^{p+1})v^{j}\right](\overline{v^{k} - v^{j}})dx \nonumber \\
\label{ecu7}\leq\ &  \widetilde{M}_{6}\left(\int_{\mathbb{R}}|u^k-u^j|^2dx+\int_{\mathbb{R}} |v^k-v^j|^2dx\right)
\end{align}
Gathering \eqref{506}, \eqref{509}, \eqref{510},  \eqref{ecu6} and \eqref{ecu7}
\begin{eqnarray}
\lefteqn{\frac{d}{dt}\|u^{k} - u^{j}\|_{L^{2}(\mathbb{R})}^{2} } \nonumber \\
& \leq & M_{3}\|u^{k} - u^{j}\|_{L^{2}(\mathbb{R})}^{2} \nonumber
+\ M_6\left(\|v^{k} - v^{j}\|_{L^{2}(\mathbb{R})}^{2} + \|u^{k} - u^{j}\|_{L^{2}(\mathbb{R})}^{2}\right)\nonumber\\
\label{515}& \leq &\ M_{7}\left(\|v^{k} - v^{j}\|_{L^{2}(\mathbb{R})}^{2} + \|u^{k} - u^{j}\|_{L^{2}(\mathbb{R})}^{2}\right)
\end{eqnarray}
where $M_{7}=max\{M_3,M_6\}$. Similarly,
\begin{align}\label{516}
\frac{d}{dt}\|v^{k} - v^{j}\|_{L^{2}(\mathbb{R})}^{2}
\leq \ \widetilde{M}_{7}\left(\|v^{k} - v^{j}\|_{L^{2}(\mathbb{R})}^{2} + \|u^{k} - u^{j}\|_{L^{2}(\mathbb{R})}^{2}\right)
\end{align}
Adding \eqref{515} and \eqref{516} we obtain
\begin{eqnarray}
\label{517}\frac{d}{dt}\left[\|u^{k} - u^{j}\|_{L^{2}(\mathbb{R})}^{2} + \|v^{k} - v^{j}\|_{L^{2}(\mathbb{R})}^{2}\right]
\leq  C\left[\|u^{k} - u^{j}\|_{L^{2}(\mathbb{R})}^{2} + \|v^{k} - v^{j}\|_{L^{2}(\mathbb{R})}^{2}\right],
\end{eqnarray}
where $C=max\{M_{7},\widetilde{M}_{7}\}$.\\
On the other hand, differentiating \eqref{503} with respect to the $x$-variable, we have
\begin{eqnarray}
\lefteqn{i(u^{k} - u^{j})_{xt} + (u^{k} - u^{j})_{xxx} } \nonumber \\
\label{518} & = & (|u^{k}|^{2p}u^{k} -  |u^{j}|^{2p}u^{j})_{x} + \beta(|u^{k}|^{p - 1}|v^{k}|^{p + 1}u^{k} - |u^{j}|^{p - 1}|v^{j}|^{p + 1}u^{j})_{x}.
\end{eqnarray}
Multiplying \eqref{518} by $(\overline{u^{k} - u^{j}})_{x}$ it follows that
\begin{eqnarray*}
\lefteqn{i(\overline{u^{k} - u^{j}})_{x}(u^{k} - u^{j})_{xt} + (\overline{u^{k} - u^{j}})_{x}(u^{k} - u^{j})_{xxx} } \nonumber \\
& = & (|u^{k}|^{2p}u^{k} -  |u^{j}|^{2p}u^{j})_{x}(\overline{u^{k} - u^{j}})_{x} \nonumber \\
& & +\ \beta(|u^{k}|^{p - 1}|v^{k}|^{p + 1}u^{k} - |u^{j}|^{p - 1}|v^{j}|^{p + 1}u^{j})_{x}(\overline{u^{k} - u^{j}})_{x}
\end{eqnarray*}
and its conjugate
\begin{eqnarray*}
\lefteqn{-i(u^{k} - u^{j})_{x}(\overline{u^{k} - u^{j}})_{xt} + (u^{k} - u^{j})_{x}(\overline{u^{k} - u^{j}})_{xxx}}  \nonumber \\
& = & \overline{(|u^{k}|^{2p}u^{k} -  |u^{j}|^{2p}u^{j})_{x}(\overline{u^{k} - u^{j}}})_{x} \nonumber \\
& & +\ \beta\overline{(|u^{k}|^{p - 1}|v^{k}|^{p + 1}u^{k} - |u^{j}|^{p - 1}|v^{j}|^{p + 1}u^{j})_{x}(\overline{u^{k} - u^{j}})}_{x}.
\end{eqnarray*}
Subtracting and integrating over $x\in\mathbb{R}$ we have
\begin{eqnarray}
\lefteqn{\frac{d}{dt}\|(u^{k} - u^{j})_{x}\|_{L^{2}(\mathbb{R})}^{2} }\nonumber \\
& = & 2 Im\int_{\mathbb{R}}(|u^{k}|^{2p}u^{k} -  |u^{j}|^{2p}u^{j})_{x}(\overline{u^{k} - u^{j}})_{x}dx \nonumber \\
\label{519}& & +\ 2\beta Im\int_{\mathbb{R}}(|u^{k}|^{p - 1}|v^{k}|^{p + 1}u^{k} - |u^{j}|^{p - 1}|v^{j}|^{p + 1}u^{j})_{x}(\overline{u^{k} - u^{j}})_{x}dx.
\end{eqnarray}
In a similar way we obtain that
\begin{eqnarray}
\lefteqn{\frac{d}{dt}\|(v^{k} - v^{j})_{x}\|_{L^{2}(\mathbb{R})}^{2} } \nonumber \\
& = & 2Im\int_{\mathbb{R}}(|v^{k}|^{2p}v^{k} -  |v^{j}|^{2p}v^{j})_{x}(\overline{v^{k} - v^{j}})_{x}dx \nonumber \\
\label{520}& &+\ 2\beta Im\int_{\mathbb{R}}(|v^{k}|^{p - 1}|u^{k}|^{p + 1}v^{k} - |v^{j}|^{p - 1}|u^{j}|^{p + 1}v^{j})_{x}(\overline{v^{k} - v^{j}})_{x}dx.
\end{eqnarray}
We estimate the first term in \eqref{519}
\begin{align*}
\lefteqn{(|u^{k}|^{2p}u^{k} -  |u^{j}|^{2p}u^{j})_{x}(\overline{u^{k} - u^{j}})_{x} } \nonumber \\
= &\ \left[|u^{k}|^{2p}(u^{k} - u^{j}) + (|u^{k}|^{2p} - |u^{j}|^{2p})u^{j}\right]_{x}(\overline{u^{k} - u^{j}})_{x} \nonumber \\
= &\ \left[(|u^{k}|^{2p})_{x}(u^{k} - u^{j}) + |u^{k}|^{2p}(u^{k} - u^{j})_{x}\right](\overline{u^{k} - u^{j}})_{x} \nonumber \\
& +\ \left[(|u^{k}|^{2p} - |u^{j}|^{2p})_{x}u^{j} + (|u^{k}|^{2p} - |u^{j}|^{2p})u_{x}^{j}\right](\overline{u^{k} - u^{j}})_{x} \nonumber \\
= &\ (|u^{k}|^{2p})_{x}(u^{k} - u^{j})(\overline{u^{k} - u^{j}})_{x} + |u^{k}|^{2p}\,|(u^{k} - u^{j})_{x}|^{2} \nonumber \\
& +\ (|u^{k}|^{2p} - |u^{j}|^{2p})_{x}u^{j}(\overline{u^{k} - u^{j}})_{x} + (|u^{k}|^{p} + |u^{j}|^{p})(|u^{k}|^{p} - |u^{j}|^{p})u_{x}^{j}(\overline{u^{k} - u^{j}})_{x} \nonumber \\
= &\ (|u^{k}|^{2p})_{x}(u^{k} - u^{j})(\overline{u^{k} - u^{j}})_{x} + |u^{k}|^{2p}\,|(u^{k} - u^{j})_{x}|^{2} \nonumber \\
& +\ u^{j}\left[(|u^{k}|^{2p})_{x} - (|u^{j}|^{2p})_{x}\right](\overline{u^{k} - u^{j}})_{x} + (|u^{k}|^{p} + |u^{j}|^{p})(|u^{k}|^{p} - |u^{j}|^{p})u_{x}^{j}(\overline{u^{k} - u^{j}})_{x}.
\end{align*}
But, $(|u^{k}|^{2p})_{x} = 2Re[(\overline{u}^{k})^{p}((u^{k})^{p})_{x}].$ Then
\begin{align}
\lefteqn{(|u^{k}|^{2p}u^{k} -  |u^{j}|^{2p}u^{j})_{x}(\overline{u^{k} - u^{j}})_{x} } \nonumber \\
= & \  |u^{k}|^{2p}\,|(u^{k} - u^{j})_{x}|^{2} + (|u^{k}|^{p}) + |u^{j}|^{p})(|u^{k}|^{p} - |u^{j}|^{p})u_{x}^{j}(\overline{u^{k} - u^{j}})_{x} \nonumber \\
& + \ 2Re[(\overline{u}^{k})^{p}((u^{k})^{p})_{x}](u^{k} - u^{j})(\overline{u^{k} - u^{j}})_{x}  \nonumber \\
\label{521}& + \ 2u^{j}Re\left[(\overline{u}^{k})^{p}((u^{k})^{p})_{x} - (\overline{u}^{j})^{p}((u^{j})^{p})_{x}\right](\overline{u^{k} - u^{j}})_{x}.
\end{align}
Integrating and taking imaginary part over \eqref{521}, we have
\begin{align*}
\lefteqn{Im\int_{\mathbb{R}}(|u^{k}|^{2p}u^{k} -  |u^{j}|^{2p}u^{j})_{x}(\overline{u^{k} - u^{j}})_{x}dx } \nonumber \\
\leq &  \int_{\mathbb{R}}|\,|u^{k}|^{p} + |u^{j}|^{p}\,|\ |\,|u^{k}|^{p} - |u^{j}|^{p}\,|\ |u_{x}^{j}|\ |(u^{k} - u^{j})_{x}|dx \nonumber \\
& + \,2\int_{\mathbb{R}}|(u^{k})^{p}|\,|((u^{k})^{p})_{x}|\ |u^{k} - u^{j}|\,|(u^{k} - u^{j})_{x}|dx  \nonumber \\
& + \,2\,\int_{\mathbb{R}}|u^{j}|\left|(u^{k})^{p}\ ((u^{k})^{p})_{x} - (u^{j})^{p}((u^{j})^{p})_{x}\right||(u^{k} - u^{j})_{x}|dx \nonumber \\
\leq &  \left(\|\,|u^{k}|^{p} - |u^{j}|^{p}\,\|_{L^{\infty}(\mathbb{R})}\right)\left(\|\,|u^{k}|^{p}\,\|_{L^{\infty}(\mathbb{R})} + \|\,|u^{j}|^{p}\,\|_{L^{\infty}(\mathbb{R})}\right)\|u_{x}^{j}\|_{L^{2}(\mathbb{R})}\|(u^{k} - u^{j})_{x}\|_{L^{2}(\mathbb{R})} \nonumber \\
& + \,2\|(u^{k})^{p}\|_{L^{\infty}(\mathbb{R})}\,\|u^{k} - u^{j}\|_{L^{\infty}(\mathbb{R})}\|((u^{k})^{p})_{x}\|_{L^{2}(\mathbb{R})}\ \|(u^{k} - u^{j})_{x}\|_{L^{2}(\mathbb{R})}  \nonumber \\
& + \,2\|u^{j}\|_{L^{\infty}(\mathbb{R})}\|(u^{k})^{p}\|_{L^{\infty}(\mathbb{R})}\|((u^{k})^{p})_{x}\|_{L^{2}(\mathbb{R})}\,\|(u^{k} - u^{j})_{x}\|_{L^{2}(\mathbb{R})} \nonumber \\
& + \,2\|u^{j}\|_{L^{\infty}(\mathbb{R})}\|(u^{j})^{p}\|_{L^{\infty}(\mathbb{R})}\|((u^{j})^{p})_{x}\|_{L^{2}(\mathbb{R})}\,\|(u^{k} - u^{j})_{x}\|_{L^{2}(\mathbb{R})}.
\end{align*}
Using that $2ab \leq a^{2} + b^{2}$ it follows that
\begin{align}
&Im\int_{\mathbb{R}}(|u^{k}|^{2p}u^{k} -  |u^{j}|^{2p}u^{j})_{x}(\overline{u^{k} - u^{j}})_{x}dx
\leq   N_{1}\left(\|u_{x}^{j}\|_{L^{2}(\mathbb{R})}^{2} + \|(u^{k} - u^{j})_{x}\|_{L^{2}(\mathbb{R})}^{2}\right) \nonumber \\
&+ N_{2}\left(\|((u^{k})^{p})_{x}\|_{L^{2}(\mathbb{R})}^{2} + \|(u^{k} - u^{j})_{x}\|_{L^{2}(\mathbb{R})}^{2}\right) + N_{3}\left(\|((u^{k})^{p})_{x}\|_{L^{2}(\mathbb{R})}^{2} + \|(u^{k} - u^{j})_{x}\|_{L^{2}(\mathbb{R})}^{2}\right) \nonumber \\
\label{522}&+ N_{4}\left(\|((u^{j})^{p})_{x}\|_{L^{2}(\mathbb{R})}^{2} + \|(u^{k} - u^{j})_{x}\|_{L^{2}(\mathbb{R})}\right)
\end{align}
where
\begin{align*}
& N_{1} = \left(\|\,|u^{k}|^{p} - |u^{j}|^{p}\,\|_{L^{\infty}(\mathbb{R})}\right)\left(\|\,|u^{k}|^{p}\,\|_{L^{\infty}(\mathbb{R})} + \|\,|u^{j}|^{p}\,\|_{L^{\infty}(\mathbb{R})}\right),  \\
& N_{2} = \|(u^{k})^{p}\|_{L^{\infty}(\mathbb{R})}\,\|u^{k} - u^{j}\|_{L^{\infty}(\mathbb{R})}, \\
& N_{3} = \|u^{j}\|_{L^{\infty}(\mathbb{R})}\|(u^{k})^{p}\|_{L^{\infty}(\mathbb{R})}, \qquad  N_{4} = \|u^{j}\|_{L^{\infty}(\mathbb{R})}\|(u^{j})^{p}\|_{L^{\infty}(\mathbb{R})}.
\end{align*}
In a similar way, we estimate the first term in \eqref{520}, that is,
\begin{align}
&Im\int_{\mathbb{R}}(|v^{k}|^{2p}v^{k} -  |v^{j}|^{2p}v^{j})_{x}(\overline{v^{k} - v^{j}})_{x}dx
\leq   \widetilde{N}_{1}\left(\|v_{x}^{j}\|_{L^{2}(\mathbb{R})}^{2} + \|(v^{k} - v^{j})_{x}\|_{L^{2}(\mathbb{R})}^{2}\right) \nonumber \\
&+ \widetilde{N}_{2}\left(\|((v^{k})^{p})_{x}\|_{L^{2}(\mathbb{R})}^{2} + \|(v^{k} - v^{j})_{x}\|_{L^{2}(\mathbb{R})}^{2}\right) + \widetilde{N}_{3}\left(\|((v^{k})^{p})_{x}\|_{L^{2}(\mathbb{R})}^{2} + \|(v^{k} - v^{j})_{x}\|_{L^{2}(\mathbb{R})}^{2}\right) \nonumber \\
\label{523}&+ \widetilde{N}_{4}\left(\|((v^{j})^{p})_{x}\|_{L^{2}(\mathbb{R})}^{2} + \|(v^{k} - v^{j})_{x}\|_{L^{2}(\mathbb{R})}\right)
\end{align}
where
\begin{align*}
& \widetilde{N}_{1} = \left(\|\,|v^{k}|^{p} - |v^{j}|^{p}\,\|_{L^{\infty}(\mathbb{R})}\right)\left(\|\,|v^{k}|^{p}\,\|_{L^{\infty}(\mathbb{R})} + \|\,|v^{j}|^{p}\,\|_{L^{\infty}(\mathbb{R})}\right),  \\
& \widetilde{N}_{2} = \|(v^{k})^{p}\|_{L^{\infty}(\mathbb{R})}\,\|v^{k} - v^{j}\|_{L^{\infty}(\mathbb{R})}, \\
& \widetilde{N}_{3} = \|v^{j}\|_{L^{\infty}(\mathbb{R})}\|(v^{k})^{p}\|_{L^{\infty}(\mathbb{R})}, \qquad  \widetilde{N}_{4} = \|v^{j}\|_{L^{\infty}(\mathbb{R})}\|(v^{j})^{p}\|_{L^{\infty}(\mathbb{R})}.
\end{align*}
Now we estimate the second term in \eqref{519}
\begin{align}
\lefteqn{(|u^{k}|^{p - 1}|v^{k}|^{p + 1}u^{k} - |u^{j}|^{p - 1}|v^{j}|^{p + 1}u^{j})_{x}(\overline{u^{k} - u^{j}})_{x} } \nonumber \\
= &\left[|u^{k}|^{p - 1}|v^{k}|^{p + 1}(u^{k} - u^{j}) + \left(|u^{k}|^{p - 1}|v^{k}|^{p + 1} - |u^{j}|^{p - 1}|v^{j}|^{p + 1}\right)u^{j}\right]_{x}(\overline{u^{k} - u^{j}})_{x} \nonumber \\
= & \left[|u^{k}|^{p - 1}|v^{k}|^{p + 1}(u^{k} - u^{j})_{x} + \left(|u^{k}|^{p - 1}|v^{k}|^{p + 1}\right)_{x}(u^{k} - u^{j})\right](\overline{u^{k} - u^{j}})_{x} \nonumber \\
& + \left[\left(|u^{k}|^{p - 1}|v^{k}|^{p + 1} - |u^{j}|^{p - 1}|v^{j}|^{p + 1}\right)u_{x}^{j}\right](\overline{u^{k} - u^{j}})_{x} \nonumber \\
& + \left[\left(|u^{k}|^{p - 1}|v^{k}|^{p + 1} - |u^{j}|^{p - 1}|v^{j}|^{p + 1}\right)_{x}u^{j}\right](\overline{u^{k} - u^{j}})_{x} \nonumber \\
= & |u^{k}|^{p - 1}|v^{k}|^{p + 1}|(u^{k} - u^{j})_{x}|^{2} \nonumber \\
& + \left(|u^{k}|^{p - 1}\right)_{x}|v^{k}|^{p + 1}(u^{k} - u^{j})(\overline{u^{k} - u^{j}})_{x} + |u^{k}|^{p - 1}\left(|v^{k}|^{p + 1}\right)_{x}(u^{k} - u^{j})(\overline{u^{k} - u^{j}})_{x} \nonumber \\
& + \left[\left(|u^{k}|^{p - 1}|v^{k}|^{p + 1} - |u^{j}|^{p - 1}|v^{j}|^{p + 1}\right)u_{x}^{j}\right](\overline{u^{k} - u^{j}})_{x} \nonumber \\
& + \left[\left(|u^{k}|^{p - 1}|v^{k}|^{p + 1}\right)_{x} - \left(|u^{j}|^{p - 1}|v^{j}|^{p + 1}\right)_{x}\right]u^{j}(\overline{u^{k} - u^{j}})_{x} \nonumber \\
= & |u^{k}|^{p - 1}|v^{k}|^{p + 1}|(u^{k} - u^{j})_{x}|^{2} \nonumber \\
& + \left(|u^{k}|^{p - 1}\right)_{x}|v^{k}|^{p + 1}(u^{k} - u^{j})(\overline{u^{k} - u^{j}})_{x} + |u^{k}|^{p - 1}\left(|v^{k}|^{p + 1}\right)_{x}(u^{k} - u^{j})(\overline{u^{k} - u^{j}})_{x} \nonumber \\
& + |u^{k}|^{p - 1}|v^{k}|^{p + 1}u_{x}^{j}(\overline{u^{k} - u^{j}})_{x} - |u^{j}|^{p - 1}|v^{j}|^{p + 1}u_{x}^{j}(\overline{u^{k} - u^{j}})_{x} \nonumber \\
& + \left[\left(|u^{k}|^{p - 1}\right)_{x}|v^{k}|^{p + 1} + |u^{k}|^{p - 1}\left(|v^{k}|^{p + 1}\right)_{x}\right]u^{j}(\overline{u^{k} - u^{j}})_{x} \nonumber \\
\label{524}& - \left(|u^{j}|^{p - 1}\right)_{x}|v^{j}|^{p + 1}u^{j}(\overline{u^{k} - u^{j}})_{x} - |u^{j}|^{p - 1}\left(|v^{j}|^{p + 1}\right)_{x}u^{j}(\overline{u^{k} - u^{j}})_{x}.
\end{align}
But
\begin{align*}
\left(|u^{k}|^{p - 1}\right)_{x} = &\left(|u^{k}|^{2\left(\frac{p - 1}{2}\right)}\right)_{x} = \left[\left(u^{k}\overline{u^{k}}\right)^{\left(\frac{p - 1}{2}\right)}\right]_{x} \nonumber \\
= &\left[\left(u^{k}\right)^{\left(\frac{p - 1}{2}\right)}\right]_{x}\overline{u^{k}}^{\left(\frac{p - 1}{2}\right)} + (u^{k})^{\left(\frac{p - 1}{2}\right)}\left[\left(\overline{u^{k}}\right)^{\left(\frac{p - 1}{2}\right)}\right]_{x} \nonumber \\
= &\left[\left(u^{k}\right)^{\left(\frac{p - 1}{2}\right)}\right]_{x}\overline{u^{k}}^{\left(\frac{p - 1}{2}\right)} + \overline{\left[\left(u^{k}\right)^{\left(\frac{p - 1}{2}\right)}\right]_{x}\overline{u^{k}}^{\left(\frac{p - 1}{2}\right)}} \nonumber \\
= &2Re\left\{\left[\left(u^{k}\right)^{\left(\frac{p - 1}{2}\right)}\right]_{x}\left(\overline{u^{k}}\right)^{\left(\frac{p - 1}{2}\right)}\right\}.
\end{align*}
Similarly
\begin{align*}
\left(|u^{k}|^{p + 1}\right)_{x} = 2Re\left\{\left[\left(u^{k}\right)^{\left(\frac{p + 1}{2}\right)}\right]_{x}\left(\overline{u^{k}}\right)^{\left(\frac{p + 1}{2}\right)}\right\}.
\end{align*}
Replacing into \eqref{524}
\begin{align}
\lefteqn{(|u^{k}|^{p - 1}|v^{k}|^{p + 1}u^{k} - |u^{j}|^{p - 1}|v^{j}|^{p + 1}u^{j})_{x}(\overline{u^{k} - u^{j}})_{x} } \nonumber \\
= & |u^{k}|^{p - 1}|v^{k}|^{p + 1}|(u^{k} - u^{j})_{x}|^{2} \nonumber \\
& + 2Re\left\{\left[\left(u^{k}\right)^{\left(\frac{p - 1}{2}\right)}\right]_{x}\left(\overline{u^{k}}\right)^{\left(\frac{p - 1}{2}\right)}\right\}|v^{k}|^{p + 1}(u^{k} - u^{j})(\overline{u^{k} - u^{j}})_{x} \nonumber \\
& + |u^{k}|^{p - 1}2Re\left\{\left[\left(v^{k}\right)^{\left(\frac{p + 1}{2}\right)}\right]_{x}\left(\overline{v^{k}}\right)^{\left(\frac{p + 1}{2}\right)}\right\}(u^{k} - u^{j})(\overline{u^{k} - u^{j}})_{x} \nonumber \\
& + |u^{k}|^{p - 1}|v^{k}|^{p + 1}u_{x}^{j}(\overline{u^{k} - u^{j}})_{x} - |u^{j}|^{p - 1}|v^{j}|^{p + 1}u_{x}^{j}(\overline{u^{k} - u^{j}})_{x} \nonumber \\
& + 2Re\left\{\left[\left(u^{k}\right)^{\left(\frac{p - 1}{2}\right)}\right]_{x}\left(\overline{u^{k}}\right)^{\left(\frac{p - 1}{2}\right)}\right\}|v^{k}|^{p + 1}u^{j}(\overline{u^{k} - u^{j}})_{x} \nonumber \\
& + 2|u^{k}|^{p - 1}Re\left\{\left[\left(v^{k}\right)^{\left(\frac{p + 1}{2}\right)}\right]_{x}\left(\overline{v^{k}}\right)^{\left(\frac{p + 1}{2}\right)}\right\}u^{j}(\overline{u^{k} - u^{j}})_{x} \nonumber \\
& - 2Re\left\{\left[\left(u^{k}\right)^{\left(\frac{p - 1}{2}\right)}\right]_{x}\left(\overline{u^{k}}\right)^{\left(\frac{p - 1}{2}\right)}\right\}|v^{j}|^{p + 1}u^{j}(\overline{u^{k} - u^{j}})_{x} \nonumber \\
& - 2|u^{j}|^{p - 1}Re\left\{\left[\left(v^{j}\right)^{\left(\frac{p + 1}{2}\right)}\right]_{x}\left(\overline{v^{j}}\right)^{\left(\frac{p + 1}{2}\right)}\right\}u^j(\overline{u^{k} - u^{j}})_{x}.
\end{align}
Then
\begin{align*}
2\beta\lefteqn{Im\int_{\mathbb{R}}(|u^{k}|^{p - 1}|v^{k}|^{p + 1}u^{k} - |u^{j}|^{p - 1}|v^{j}|^{p + 1}u^{j})_{x}(\overline{u^{k} - u^{j}})_{x}dx } \nonumber \\
\leq &\ 2\beta\int_{\mathbb{R}}2\left|[(u^{k})^{(p - 1)/2}]_{x}\right|\,|u^{k}|^{(p - 1)/2)}\ |v^{k}|^{p + 1}\,|u^{k} - u^{j}|\ |(u^{k} - u^{j})_{x}|dx \nonumber \\
& + 2\beta\int_{\mathbb{R}}2|u^{k}|^{p - 1}\left|[(v^{k})^{(p + 1)/2}]_{x}\right|\,|v^{k}|^{(p + 1)/2}\,|u^{k} - u^{j}|\ |(u^{k} - u^{j})_{x}|dx \nonumber \\
& +2\beta \int_{\mathbb{R}}|u^{k}|^{p - 1}|v^{k}|^{p + 1}\,|u_{x}^{j}|\ |(u^{k} - u^{j})_{x}|dx +2\beta \int_{\mathbb{R}}|u^{j}|^{p - 1}|v^{j}|^{p + 1}\,|u_{x}^{j}|\ |(\overline{u^{k} - u^{j}})_{x}|dx \nonumber \\
& + 2\beta\int_{\mathbb{R}}2|[(u^{k})^{(p - 1)/2}]_{x}|\ |u^{k}|^{(p - 1)/2}\ |v^{k}|^{p + 1}\,|u^{j}|\ |(u^{k} - u^{j})_{x}|dx \nonumber \\
& + 2\beta\int_{\mathbb{R}}2|u^{k}|^{p - 1}|[(v^{k})^{(p + 1)/2}]_{x}|\ |v^{k}|^{(p + 1)/2}|\ |u^{j}|\ |(u^{k} - u^{j})_{x}|dx \nonumber \\
& + 2\beta\int_{\mathbb{R}}2|[(u^{k})^{(p - 1)/2)}]_{x}|\ |u^{k}|^{(p - 1)/2}|v^{j}|^{p + 1}|u^{j}|\ |(u^{k} - u^{j})_{x}|dx \nonumber \\
& + 2\beta\int_{\mathbb{R}}2|u^{j}|^{p - 1}\,|[(v^{j})^{(p + 1)/2}]_{x}|\ |v^{j}|^{(p + 1)/2}\ |u^j|\ |(u^{k} - u^{j})_{x}|dx.
\end{align*}
Hence
\begin{align*}
2\beta\lefteqn{Im\int_{\mathbb{R}}(|u^{k}|^{p - 1}|v^{k}|^{p + 1}u^{k} - |u^{j}|^{p - 1}|v^{j}|^{p + 1}u^{j})_{x}(\overline{u^{k} - u^{j}})_{x}dx } \nonumber \\
\leq &\ 4\beta\|u^{k}\|_{L^{\infty}(\mathbb{R})}^{(p - 1)/2)}\|v^{k}\|_{L^{\infty}(\mathbb{R})}^{p + 1}\|u^{k} - u^{j}\|_{L^{\infty}(\mathbb{R})}\int_{\mathbb{R}}\left|[(u^{k})^{(p - 1)/2}]_{x}\right|\ |(u^{k} - u^{j})_{x}|dx \nonumber \\
& + 4\beta\|u^{k}\|_{L^{\infty}(\mathbb{R})}^{p - 1}\|v^{k}\|_{L^{\infty}(\mathbb{R})}^{(p + 1)/2}\,\|u^{k} - u^{j}\|_{L^{\infty}(\mathbb{R})}\int_{\mathbb{R}}\left|[(v^{k})^{(p + 1)/2}]_{x}\right|\ |(u^{k} - u^{j})_{x}|dx \nonumber \\
& + 2\beta\|u^{k}\|_{L^{\infty}(\mathbb{R})}^{p - 1}\|v^{k}\|_{L^{\infty}(\mathbb{R})}^{p + 1}\int_{\mathbb{R}}|u_{x}^{j}||(u^{k} - u^{j})_{x}|dx + 2\beta\|u^{j}\|_{L^{\infty}(\mathbb{R})}^{p - 1}\|v^{j}\|_{L^{\infty}(\mathbb{R})}^{p + 1}\int_{\mathbb{R}}|u_{x}^{j}| |(\overline{u^{k} - u^{j}})_{x}|dx \nonumber \\
& + 4\beta\|u^{k}\|_{L^{\infty}(\mathbb{R})}^{(p - 1)/2}\|v^{k}\|_{L^{\infty}(\mathbb{R})}^{p + 1}\|u^{j}\|_{L^{\infty}(\mathbb{R})}\int_{\mathbb{R}}|[(u^{k})^{(p - 1)/2}]_{x}|\ |(u^{k} - u^{j})_{x}|dx \nonumber \\
& + 4\beta\|u^{k}\|_{L^{\infty}(\mathbb{R})}^{p - 1}\|v^{k}\|_{L^{\infty}(\mathbb{R})}^{(p + 1)/2}|\|u^{j}|\|_{L^{\infty}(\mathbb{R})}\int_{\mathbb{R}}|[(v^{k})^{(p + 1)/2}]_{x}|\ |(u^{k} - u^{j})_{x}|dx \nonumber \\
& + 4\beta\|u^{j}\|_{L^{\infty}(\mathbb{R})}^{(p - 1)/2}\|v^{j}\|_{L^{\infty}(\mathbb{R})}^{p + 1}\|u^{j}\|_{L^{\infty}(\mathbb{R})}\int_{\mathbb{R}}|[(u^{k})^{(p - 1)/2)}]_{x}|\,|(u^{k} - u^{j})_{x}|dx \nonumber \\
& + 4\beta\|u^{j}\|_{L^{\infty}(\mathbb{R})}^{p - 1}\|v^{j}\|_{L^{\infty}(\mathbb{R})}^{(p + 1)/2}\|u^j\|_{L^\infty(\mathbb{R})}\int_{\mathbb{R}}|[(v^{j})^{(p + 1)/2}]_{x}|\ |(u^{k} - u^{j})_{x}|dx.
\end{align*}
Thus
\begin{align}
2\beta\lefteqn{Im\int_{\mathbb{R}}(|u^{k}|^{p - 1}|v^{k}|^{p + 1}u^{k} - |u^{j}|^{p - 1}|v^{j}|^{p + 1}u^{j})_{x}(\overline{u^{k} - u^{j}})_{x}dx } \nonumber \\
\leq &\ 4\beta\|u^{k}\|_{L^{\infty}(\mathbb{R})}^{(p - 1)/2)}\|v^{k}\|_{L^{\infty}(\mathbb{R})}^{p + 1}\|u^{k} - u^{j}\|_{L^{\infty}(\mathbb{R})}\left\|[(u^{k})^{(p - 1)/2}]_{x}\right\|_{L^{2}(\mathbb{R})} \|(u^{k} - u^{j})_{x}\|_{L^{2}(\mathbb{R})} \nonumber \\
& + 4\beta\|u^{k}\|_{L^{\infty}(\mathbb{R})}^{p - 1}\|v^{k}\|_{L^{\infty}(\mathbb{R})}^{(p + 1)/2}\|u^{k} - u^{j}\|_{L^{\infty}(\mathbb{R})}\left\|[(v^{k})^{(p + 1)/2}]_{x}\right\|_{L^{2}(\mathbb{R})}\|(u^{k} - u^{j})_{x}\|_{L^{2}(\mathbb{R})} \nonumber \\
& + 2\beta\|u^{k}\|_{L^{\infty}(\mathbb{R})}^{p - 1}\|v^{k}\|_{L^{\infty}(\mathbb{R})}^{p + 1}\|u_{x}^{j}\|_{L^{2}(\mathbb{R})}\|(u^{k} - u^{j})_{x}\|_{L^{2}(\mathbb{R})} \nonumber \\
&+ 2\beta\|u^{j}\|_{L^{\infty}(\mathbb{R})}^{p - 1}\|v^{j}\|_{L^{\infty}(\mathbb{R})}^{p + 1}\|u_{x}^{j}\|_{L^{2}(\mathbb{R})}\|(\overline{u^{k} - u^{j}})_{x}\|_{L^{2}(\mathbb{R})} \nonumber \\
& + 4\beta\|u^{k}\|_{L^{\infty}(\mathbb{R})}^{(p - 1)/2}\|v^{k}\|_{L^{\infty}(\mathbb{R})}^{p + 1}\|u^{j}\|_{L^{\infty}(\mathbb{R})}\|[(u^{k})^{(p - 1)/2}]_{x}\|_{L^{2}(\mathbb{R})}\|(u^{k} - u^{j})_{x}\|_{L^{2}(\mathbb{R})} \nonumber \\
& + 4\beta\|u^{k}\|_{L^{\infty}(\mathbb{R})}^{p - 1}\|v^{k}\|_{L^{\infty}(\mathbb{R})}^{(p + 1)/2}\|u^{j}\|_{L^{\infty}(\mathbb{R})}\|[(v^{k})^{(p + 1)/2}]_{x}\|_{L^{2}(\mathbb{R})}\|(u^{k} - u^{j})_{x}\|_{L^{2}(\mathbb{R})} \nonumber \\
& + 4\beta\|u^{k}\|_{L^{\infty}(\mathbb{R})}^{(p - 1)/2}\|v^{k}\|_{L^{\infty}(\mathbb{R})}^{p + 1}\|u^{j}\|_{L^{\infty}(\mathbb{R})}\|[(u^{k})^{(p - 1)/2)}]_{x}\|_{L^{2}(\mathbb{R})}\|(u^{k} - u^{j})_{x}\|_{L^{2}(\mathbb{R})} \nonumber \\
& + 4\beta\|u^{j}\|_{L^{\infty}(\mathbb{R})}^{p - 1}\|v^{j}\|_{L^{\infty}(\mathbb{R})}^{(p + 1)/2}\|[(v^{j})^{(p + 1)/2}]_{x}\|_{L^{2}(\mathbb{R})}\|(u^{k} - u^{j})_{x}\|_{L^{2}(\mathbb{R})}.
\end{align}
Therefore
\begin{align}
2\beta\lefteqn{Im\int_{\mathbb{R}}(|u^{k}|^{p - 1}|v^{k}|^{p + 1}u^{k} - |u^{j}|^{p - 1}|v^{j}|^{p + 1}u^{j})_{x}(\overline{u^{k} - u^{j}})_{x}dx } \nonumber \\
\leq &\ R_{1}\left(\left\|[(u^{k})^{(p - 1)/2}]_{x}\right\|_{L^{2}(\mathbb{R})}^{2} + \|(u^{k} - u^{j})_{x}\|_{L^{2}(\mathbb{R})}^{2}\right) \nonumber \\
& + R_{2}\left(\left\|[(v^{k})^{(p + 1)/2}]_{x}\right\|_{L^{2}(\mathbb{R})}^{2} + \|(u^{k} - u^{j})_{x}\|_{L^{2}(\mathbb{R})}^{2}\right) \nonumber \\
& + R_{3}\left(\|u_{x}^{j}\|_{L^{2}(\mathbb{R})}^{2} + \|(u^{k} - u^{j})_{x}\|_{L^{2}(\mathbb{R})}^{2}\right) \nonumber \\
&+ R_{4}\left(\|u_{x}^{j}\|_{L^{2}(\mathbb{R})}^{2} + \|(u^{k} - u^{j})_{x}\|_{L^{2}(\mathbb{R})}^{2}\right) \nonumber \\
& + R_{5}\left(\|[(u^{k})^{(p - 1)/2}]_{x}\|_{L^{2}(\mathbb{R})}^{2} + \|(u^{k} - u^{j})_{x}\|_{L^{2}(\mathbb{R})}^{2}\right) \nonumber \\
& + R_{6}\left(\|[(v^{k})^{(p + 1)/2}]_{x}\|_{L^{2}(\mathbb{R})}^{2} + \|(u^{k} - u^{j})_{x}\|_{L^{2}(\mathbb{R})}^{2}\right) \nonumber \\
& + R_{7}\left(\|[(u^{j})^{(p - 1)/2)}]_{x}\|_{L^{2}(\mathbb{R})}^{2} + \|(u^{k} - u^{j})_{x}\|_{L^{2}(\mathbb{R})}^{2}\right) \nonumber \\
\label{527}& + R_{8}\left(\|[(v^{j})^{(p + 1)/2}]_{x}\|_{L^{2}(\mathbb{R})}^{2} + \|(u^{k} - u^{j})_{x}\|_{L^{2}(\mathbb{R})}^{2}\right),
\end{align}
where
\begin{align*}
& R_{1} = 2\beta\|u^{k}\|_{L^{\infty}(\mathbb{R})}^{(p - 1)/2)}\|v^{k}\|_{L^{\infty}(\mathbb{R})}^{p + 1}\|u^{k} - u^{j}\|_{L^{\infty}(\mathbb{R})},\  R_{2} = 2\beta\|u^{k}\|_{L^{\infty}(\mathbb{R})}^{p - 1}\|v^{k}\|_{L^{\infty}(\mathbb{R})}^{(p + 1)/2}\|u^{k} - u^{j}\|_{L^{\infty}(\mathbb{R})} \\
& R_{3} =\beta \|u^{k}\|_{L^{\infty}(\mathbb{R})}^{p - 1}\|v^{k}\|_{L^{\infty}(\mathbb{R})}^{p + 1}, \quad R_{4} = \beta\|u^{j}\|_{L^{\infty}(\mathbb{R})}^{p - 1}\|v^{j}\|_{L^{\infty}(\mathbb{R})}^{p + 1} \\
& R_{5} = 2\beta\|u^{k}\|_{L^{\infty}(\mathbb{R})}^{(p - 1)/2}\|v^{k}\|_{L^{\infty}(\mathbb{R})}^{p + 1}\|u^{j}\|_{L^{\infty}(\mathbb{R})}, \quad R_{6} = 2\beta\|u^{k}\|_{L^{\infty}(\mathbb{R})}^{p - 1}\|v^{k}\|_{L^{\infty}(\mathbb{R})}^{(p + 1)/2}\|u^{j}\|_{L^{\infty}(\mathbb{R})}  \\
& R_{7} = 2\beta\|u^{j}\|_{L^{\infty}(\mathbb{R})}^{(p - 1)/2}\|v^{j}\|_{L^{\infty}(\mathbb{R})}^{p + 1}\|u^{j}\|_{L^{\infty}(\mathbb{R})}, \quad R_{8} = 2\beta\|u^{j}\|_{L^{\infty}(\mathbb{R})}^{p - 1}\|v^{j}\|_{L^{\infty}(\mathbb{R})}^{(p + 1)/2}\|u^j\|_{L^\infty(\mathbb{R})}.
\end{align*}
Therefore, from (\ref{527})
\begin{align}
&2\beta Im\int_{\mathbb{R}}(|u^{k}|^{p - 1}|v^{k}|^{p + 1}u^{k} - |u^{j}|^{p - 1}|v^{j}|^{p + 1}u^{j})_{x}(\overline{u^{k} - u^{j}})_{x}dx  \nonumber \\
\label{ecu8}\leq &\ K_1+K_2\|(u^k-u^j)_x\|_{L^2(\mathbb{R})}^2.
\end{align}
Similarly
\begin{align}
&2\beta Im\int_{\mathbb{R}}(|v^{k}|^{p - 1}|u^{k}|^{p + 1}v^{k} - |v^{j}|^{p - 1}|u^{j}|^{p + 1}v^{j})_{x}(\overline{v^{k} - v^{j}})_{x}dx  \nonumber \\
\label{ecu9}\leq &\ \widetilde{K}_1+\widetilde{K}_2\|(v^k-v^j)_x\|_{L^2(\mathbb{R})}^2.
\end{align}
Gathering (\ref{519}), (\ref{520}), (\ref{522}), (\ref{523}), (\ref{ecu8}) and (\ref{ecu9}), we obtain
\begin{align}
&\dfrac{d}{dt}\left[\|(u^k-u^j)_x\|_{L^2(\mathbb{R})}^2+\|(v^k-v^j)_x\|_{L^2(\mathbb{R})}^2\right]\nonumber\\
&\label{ecu10}\leq\ C_0+C_1\left[\|(u^k-u^j)_x\|_{L^2(\mathbb{R})}^2+\|(v^k-v^j)_x\|_{L^2(\mathbb{R})}^2 \right].
\end{align}
From (\ref{517}) and (\ref{ecu10}) the Lemma following.
\begin{lemma}
\label{lemma502}
Let $p>1$ odd integer number, for $m=1,2,3,...,n$, we have $\{(J^m u^k,J^m u^k)\}$ is Cauchy sequence in $C([0,T]: L^2(\mathbb{R}))\times C([0,T]: L^2(\mathbb{R}))$ for any $T>0$. Moreover
\begin{align*}
&\|J^m u^k-J^m u^j\|_{L^2(\mathbb{R})}^2+\|J^m v^k-J^m v^j\|_{L^2(\mathbb{R})}^2\\
&\ \leq c(T)\left[\|x^m u_0^k-x^m u_0^j\|_{L^2(\mathbb{R})}^2+\|x^m v_0^k-x^m v_0^j\|_{L^2(\mathbb{R})}^2\right]
\end{align*}
where $c(T)$ is a positive constant independent of $k$ and $j$.
\end{lemma}
\proof Let $(u^k,u^j)$ solution, then
\begin{eqnarray}
\label{ecu11} \left\lbrace
\begin{array}{l}
i(J^mu^k)_t+(J^mu^k)_{xx} = J^m[(|u^k|^{2p}+\beta|u^k|^{p-1}|v^k|^{p+1})u^k],\\
i(J^mu^j)_t+(J^mu^j)_{xx} = J^m[(|u^j|^{2p}+\beta|u^j|^{p-1}|v^j|^{p+1})u^j],\\
i(J^mv^k)_t+(J^mv^k)_{xx} = J^m[(|v^k|^{2p}+\beta|v^k|^{p-1}|u^k|^{p+1})v^k],\\
i(J^mv^j)_t+(J^mv^j)_{xx} = J^m[(|v^j|^{2p}+\beta|v^j|^{p-1}|u^j|^{p+1})v^j],
\end{array}
\right.
\end{eqnarray}
we have
\begin{align*}
i(J^mu^k-J^mu^j)_t+(J^mu^k-J^mu^j)_{xx}=& J^m(|u^k|^{2p}u^k)+\beta J^m(|u^k|^{p-1}|v^k|^{p+1}u^k)\\
&-J^m(|u^j|^{2p}u^j)-\beta J^m(|u^j|^{p-1}|v^j|^{p+1}u^j).
\end{align*}
Multiplying by $\overline{J^m u^k-J^m u^j}$, we obtain
\begin{align*}
&i(J^mu^k-J^mu^j)_t(\overline{J^m u^k-J^m u^j})+(J^mu^k-J^mu^j)_{xx}(\overline{J^m u^k-J^m u^j})\\
&\ =[J^m(|u^k|^{2p}u^k)+\beta J^m(|u^k|^{p-1}|v^k|^{p+1}u^k)](\overline{J^m u^k-J^m u^j})\\
&\ \ \  +[-J^m(|u^j|^{2p}u^j)-\beta J^m(|u^j|^{p-1}|v^j|^{p+1}u^j)](\overline{J^m u^k-J^m u^j}).
\end{align*}
Applying conjugate
\begin{align*}
&-i(\overline{J^mu^k-J^mu^j})_t(J^m u^k-J^m u^j)+(\overline{J^mu^k-J^mu^j})_{xx}(J^m u^k-J^m u^j)\\
&\ =\overline{[J^m(|u^k|^{2p}u^k)+\beta J^m(|u^k|^{p-1}|v^k|^{p+1}u^k)](\overline{J^m u^k-J^m u^j})}\\
&\ \ \  +\overline{[-J^m(|u^j|^{2p}u^j)-\beta J^m(|u^j|^{p-1}|v^j|^{p+1}u^j)](\overline{J^m u^k-J^m u^j})}.
\end{align*}
Subtracting and integrating over $x\in\mathbb{R}$, we have
\begin{align*}
&\frac{d}{dt}\|J^mu^k-J^mu^j\|_{L^2(\mathbb{R})}^2
=2\ Im\left[\int_{\mathbb{R}}\left[J^m(|u^k|^{2p}u^k)-J^m(|u^j|^{2p}u^j)\right](\overline{J^m u^k-J^m u^j})\ dx\right]\\
&\ +2\beta \ Im\left[\int_{\mathbb{R}}[J^m(|u^k|^{p-1}|v^k|^{p+1}u^k)-J^m(|u^j|^{p-1}|v^j|^{p+1}u^j)]
(\overline{J^m u^k-J^m u^j})\ dx\right].
\end{align*}
By the H\"{o}lder inequality
\begin{align}\label{ecu12}
&\frac{d}{dt}\|J^mu^k-J^mu^j\|_{L^2(\mathbb{R})}^2
\leq C \left[\|J^m(|u^k|^{2p}u^k)\|_{L^2(\mathbb{R})}+\|J^m(|u^j|^{2p}u^j)\|_{L^2(\mathbb{R})}\right]\|J^m u^k-J^m u^j\|_{L^2(\mathbb{R})}\nonumber\\
&\ +C\left[\|J^m(|u^k|^{p-1}|v^k|^{p+1}u^k)\|_{L^2(\mathbb{R})}+\|J^m(|u^j|^{p-1}|v^j|^{p+1}u^j)\|_{L^2(\mathbb{R})}\right]
\|J^m u^k-J^m u^j\|_{L^2(\mathbb{R})}.
\end{align}
Using the  Lemma (\ref{lema1}) and  (\ref{lema2}), of (\ref{ecu12}), we have
\begin{align}\label{ecu13}
\frac{d}{dt}\|J^mu^k-J^mu^j\|_{L^2(\mathbb{R})}^2
\leq& C_m
\left[\|J^m u^k\|_{L^2(\mathbb{R})}\|u^k\|_{L^\infty(\mathbb{R})}^{2p}+\|J^m u^j\|_{L^2(\mathbb{R})}\|u^j\|_{L^\infty(\mathbb{R})}^{2p}\right]\nonumber\\
&\times\|J^m u^k-J^m u^j\|_{L^2(\mathbb{R})}\nonumber\\
&+C_m\left[\|J^m u^k\|_{L^2(\mathbb{R})}+\|J^m v^k\|_{L^2(\mathbb{R})}+\|J^m u^j\|_{L^2(\mathbb{R})}+\|J^m v^j\|_{L^2(\mathbb{R})}\right]\nonumber\\
&\times\|J^m u^k-J^m u^j\|_{L^2(\mathbb{R})}\nonumber\\
&\leq C_m\|J^m u^k-J^m u^j\|_{L^2(\mathbb{R})}.
\end{align}
From (\ref{ecu13}) by the Young inequality, we obtain
\begin{align}\label{ecu14}
\frac{d}{dt}\|J^mu^k-J^mu^j\|_{L^2(\mathbb{R})}^2
\leq C_m+\|J^m u^k-J^m u^j\|_{L^2(\mathbb{R})}^2.
\end{align}
Similar way we obtain for $v$ that
\begin{align}\label{ecu15}
\frac{d}{dt}\|J^mv^k-J^mv^j\|_{L^2(\mathbb{R})}^2
\leq \widetilde{C_m}+\|J^m v^k-J^m v^j\|_{L^2(\mathbb{R})}^2.
\end{align}
Adding (\ref{ecu14}) and (\ref{ecu15}) we obtain
\begin{align}\label{ecu16}
\frac{d}{dt}\left[\|J^mu^k-J^mu^j\|_{L^2(\mathbb{R})}^2+\|J^mv^k-J^mv^j\|_{L^2(\mathbb{R})}^2\right]
\leq& C_m+\|J^m u^k-J^m u^j\|_{L^2(\mathbb{R})}^2\nonumber\\
&+\|J^m v^k-J^m v^j\|_{L^2(\mathbb{R})}^2.
\end{align}
From (\ref{ecu16}) the Lemma following.\\
{\bf Remark.}
If the assumption  \eqref{ecudoc3} holds, then
$$e^{\frac{ix^2}{4t}}u\in C(\mathbb{R}-\{0\}:H^m(\mathbb{R})),\ \ e^{\frac{ix^2}{4t}}v\in C(\mathbb{R}-\{0\}:H^m(\mathbb{R}))$$
$\textbf{Proof of the main theorem.}$ From Lemma (\ref{lemma501}) and (\ref{lemma502}) we obtain that there exists $u=u(x,t)$ and $v=v(x,t)$
satisfying \eqref{ecudoc1}-\eqref{ecudoc2} and such that for any $T>0$ we have
\begin{align*}
&u^k\ \longrightarrow u\ strongly\ in\ C(\mathbb{R}:H^1(\mathbb{R}))\\
&v^k\ \longrightarrow v\ strongly\ in\ C(\mathbb{R}:H^1(\mathbb{R}))
\end{align*}
and
\begin{align*}
&J^mu^k\ \longrightarrow J^mu\ strongly\ in\ C(\mathbb{R}:L^2(\mathbb{R}))\\
&J^mv^k\ \longrightarrow J^mv\ strongly\ in\ C(\mathbb{R}:L^2(\mathbb{R})).
\end{align*}
It is easily verified that $(u, v)$ solves (\ref{main2}). The proof
then follows.
\begin{corollary}
If the hypotheses in Lemma \ref{lemma502} are satisfied, then
\begin{align*}
u\in\bigcup_{m=0}^{[\frac{n}{2}]}C^m(\mathbb{R}-\{0\}:C^{n-2m-1}(\mathbb{R}))\ \ and\ \ v\in\bigcup_{m=0}^{[\frac{n}{2}]}C^m(\mathbb{R}-\{0\}:C^{n-2m-1}(\mathbb{R})).
\end{align*}
\end{corollary}
\begin{corollary}
If $(x^nu_0,x^nv_0)\in L^2(\mathbb{R})\times L^2(\mathbb{R})$ for all $n\in\mathbb{N}$ in Lemma \ref{lemma502}
then the solution $(u,\ v)$ of \eqref{main2} is infinitely differentiable in $x$ and $t$ for $t\neq0$.
\end{corollary}

\section{Numerical Experiments}

In this section, we will present some numerical simulations which replicate the decay rates proved recently. We will find approximate solutions of \eqref{main2} using a finite differences approach. Due to computational limitations, we will assume from now on that $\Omega$ is bounded. The derivatives will be approximated using second-order centered finite differences, and the nonlinear term will be treated as proposed in Delfour \cite{delfour}. \\

We will briefly describe the numerical scheme used to obtain our results. Due to computational limitations, we will consider a bounded domain sufficiently large so that the boundary doest not interfere with our results. Let us denote it by $\Omega := [x_0,x_f]\subset \mathbb{R}$. For a given $M\in\mathbb{N}$, we will introduce the vector space
  \begin{equation*}
    X_M := \big\{y = [y_0\:  y_1 \: \dots \: y_M]^T \in \mathbb{C}^{M+1} : y_0 = y_{M-1} = y_{M} = 0  \big\}
  \end{equation*}
  This condition mimicks the boundary conditions $y = y(x,t) = 0$ for $x\in \{x_0,x_f\}$ and $y_x = 0$ for $x = x_f$. For $\Delta x := \frac{x_f-x_0}{M-1}$, we introduce the following classical finite differences operators for complex-valued arrays:  
  \begin{align*}
    \big[\boldsymbol{D}^+ u\big]_j = D^+y_j &:= \dfrac{y_{j+1} - y_j}{\Delta x} \\
    \big[\boldsymbol{D}^- y\big]_j = D^-y_j &:= \dfrac{y_{j} - y_{j-1}}{\Delta x} \\
    \boldsymbol{D}_x y  &:= \dfrac{1}{2}\Big(\boldsymbol{D}^+y + \boldsymbol{D}^- y \Big) \\
    \boldsymbol{D}_x^2 y &:= \boldsymbol{D}^+\boldsymbol{D}^- y
  \end{align*}
  Our numerical scheme is defined as follows: for $u^n,\: v^n \in X_M$ and $t_n = n\Delta t,\:\Delta t<1$:
  \begin{align*}
    i\dfrac{u^{n+1}-u^n}{\Delta t}  + \boldsymbol{D}_x^2u^{n+\frac{1}{2}} &= \dfrac{|u^{n+1}|^{2p+2} - |u^n|^{2p+2}}{|u^{n+1}|^2 - |u^n|^2}\dfrac{u^{n+1}+u^n}{2p+2} \\
    &\quad + \beta \dfrac{|v^{n+1}|^{p+1} + |v^n|^{p+1}}{2}\dfrac{|u^{n+1}|^{p+1} - |u^n|^{p+1}}{|u^{n+1}|^2 - |u^n|^2}\dfrac{u^{n+1}+u^n}{p+1} \\
        i\dfrac{v^{n+1}-v^n}{\Delta t}  + \boldsymbol{D}_x^2v^{n+\frac{1}{2}} &= \dfrac{|v^{n+1}|^{2p+2} - |v^n|^{2p+2}}{|v^{n+1}|^2 - |v^n|^2}\dfrac{v^{n+1}+v^n}{2p+2} \\
    &\quad + \beta \dfrac{|u^{n+1}|^{p+1} + |u^n|^{p+1}}{2}\dfrac{|v^{n+1}|^{p+1} - |v^n|^{p+1}}{|v^{n+1}|^2 - |v^n|^2}\dfrac{v^{n+1}+v^n}{p+1}
  \end{align*}
  This is a nonlinear problem; thus a Picard fixed point iteration is used in each timestep. This leads to the solving of the same linear system of equations many times per step. As the coefficients matrix is tridiagonal, we use the {\tt LAPACK} package to solve the linear system via LU factorization. In order to avoid overflow/underflow errors, each code is programmed to include the exact factorization associated to each value of $p$. For example, when $p = 1$ we get
    \begin{align*}
      i\dfrac{u^{n+1}-u^n}{\Delta t}  + \boldsymbol{D}_x^2u^{n+\frac{1}{2}} &= \dfrac{|u^{n+1}|^{2} + |u^n|^{2}}{2}\dfrac{u^{n+1}+u^n}{2} + \beta \dfrac{|v^{n+1}|^{2} + |v^n|^{2}}{2}\dfrac{u^{n+1}+u^n}{2} \\
      i\dfrac{v^{n+1}-v^n}{\Delta t}  + \boldsymbol{D}_x^2v^{n+\frac{1}{2}} &= \dfrac{|v^{n+1}|^{2} + |v^n|^{2}}{2}\dfrac{v^{n+1}+v^n}{2} + \beta \dfrac{|u^{n+1}|^{2} + |u^n|^{2}}{2}\dfrac{v^{n+1}+v^n}{2},   \end{align*}
    which is a classical discretization for the nonlinear term, as proposed in Delfour, Fortin and Payre \cite{delfour}. Regarding the norms to be studied, let us define the discrete $L^p$ norm as follows
    \[ ||u||_p^p := \Delta x \sum\limits_{i=1}^M|u_i|^p,\quad  \text{ for }u\in \mathbb{C}^M,\: p\in \mathbb{N}. \]
    The discrete $L^\infty$ norm will be defined as
    \[||u||_\infty = \max\limits_{1\leq i \leq M}|u_i| \]
    Thus, the numerical Energy  will be given by
    \[E(u,v) := ||\boldsymbol{D}^+_xu||^2_2 + ||\boldsymbol{D}^+_x v||_2^2 + \dfrac{1}{p+1}||u||_{2p+2}^{2p+2} + \dfrac{1}{p+1}||v||_{2p+2}^{2p+2} + \dfrac{2\beta}{p+1}\Delta x \sum_{i=1}^N|u_i|^{p+1}|v_i|^{p+1}. \]
    
    \subsection{Example 1: $p = 3$.} Starting with $p = 3$, and considering $u^{n+\frac{1}{2}} = \dfrac{u^{n+1}+u^n}{2}$, the scheme will be
        \begin{align*}
          i\dfrac{u^{n+1}-u^n}{\Delta t}  + \boldsymbol{D}_x^2u^{n+\frac{1}{2}} &= \dfrac{|u^{n+1}|^{2} + |u^n|^{2}}{2}\dfrac{|u^{n+1}|^{4} + |u^n|^{4}}{2}u^{n+\frac{1}{2}} \\
          &\quad + \beta \dfrac{|v^{n+1}|^{4} + |v^n|^{4}}{2}\dfrac{|u^{n+1}|^{2} + |u^n|^{2}}{2}u^{n+\frac{1}{2}} \\
      i\dfrac{v^{n+1}-v^n}{\Delta t}  + \boldsymbol{D}_x^2v^{n+\frac{1}{2}} &= \dfrac{|v^{n+1}|^{2} + |v^n|^{2}}{2}\dfrac{|v^{n+1}|^{4} + |v^n|^{4}}{2}v^{n+\frac{1}{2}} \\
          &\quad + \beta \dfrac{|u^{n+1}|^{4} + |u^n|^{4}}{2}\dfrac{|v^{n+1}|^{2} + |v^n|^{2}}{2}v^{n+\frac{1}{2}}  \end{align*}

\begin{center}
        \begin{figure}
          \centering
          \includegraphics[scale=0.4]{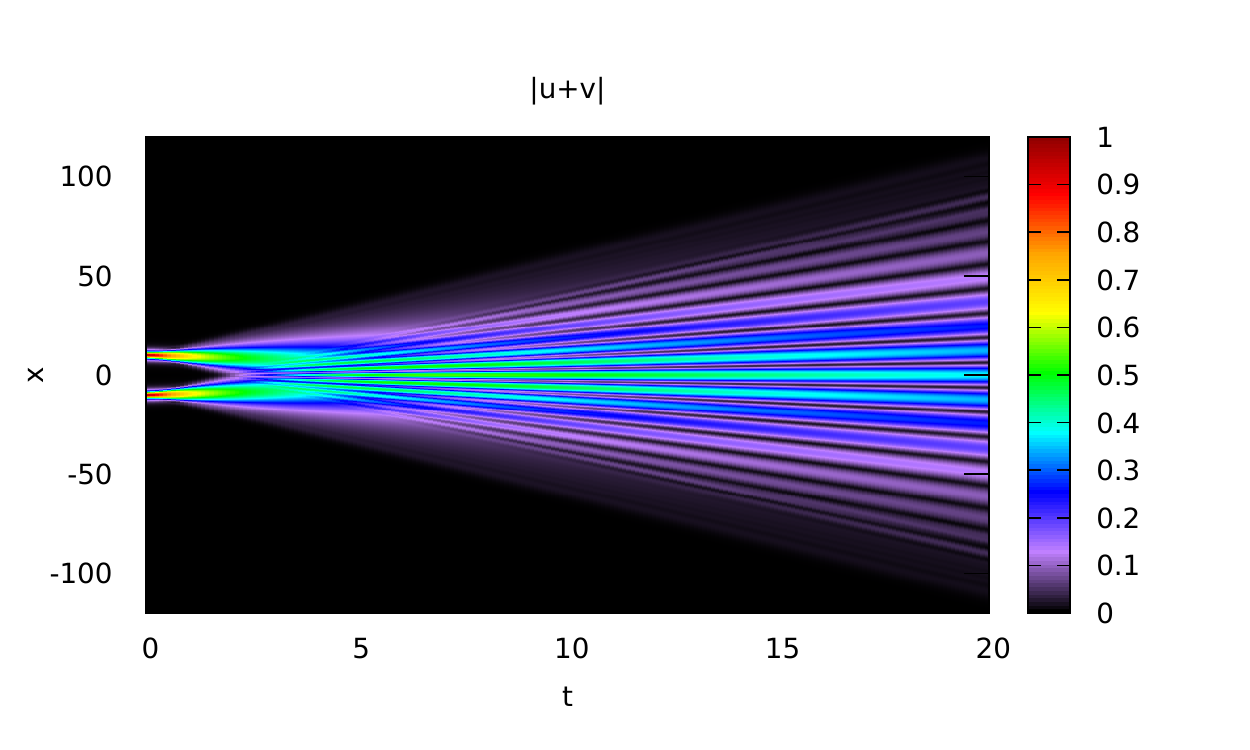}
          \includegraphics[scale=0.4]{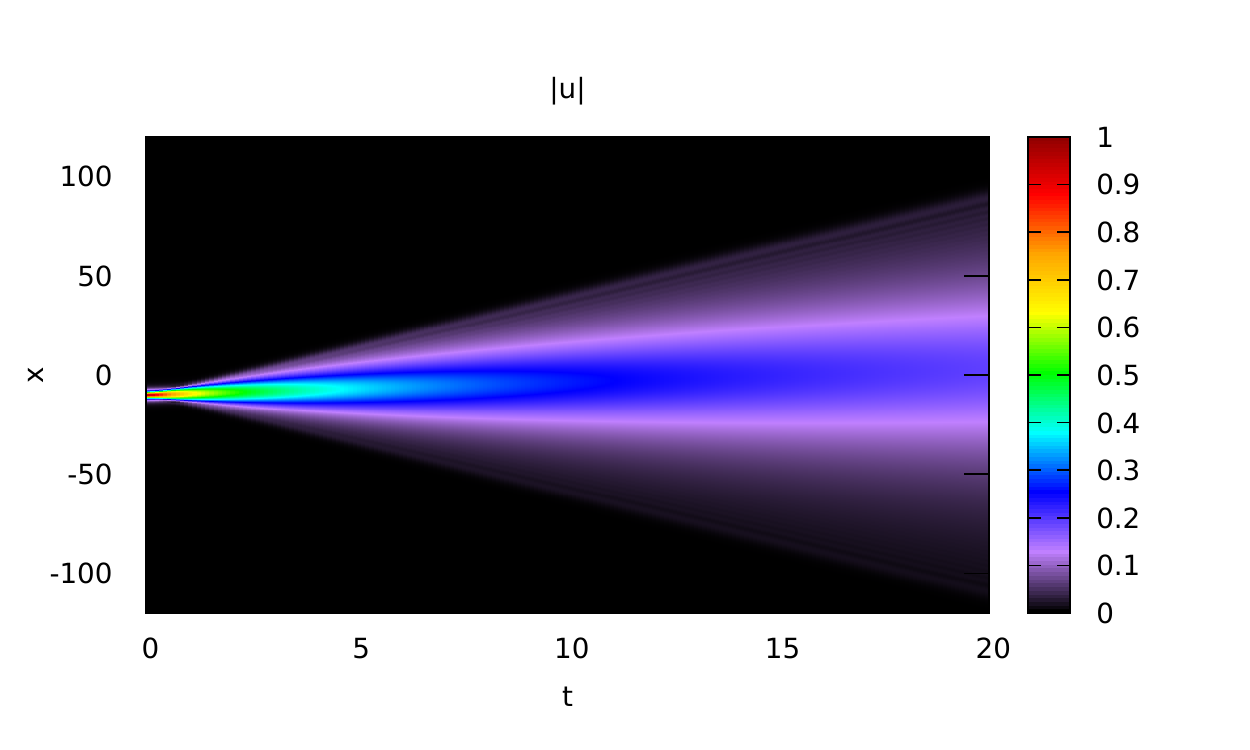}
          \includegraphics[scale=0.4]{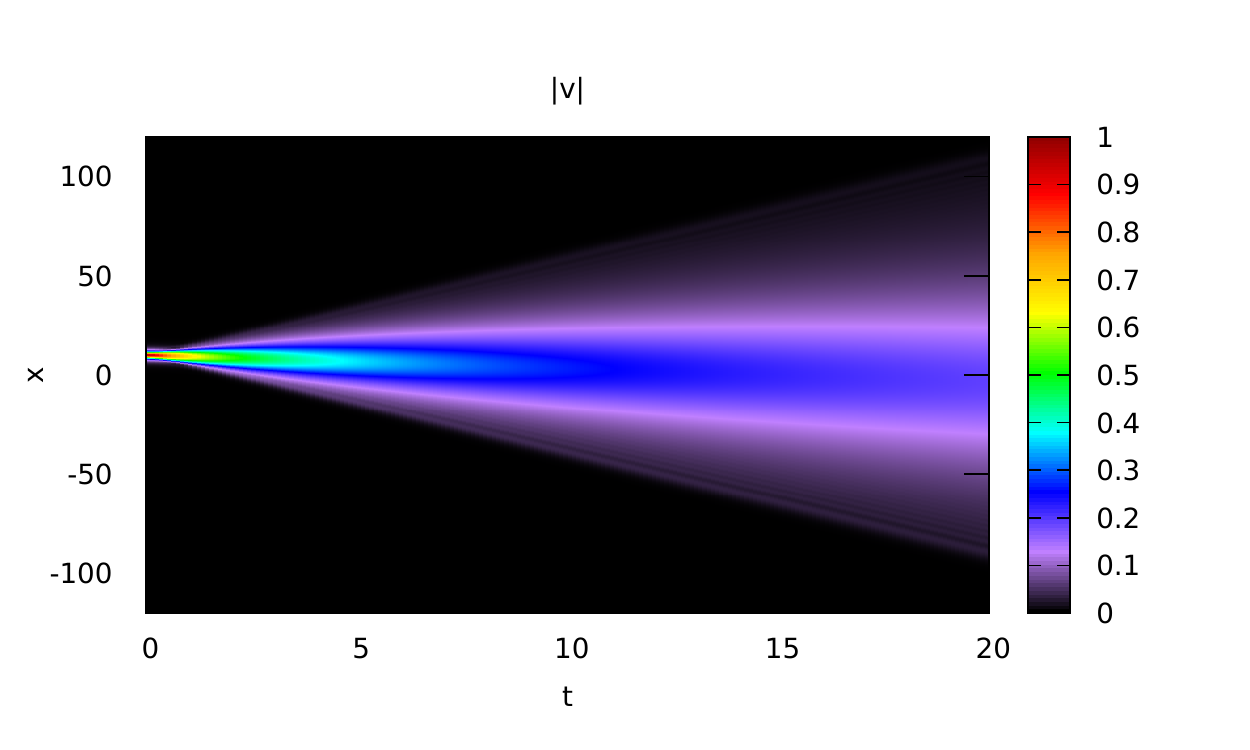}
          \caption{Modulus of the numerical solution, Example 1.}
          \label{fig_abs_1}
        \end{figure}
      \end{center}

\begin{center}
        \begin{figure}
          \centering
          \includegraphics[scale=0.62]{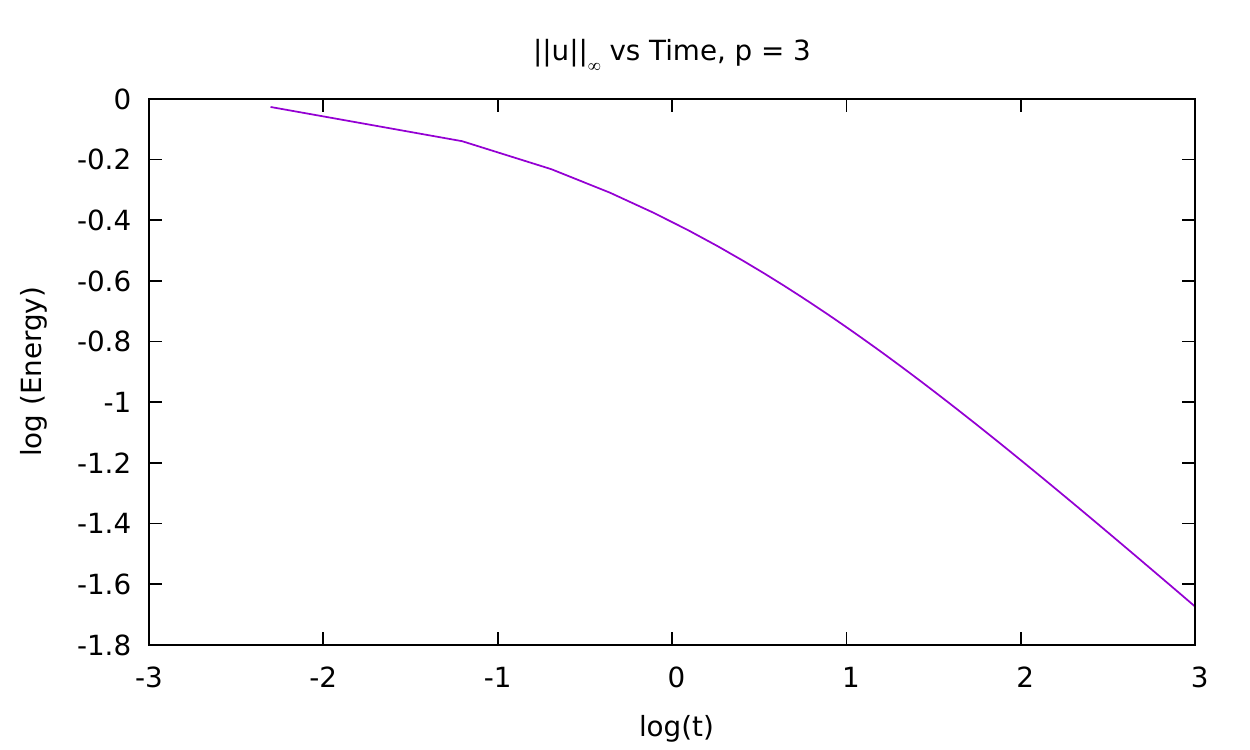}
          \includegraphics[scale=0.62]{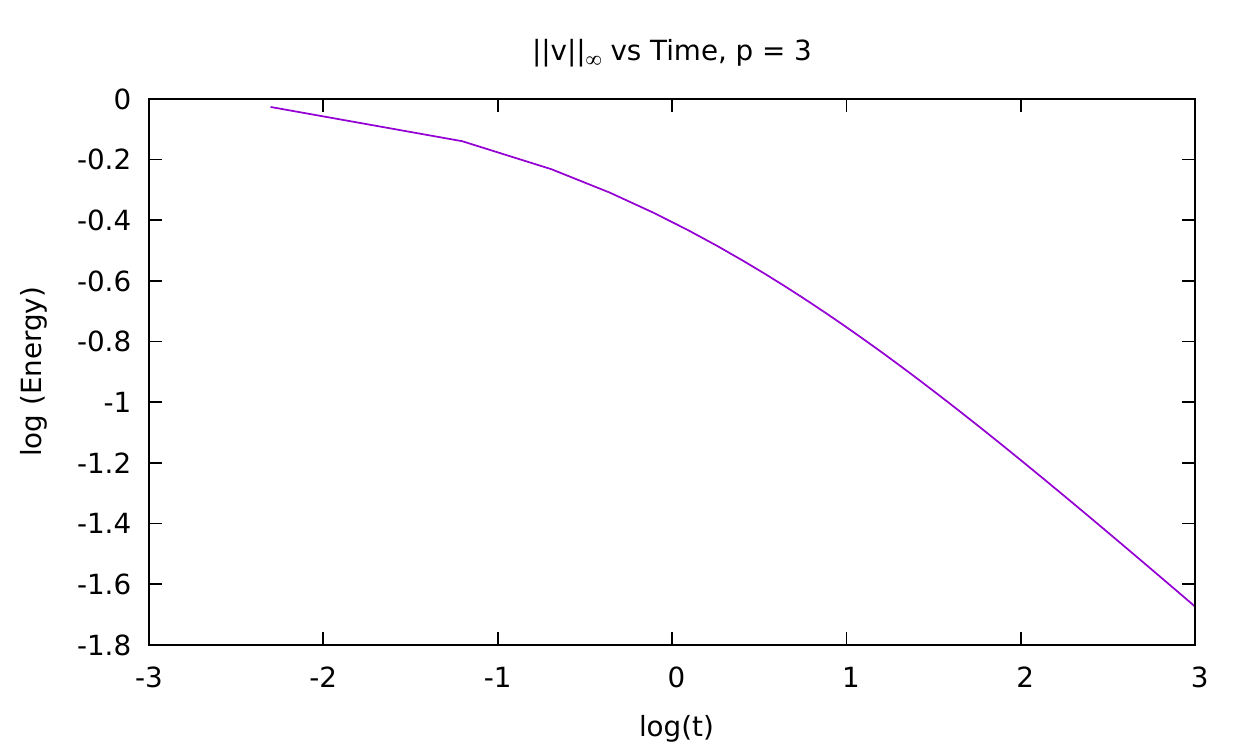}
          \caption{$L^\infty$ norm evolution, Example 1.}
          \label{fig_ej1_inf}
          \includegraphics[scale=0.62]{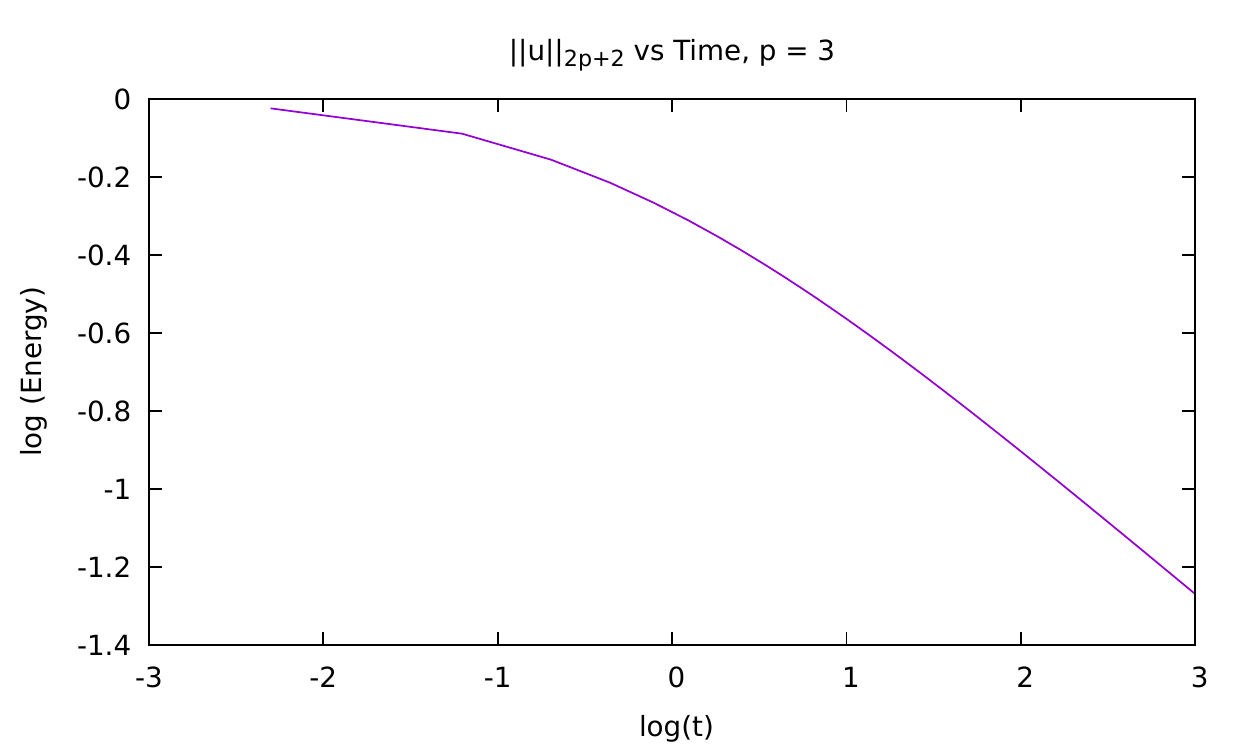}
          \includegraphics[scale=0.62]{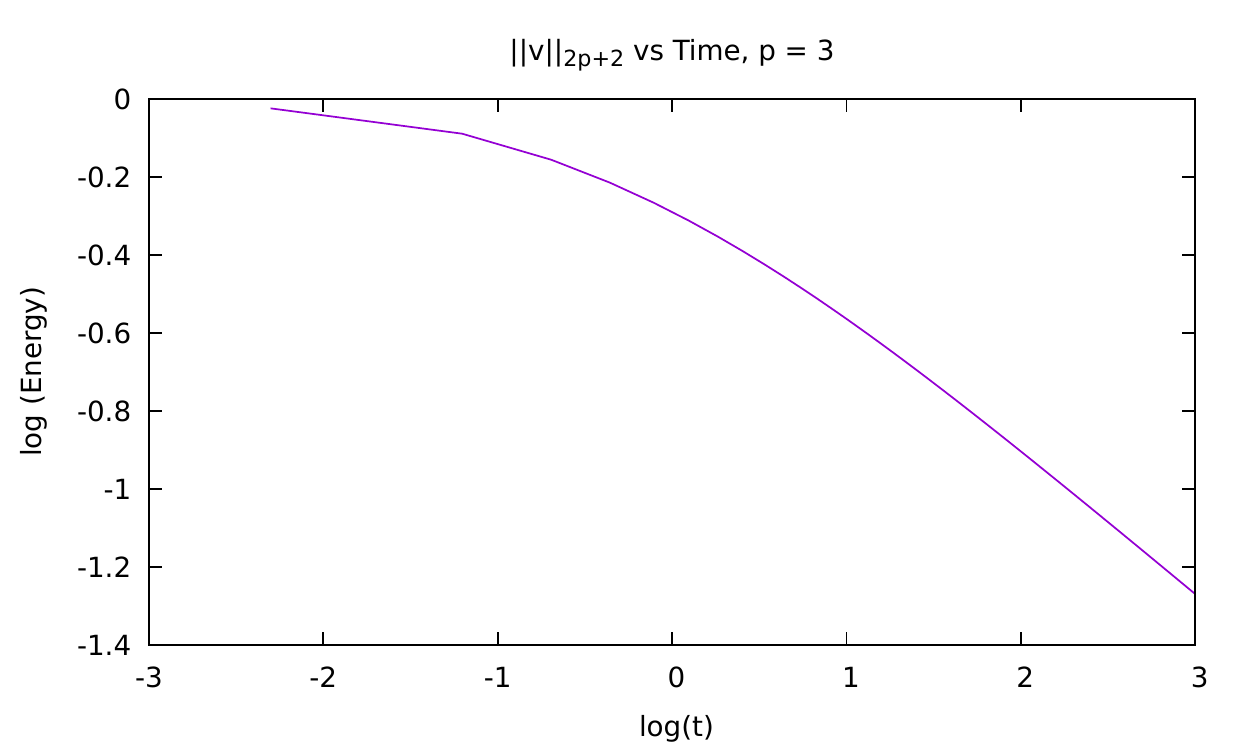}
          \caption{$L^{2p+2}$ norm evolution, Example 1.}
          \label{fig_ej1_pmas2}
          \includegraphics[scale=0.62]{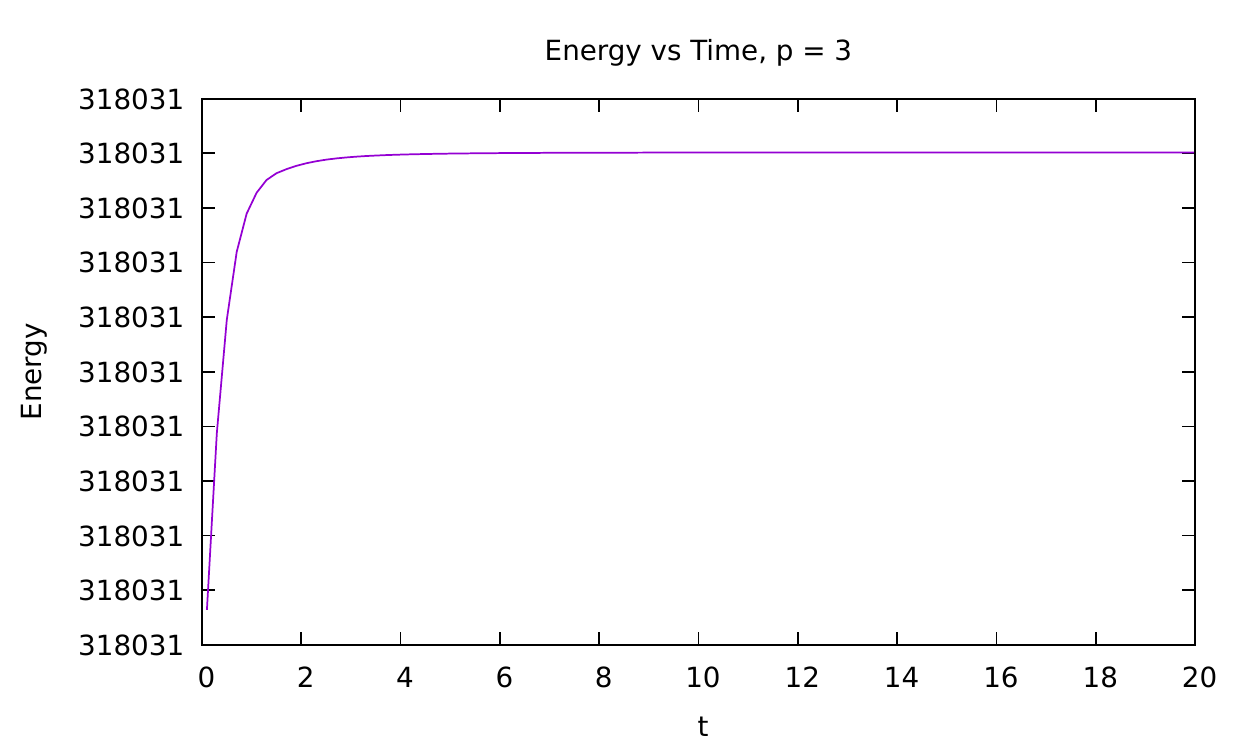}
          \includegraphics[scale=0.62]{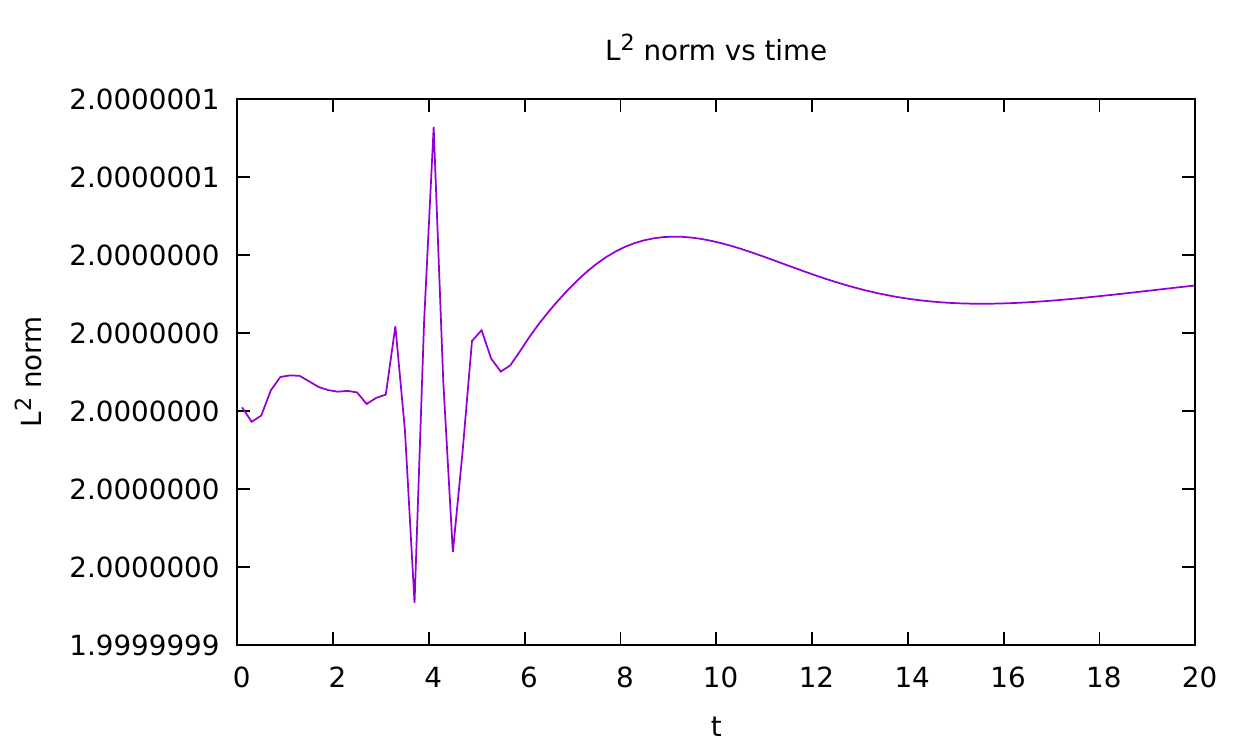}
          \caption{Energy and $L^2$-norm time evolution, Example 1.}
          \label{energia_ej1}
        \end{figure}
      \end{center}

    The following result is obtained from the initial condition
    \[u(x,0) = 1.2\sqrt{2}e^{1.3i\frac{x}{4}}\operatorname{sech}\Big(1.2x + 10\Big),\qquad v(x,0) = \sqrt{2}e^{-1.3i\frac{x}{4}}\operatorname{sech}\Big(x - 10\Big), \]
   which is a classical solution of the NLS equation \cite{polyanin} for $p = 2$, and in practice it corresponds to a soliton collision due to the differences in their velocities. Our computations consider the domain $\Omega = [-100,100]$, $\Delta t = 0.1$, $\Delta x = \dfrac{200}{2^{13}} \approx 0.024414$, ans $\beta = 1$. Figure \ref{fig_ej1_inf} and \ref{fig_ej1_pmas2} illustrates the time evolution of the $L^\infty$ and $L^{2p+2}$ norms respectively. Figure \ref{energia_ej1} shows the evolution of the preserved quantities. The energy and the $L^2$ norm remain numerically preserved, while the $L^\infty$ and $L^{2p+2}$ norms decay following a power law but with a steeper slope than the one predicted by theory for the $L^\infty$ norm. This might be due to numerical dispersion effects that contribute to its decay (for a detailed study, see \cite{fibich2003}). Figure \ref{fig_abs_1} illustrates the numerical solutions obtained. 

    \subsection{Example 2: $p = 5$.} 
\begin{center}
        \begin{figure}
          \centering
          \includegraphics[scale=0.4]{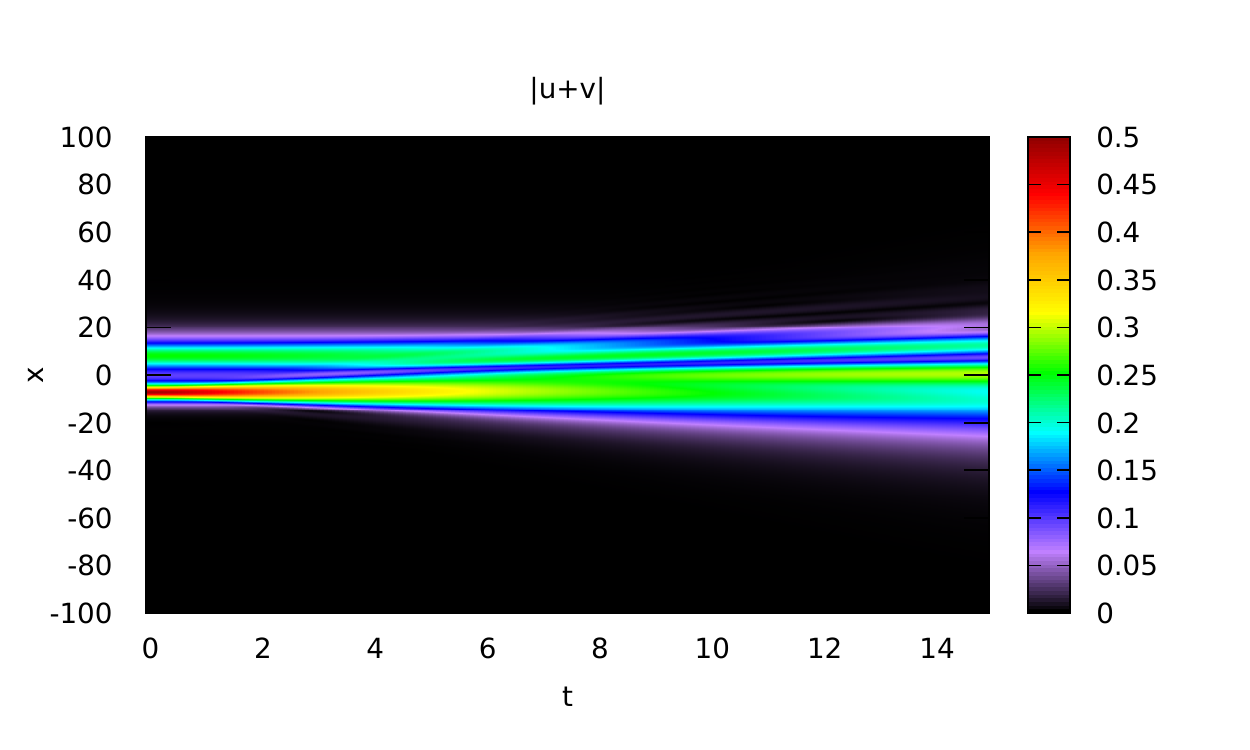}
          \includegraphics[scale=0.4]{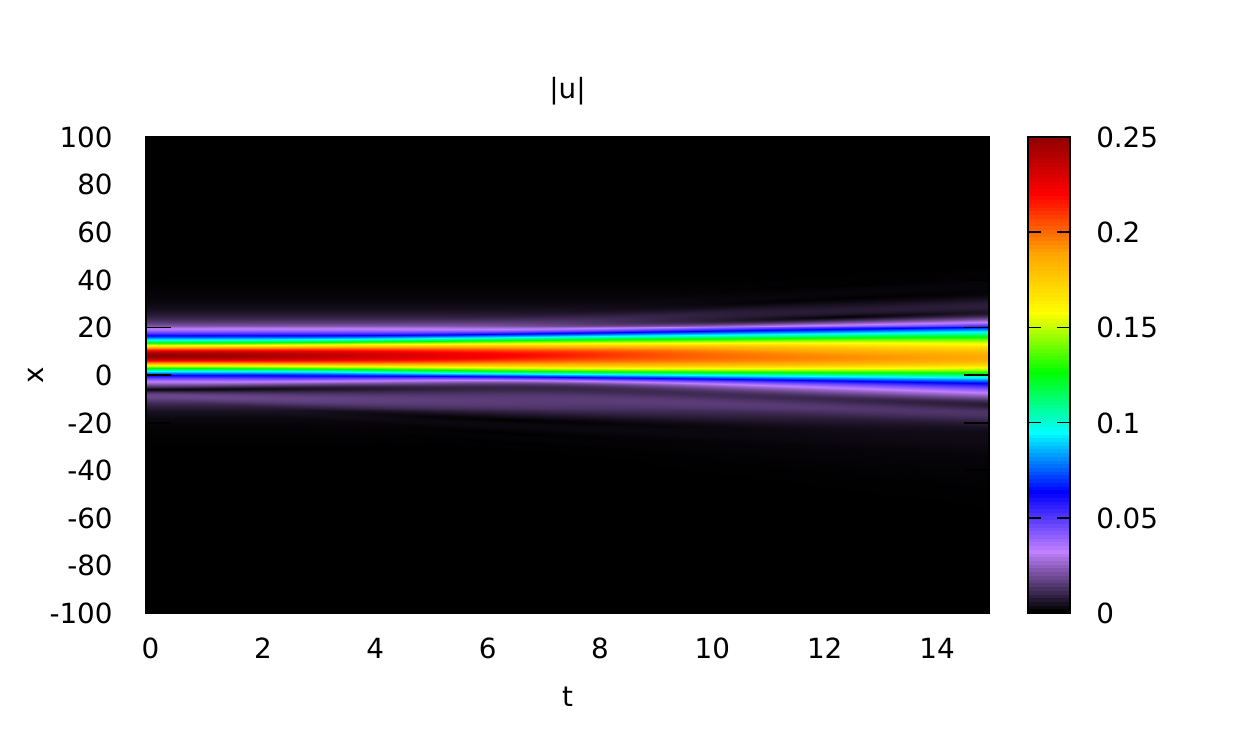}
          \includegraphics[scale=0.4]{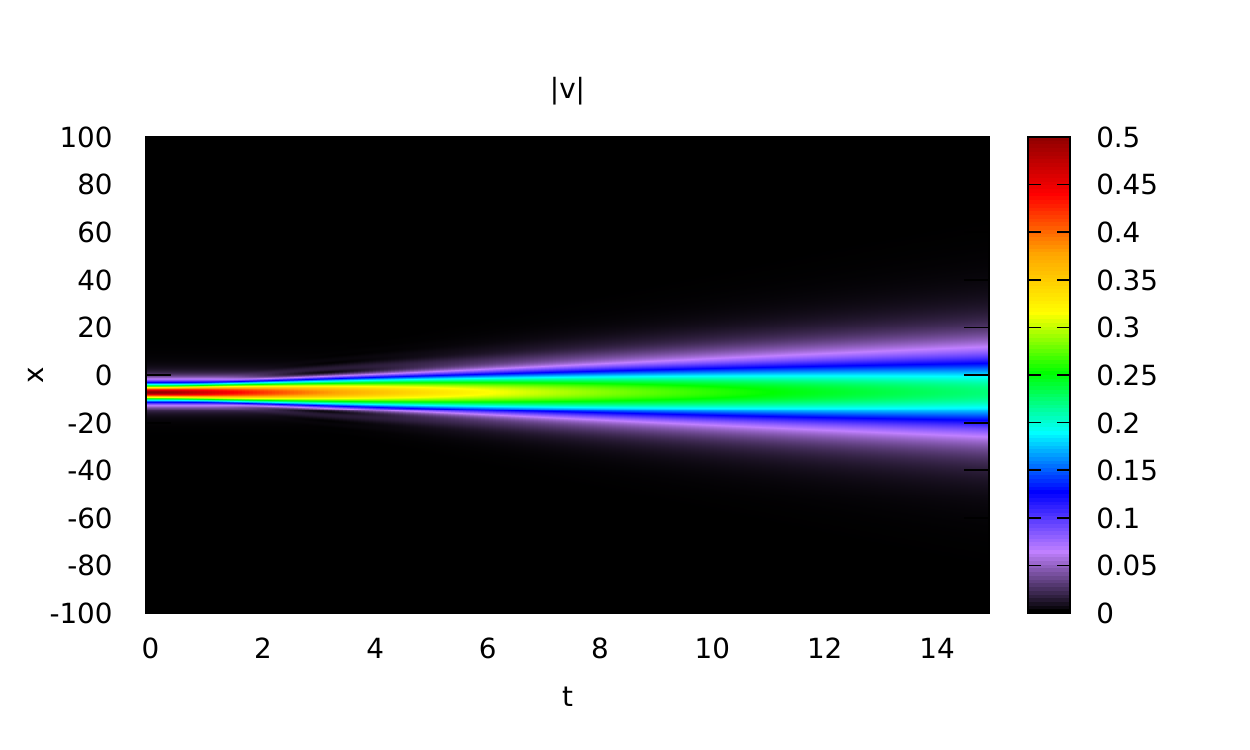}
          \caption{Modulus of the numerical solution, Example 2.}
          \label{fig_abs_2}
        \end{figure}
      \end{center}

\begin{center}
        \begin{figure}
          \centering
          \includegraphics[scale=0.62]{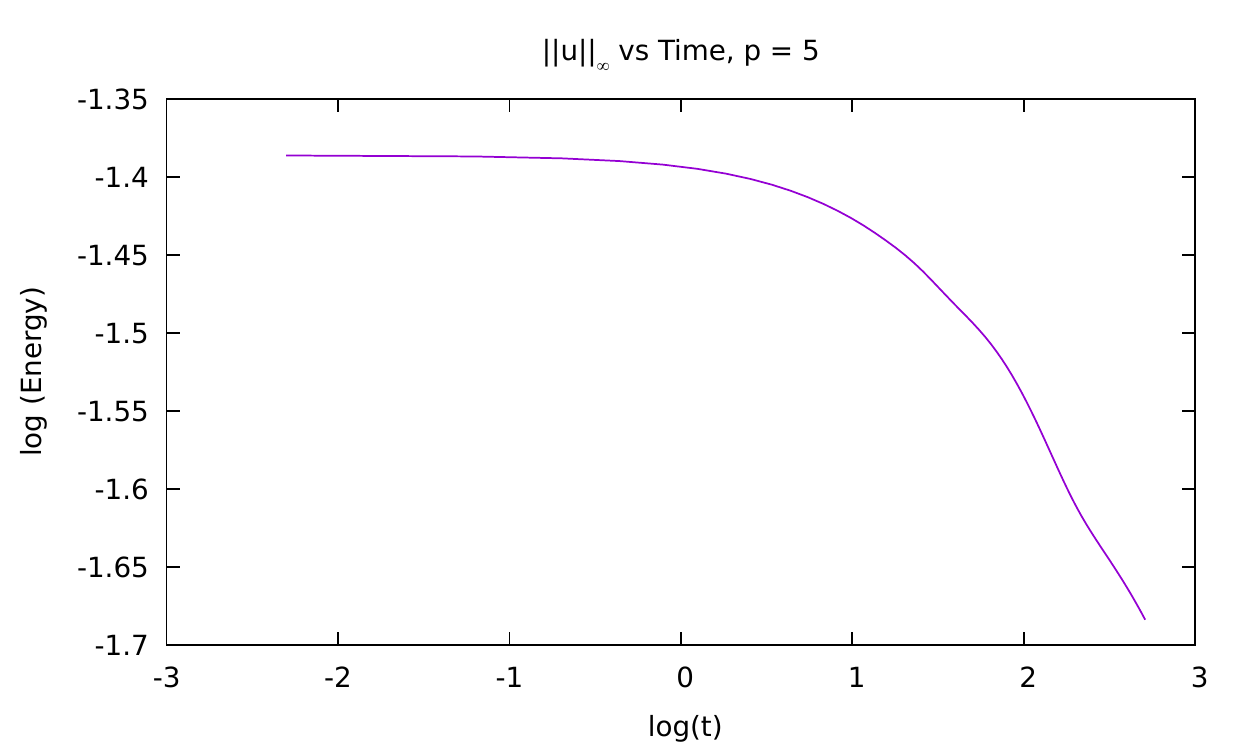}
          \includegraphics[scale=0.62]{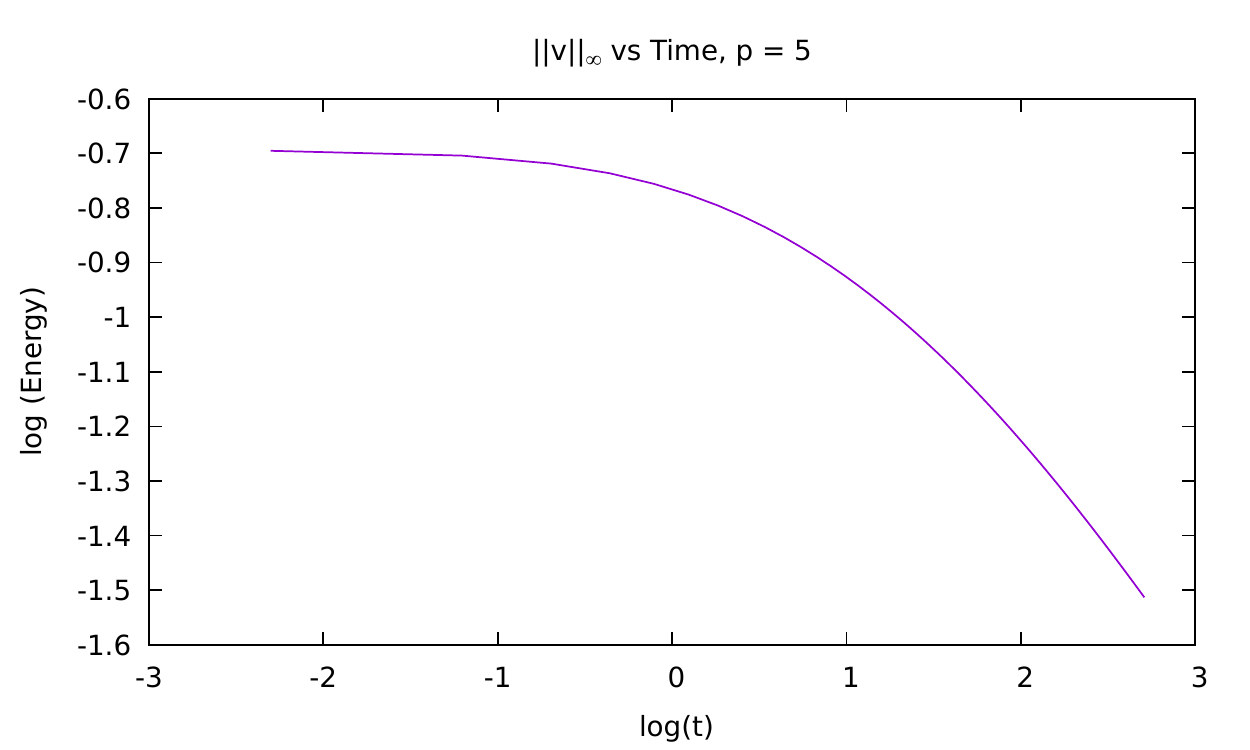}
          \caption{$L^\infty$ norm evolution, Example 2.}
          \label{fig_ej2_inf}
          \includegraphics[scale=0.62]{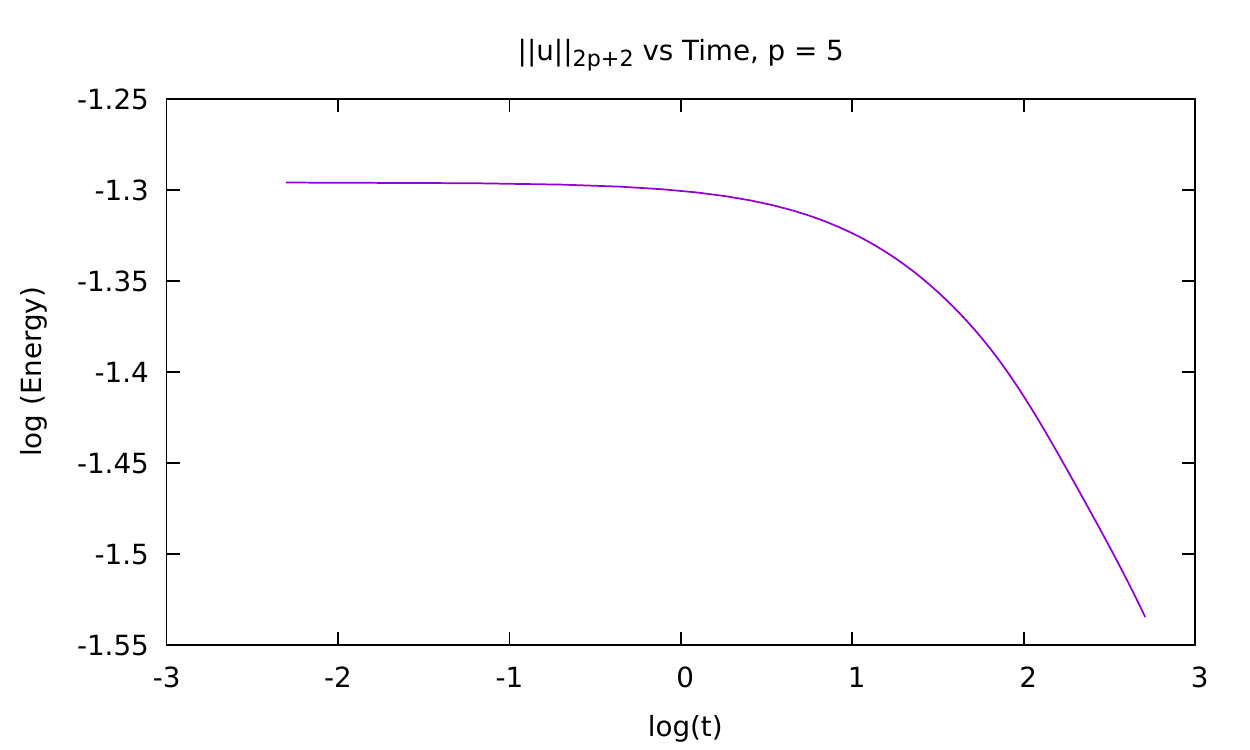}
          \includegraphics[scale=0.62]{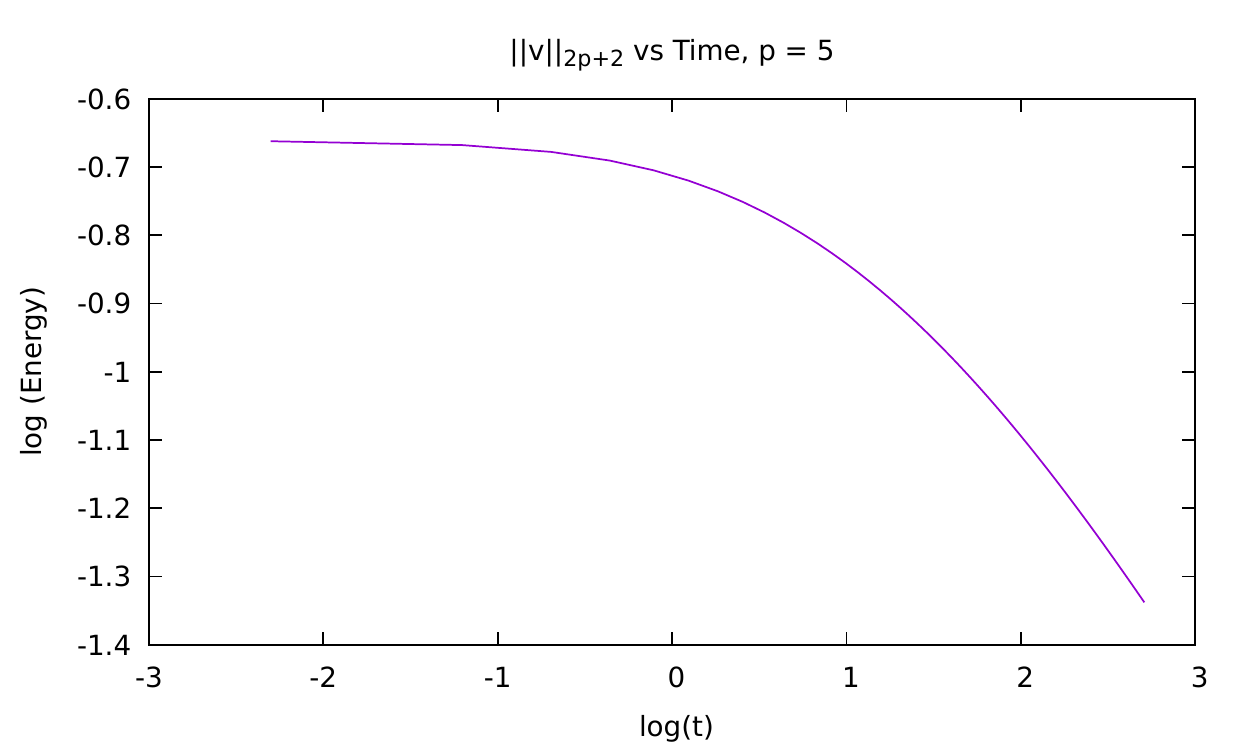}
          \caption{$L^{2p+2}$ norm evolution, Example 2.}
          \label{fig_ej2_pmas2}
          \includegraphics[scale=0.62]{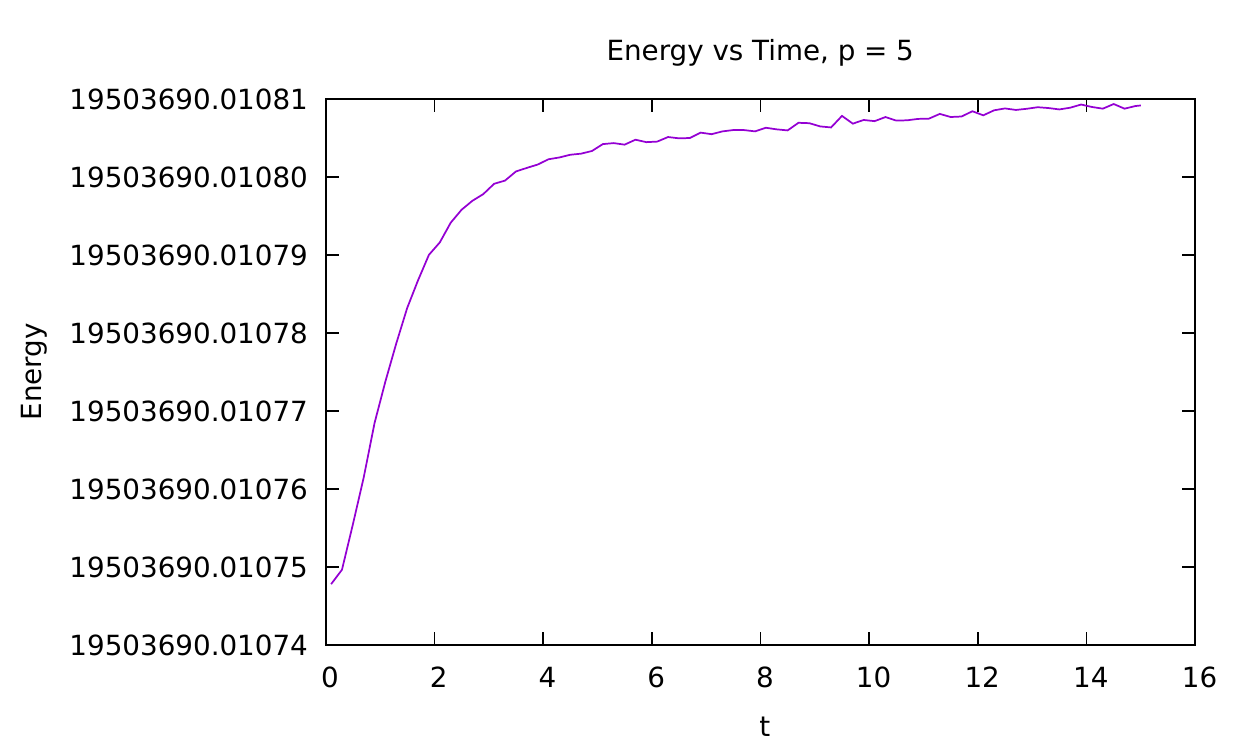}
          \includegraphics[scale=0.62]{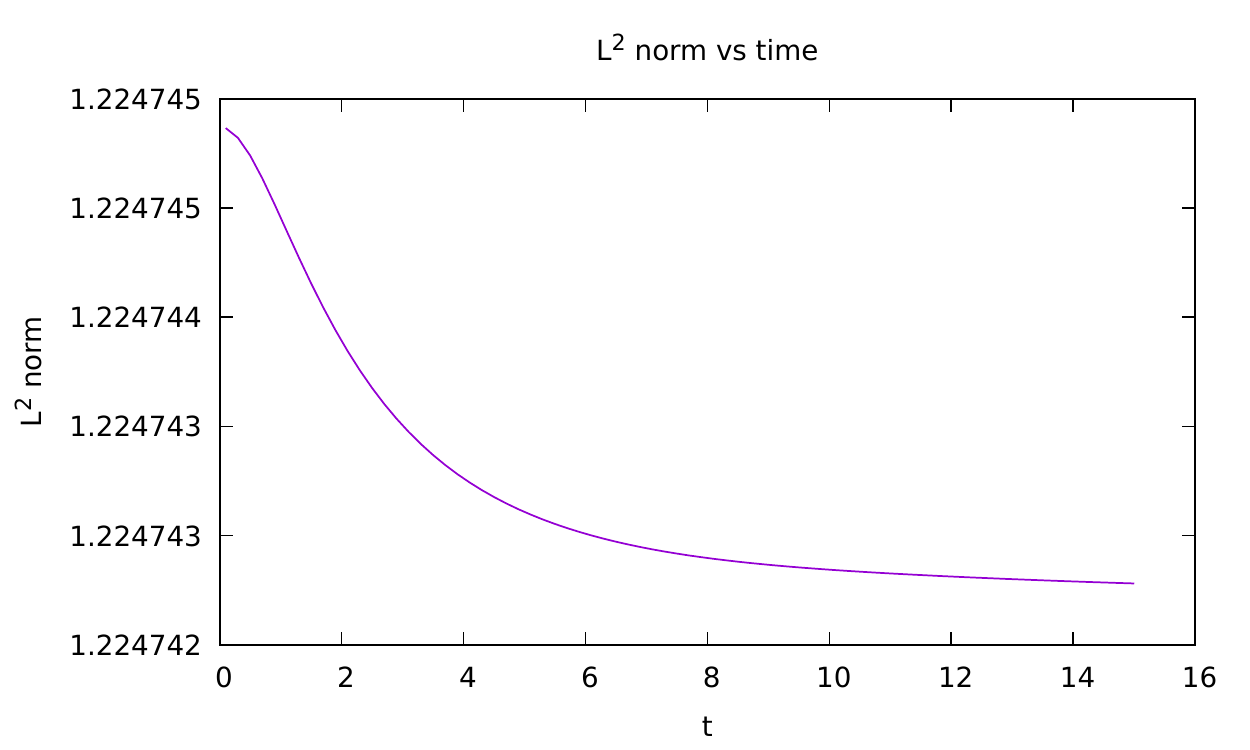}
          \caption{Energy and $L^2$-norm time evolution, Example 2.}
          \label{energia_ej2}
        \end{figure}
      \end{center}

With that value for $p$, the scheme becomes
      \begin{align*}
    i\dfrac{u^{n+1}-u^n}{\Delta t}  + \boldsymbol{D}_x^2u^{n+\frac{1}{2}} &= \dfrac{\big(|u^{n+1}|^{4} + |u^{n+1}|^2|u^n|^2+ |u^n|^{4}\big)\big(|u^{n+1}|^{6} + |u^n|^{6}\big)}{6}\dfrac{u^{n+1}+u^n}{2} \\
    &\quad + \beta \dfrac{|v^{n+1}|^{6} + |v^n|^{6}}{2}\dfrac{|u^{n+1}|^{4} +|u^{n+1}|^2|u^n|^2+ |u^n|^{4}}{3}\dfrac{u^{n+1}+u^n}{2} \\
        i\dfrac{v^{n+1}-v^n}{\Delta t}  + \boldsymbol{D}_x^2v^{n+\frac{1}{2}} &= \dfrac{\big(|v^{n+1}|^{4} + |v^{n+1}|^2|v^n|^2+ |v^n|^{4}\big)\big(|v^{n+1}|^{6} + |v^n|^{6}\big)}{6}\dfrac{v^{n+1}+v^n}{2} \\
    &\quad + \beta \dfrac{|u^{n+1}|^{6} + |u^n|^{6}}{2}\dfrac{|v^{n+1}|^{4} +|v^{n+1}|^2|v^n|^2+ |v^n|^{4}}{3}\dfrac{v^{n+1}+v^n}{2}
  \end{align*}
      For this example we will consider as initial condition the expressions (60a) and (60b) from Menyuk \cite{menyuk2}, and corresponds to another soliton collision. As parameters, we have chosen $A_1 = \frac{1}{4}$, $A_2 = \frac{1}{2}$, $s_1 = 8$, $s_2 = -5$, and $\delta = 0$. The domain, timestep and mesh size are the same from the previous example. Figure \ref{fig_ej2_inf} illustrates the evolution of the $L^\infty$-norm, while Figure \ref{fig_ej2_pmas2} does for $L^{2p+2}$-norm. They decay acording to the a power law, as expected. Figure \ref{energia_ej2} shows the time evolution of the preserved quantities. We can see that, at a numerical level, they are preserved as well. Figure \ref{fig_abs_2} shows the numerical solutions obtained.

\section*{Acknowledgments}

\noindent R. Nina-Mollisaca has been supported by postgraduate scholarship {\it Apoyo Institucional para Doctorado Acad\'emicos de la Universidad de Tarapac\'a}, Universidad de Tarapac\'a, Chile. O. Vera-Villagran has been supported by Project Fondecyt 1191137 (Universidad de Tarapac\'a, Chile). \\

\noindent M. Sep\'ulveda was supported FONDECYT Grant No. 1220869,
ECOS-ANID project C20E03, and by Centro de Modelamiento Matem\'atico (CMM), ACE210010 and FB210005, BASAL funds for centers of excellence from ANID-Chile. \\

\noindent  Authors must disclose all relationships or interests that
 could have direct or potential influence or impart bias on
 the work:

 \section*{Conflict of interest}

 The authors declare that they have no conflict of interest.

\end{document}